\documentclass[10pt, reqno]{amsart}
\usepackage{amssymb,latexsym,calligra}
\usepackage{amsthm}
\usepackage{amsmath,mathrsfs, amsfonts}

\usepackage{color}
\usepackage[OT2, T1]{fontenc}
\usepackage{url}
\usepackage{array}
\usepackage{hyperref}
\usepackage{paralist}
\usepackage{mathabx}			
\usepackage{amscd}
\usepackage[all,cmtip]{xy}			
\usepackage{sseq}	
\usepackage{verbatim}
\usepackage{parskip}
\usepackage{microtype}
\usepackage[in]{fullpage}
\usepackage{mathrsfs} 			
\usepackage{cleveref}
\usepackage{enumitem}
\usepackage{moreenum}			
\usepackage[alphabetic,lite]{amsrefs} 	
\usepackage{colonequals}	
\usepackage{fontenc}
\usepackage{amsmath}
\usepackage{geometry}
\usepackage{graphicx}
\usepackage{stmaryrd}
\usepackage[all]{xy}
\usepackage[french]{babel}
\usepackage[latin1]{inputenc}
\usepackage[T1]{fontenc}
\usepackage[foot]{amsaddr}
\usepackage{subfiles}
\usepackage{ragged2e}
\usepackage[bbgreekl]{mathbbol}     
\DeclareSymbolFontAlphabet{\mathbb}{AMSb}   

\DeclareMathOperator{\Mod}{Mod}

\DeclareMathOperator{\Sym}{Sym}

\DeclareMathOperator{\dAff}{dAff}
\DeclareMathOperator{\dPrStk}{dPrStk}

\DeclareMathOperator{\bD}{\mathbb{D}}

\DeclareMathOperator{\cL}{\mathcal{L}}

\DeclareMathOperator{\Spec}{Spec}

\DeclareMathOperator{\cHom}{\underline{Hom}}
\DeclareMathOperator{\supp}{supp}
\DeclareMathOperator{\Id}{Id}
\DeclareMathOperator{\Hom}{Hom}

\DeclareMathOperator{\Fact}{Fact}
\DeclareMathOperator{\Cat}{Cat}
\DeclareMathOperator{\End}{End}

\DeclareMathOperator{\Ima}{Im}
\DeclareMathOperator{\Bun}{Bun}

\DeclareMathOperator{\Gr}{Gr}
\DeclareMathOperator{\Spc}{Spc}
\DeclareMathOperator{\Aff}{Aff}
\DeclareMathOperator{\Sh}{Sh}

\DeclareMathOperator{\Rhom}{R\mathscr{H}\text{\kern -3pt {\calligra\large om}}\,}
\DeclareMathOperator{\EExt}{\mathscr{E}\text{\kern -3pt {\calligra\large xt}}\,}
\DeclareMathOperator{\Ind}{Ind}

\DeclareMathOperator{\ad}{ad}

\DeclareMathOperator{\Ens}{Ens}

\DeclareMathOperator{\Pro}{Pro}
\DeclareMathOperator{\TopMod}{TopMod}

\DeclareMathOperator{\sinf}{<}

\DeclareMathOperator{\ssup}{>}

\DeclareMathOperator{\Loc}{Loc}

\DeclareMathOperator{\coker}{Coker}
\DeclareMathOperator{\Ker}{Ker}

\DeclareMathOperator{\Cofib}{Cofib}
\DeclareMathOperator{\Map}{Map}

\DeclareMathOperator{\car}{car}
\DeclareMathOperator{\Sch}{Sch}

\DeclareMathOperator{\QCoh}{QCoh}
\DeclareMathOperator{\colim}{\varinjlim}

\DeclareMathOperator{\Ext}{Ext}
\DeclareMathOperator{\Res}{Res}
\DeclareMathOperator{\Tor}{Tor}

\DeclareMathOperator{\X}{X_{\text{proét}}}

\DeclareMathOperator{\St}{St}
\DeclareMathOperator{\PrSh}{PrSh}
\DeclareMathOperator{\PreSh}{PrSh}
\DeclareMathOperator{\AlgSp}{AlgSp}

\DeclareMathOperator{\bzl}{\overline{\mathbb{Z}}_{\ell}}

\DeclareMathOperator{\cT}{\mathcal{T}}

\DeclareMathOperator{\Lie}{Lie}

\DeclareMathOperator{\Fonct}{Fonct}
\DeclareMathOperator{\Frac}{Frac}

\DeclareMathOperator{\PrCat}{PrCat_{\ell}}
\DeclareMathOperator{\Catst}{Cat_{st,\ell}}

\DeclareMathOperator{\Der}{Der}

\DeclareMathOperator{\clp}{\mathcal{L}^{+} }

\DeclareMathOperator{\Exalcom}{Exalcom}

\usepackage{aliascnt}
\newaliascnt{numberingbase}{subsubsection}
\numberwithin{equation}{numberingbase}

\newtheoremstyle{thms}{0.pt}{0pt}{\itshape}{}{\bfseries}{.}{ }{}
\theoremstyle{thms}
\newtheorem{conj}[numberingbase]{Conjecture}
\newtheorem{cor}[numberingbase]{Corollaire}
\newtheorem{lemma}[numberingbase]{Lemme}

\newtheorem{prop}[numberingbase]{Proposition}
\newtheorem{Q}[numberingbase]{Question}
\newtheorem{thm}[numberingbase]{Théorème}

\newtheoremstyle{claims}{0pt}{0pt}{}{}{\itshape}{.}{ }{}
\theoremstyle{claims}
\newtheorem{claim}[equation]{Claim}

\newtheoremstyle{defs}{0pt}{0pt}{}{}{\bfseries}{.}{ }{}
\theoremstyle{defs}
\newtheorem{defi}[numberingbase]{Définition}

\newtheorem{exa}[numberingbase]{Exemple}
\newtheorem*{exas}{Exemples}
\newtheorem{rmq}[numberingbase]{Remarque}
\newtheorem*{rmqs}{Remarques}

\Crefname{claim}{Claim}{Claims}
\Crefname{conj}{Conjecture}{Conjectures}
\Crefname{cor}{Corollaire}{Corollaires}
\Crefname{defn}{Définition}{Définitions}
\Crefname{eg}{Exemple}{Exemples}
\Crefname{prop}{Proposition}{Propositions} 
\Crefname{Q}{Question}{Questions}
\Crefname{rem}{Remarques}{Remarques}
\Crefname{theorem}{Théorème}{Théorèmes}
\Crefname{thm}{Théorème}{Théorèmes}
\Crefname{variant}{Variant}{Variants}

\theoremstyle{thms}
\newtheorem{thm-tweak}[subsection]{Théorème}
\Crefname{thm-tweak}{théorème}{théorèmes}
\newtheorem{lemma-tweak}[subsection]{Lemme}
\Crefname{lemma-tweak}{lemme}{Lemmes}
\newtheorem{cor-tweak}[subsection]{Corollaire}
\Crefname{cor-tweak}{corollaire}{corollaires}
\newtheorem{prop-tweak}[subsection]{Proposition}
\Crefname{prop-tweak}{proposition}{propositions} 
\newtheorem{conj-tweak}[subsection]{Conjecture}
\Crefname{conj-tweak}{Conjecture}{Conjectures} 

\theoremstyle{defs}
\newtheorem{defn-tweak}[subsection]{Définition}
\Crefname{defn-tweak}{définition}{définitions}
\newtheorem{eg-tweak}[subsection]{Exemple}
\Crefname{eg-tweak}{exemple}{exemples}
\newtheorem*{rmqs-tweak}{Remarques}
\newtheorem{rmq-tweak}[subsection]{Remarque}
\Crefname{rmq-tweak}{remarque}{remarques}

\newtheoremstyle{subsection-tweak}
   {0}
   {0pt}%
   {}
   {}%
   {\bfseries}
   {}%
   {.5em}
   {\thmnumber{\@{#1}{}\@{#2}.}%
    \thmnote{~{\bfseries#3.}}}

\theoremstyle{subsection-tweak}
\newtheorem{pp}[numberingbase]{}
\newcommand{\bpp}{\begin{pp}}
\newcommand{\epp}{\end{pp}}

\theoremstyle{subsection-tweak}
\newtheorem{pp-tweak}[subsection]{}

\renewcommand{\b}{\textbf}

\newcommand{\brems}{\begin{rmqs} \hfill \begin{enumerate}[label=\b{\thenumberingbase.},ref=\thenumberingbase]}
\newcommand{\remi}{\addtocounter{numberingbase}{1} \item}
\newcommand{\erems}{\end{enumerate}\end{rmqs}}
\newcommand{\bexas}{\begin{exas} \hfill \begin{enumerate}[label=\b{\thenumberingbase.},ref=\thenumberingbase]}

\newcommand{\eexas}{\end{enumerate} \end{exas}}

\newcommand{\bsm}{\begin{smallmatrix}}
\newcommand{\esm}{\end{smallmatrix}}
\newcommand{\blem}{\begin{lemma}}
\newcommand{\elem}{\end{lemma}}

\newcommand{\bconj}{\begin{conj}}
\newcommand{\econj}{\end{conj}}
\newcommand{\bprob}{\begin{Problem}}
\newcommand{\eprob}{\end{Problem}}

\newcommand{\bq}{\begin{Q}}
\newcommand{\eq}{\end{Q}}
\newcommand{\benum}{\begin{enumerate}[label={{\upshape(\alph*)}}]}
\newcommand{\benuma}{\begin{enumerate}[label={{\upshape(\arabic*)}}]}
\newcommand{\benumr}{\begin{enumerate}[label={{\upshape(\roman*)}}]}
\newcommand{\eenum}{\end{enumerate}}
\newcommand{\bitem}{\begin{itemize}}
\newcommand{\eitem}{\end{itemize}}
\newcommand{\bc}{}
\newcommand{\bdefi}{\begin{defi}}
\newcommand{\edefi}{\end{defi}}
\newcommand{\ov}{\overline}

\newcommand{\bexa}{\begin{exa}}
\newcommand{\eexa}{\end{exa}}

\newcommand{\bcl}{\begin{claim}}
\newcommand{\ecl}{\end{claim}}

\newcommand{\ba}{\begin{aligned}}
\newcommand{\ea}{\end{aligned}}
\newcommand{\be}{\begin{equation}}
\newcommand{\ee}{\end{equation}}
\newcommand{\ti}{\widetilde}
\newcommand{\bpf}{\begin{proof}}
\newcommand{\epf}{\end{proof}}
\newcommand{\bthm}{\begin{thm}}
\newcommand{\ethm}{\end{thm}}
\newcommand{\bprop}{\begin{prop}}
\newcommand{\eprop}{\end{prop}}
\newcommand{\bcor}{\begin{cor}}
\newcommand{\ecor}{\end{cor}}
\newcommand{\brem}{\begin{rmq}}
\newcommand{\erem}{\end{rmq}}
\newcommand{\Gm}{\Gamma}

\newcommand{\ra}{\rightarrow}
\newcommand{\hra}{\hookrightarrow}
\newcommand{\cH}{\mathcal{H}}

\newcommand{\cU}{\mathcal{U}}
\newcommand{\cC}{\mathcal{C}}

\newcommand{\cJ}{\mathcal{J}}

\newcommand{\cY}{\mathcal{Y}}

\newcommand{\bP}{\mathbb{P}}
\newcommand{\co}{\mathcal{O}}
\newcommand{\cV}{\mathcal{V}}
\newcommand{\cO}{\mathcal{O}}
\newcommand{\cS}{\mathcal{S}}
\newcommand{\cW}{\mathcal{W}}

\newcommand{\g}{\gamma}
\newcommand{\la}{\lambda}

\newcommand{\km}{\mathfrak{m}}
\newcommand{\kt}{\mathfrak{t}}

\newcommand{\cI}{\mathcal{I}}
\newcommand{\cZ}{\mathcal{Z}}
\newcommand{\kg}{\mathfrak{g}}

\newcommand{\kD}{\mathfrak{D}}

\newcommand{\bZ}{\mathbb{Z}}

\newcommand{\bql}{\overline{\mathbb{Q}}_{\ell}}

\newcommand{\kc}{\mathfrak{c}}
\newcommand{\kn}{\mathfrak{n}}
\newcommand{\kC}{\mathfrak{C}}
\newcommand{\wkC}{\widetilde{\mathfrak{C}}}
\newcommand{\ab}{\mathbb{A}}

\newcommand{\NN}{\mathbb{N}}
\newcommand{\al}{\alpha}
\newcommand{\bG}{\mathbb{G}}

\newcommand{\bC}{\mathbb{C}}
\newcommand{\bQ}{\mathbb{Q}}

\newcommand{\dhra}{\twoheadrightarrow}

\newcommand{\cF}{\mathcal{F}}
\newcommand{\cD}{\mathcal{D}}
\newcommand{\cX}{\mathcal{X}}
\newcommand{\La}{\Lambda}
\newcommand{\thra}{\twoheadrightarrow}
\title{Support singulier et homologie des fibres de Springer affines}
\author{Alexis Bouthier}
\begin{document}

\maketitle
\textbf{Abstract:}
We develop a theory of singular support for various infinite dimensional stacks and establish several functoriality properties. Then we apply this theory to compute the singular support of various affine character sheaves. We also explain how it can be used to study the local constancy conjecture of Goresky-Kottwitz-McPherson on the homology of affine Springer fibers along the root valuation stratification.
\hypersetup{
    linktoc=page,     
}
\tableofcontents

\section*{Introduction}
\subsection{Sur la théorie des faisceaux caractères}

Dans les années 80, Lusztig  introduit dans une série d'articles (\cite{Lus1}, \cite{Lus2}), la théorie des faisceaux caractères. Ce sont des faisceaux pervers $\ell$-adiques $G$- équivariants  sur $G$, où $G$ est un groupe connexe réductif sur un corps.
Dans \cite{MV}, Mirkovic et Vilonen calculent sur $\bC$ le support singulier des faisceaux caractères et obtiennent une caractérisation de ceux-ci à l'aide du cône nilpotent. Ils répondaient alors à une conjecture de Laumon et Lusztig.
La situation sur un corps fini a un analogue global, incarné par les faisceaux automorphes sur $\Bun_{G}$, le champ des $G$-torseurs sur une courbe projective lisse. Motivé par le cas fini, Laumon  conjecture alors dans (\cite{Lau}) que le support singulier des faisceaux automorphes est contenu dans le cône global nilpotent.
Cette conjecture a été  récemment démontrée par Arinkin-Gaitsgory-Kazhdan-Raskin-Rozenblyum-Varshavsky \cite{AGKRRV}.
Pour compléter le tableau, il s'agit de déterminer ce qu'il se produit dans le cas local.

Dans \cite{BKV}, l'auteur et Kazhdan-Varshavsky ont développé un formalisme général pour établir une théorie de Springer affine ainsi qu'une théorie des faisceaux pervers correspondante. De plus, dans \cite{B}, on montre qu'après application des coinvariants dérivés sous le groupe de Weyl affine $\widetilde{W}$, on obtient des faisceaux caractères sur $\Lie(G)((t))$.
La principale difficulté de ces travaux résident dans le fait que les espaces en jeu sont hautement de dimension infinie et jusqu'alors aucune théorie n'existait pour de tels objets. Une question naturelle alors qui se pose est de développer une théorie du support singulier en dimension infinie et de pouvoir calculer les supports singuliers de ces faisceaux caractères affines.
C'est l'objet de ce travail. 
De plus, on verra que la description du support singulier permet de déduire  des informations concrètes sur l'homologie des fibres de Springer affines. Enfin, on ne veut pas seulement traiter le cas des complexes mais travailler sur un corps arbitraire, puisqu'ensuite il s'agit de prendre des traces de Frobenius.

\subsection{Calcul pour des faisceaux caractères affines}
Pour étendre au cas affine, les énoncés de Mirkovic-Vilonen, on suppose momentanément dans cette section que l'on dispose d'une théorie du support singulier (cf. \ref{intro-SS} pour les fondements de cette théorie). On conserve les notations de la section précédente. La seule différence est que l'on considère $G$ défini sur un corps $k$ arbitraire algébriquement clos et l'on fixe un premier $\ell$ différent de la caractéristique. Dans la suite, on travaille en $\ell$-adique.

On considère sur  les foncteurs $k$-algèbres  $R\mapsto G(R[[t]])$ et $R\mapsto G(R((t)))$.
Ils sont respectivement représentables par un schéma affine, l'espace d'arcs $\clp G$ et par un ind-schéma ind-affine, l'espace de lacets $\cL G$.
Soit la flèche au niveau des faisceaux quotients pour la topologie étale, vus comme des $\infty$-champs:
\[f:[\cL\kg/\clp G]\ra [\cL\kg/\cL G].\]
C'est une flèche ind-fp-propre (cf. \ref{indfp} pour la définition) formellement lisse et on forme alors pour tout système local renormalisé $\cL$ sur $[\cL\kg/\clp G]$ \footnote{voir \ref{SSex-bis} pour la définition d'un système local renormalisé} :
\[\cS_{\kg,\cL}= f_{!}i_{*}\cL\] 
où $i:[\clp\kg/\clp G]\ra[\cL\kg/\clp G]$ est une immersion fermée présentation finie et $\omega_{[\clp\kg/\clp G]}$ est le faisceau dualisant.
On note alors $\kC\subset\cL\kg$ le sous-ind-schéma formé des éléments compacts, le complexe $\cS_{\kg}$ est supporté sur $[\kC/\cL G]$.
Le théorème est alors le suivant (\ref{micro-comp}):
\bthm\label{thm1}
Pour tout système local renormalisé $\cL$, le fermé cônique $C_{\kg}:=[(\kC\times\kC)\cap\{(\mu,\xi)\in\cL\kg\times\cL\kg,[\mu,\xi]=0\}/\cL G]$ est un microsupport pour $\cS_{\kg,\cL}$.
\ethm
Pour obtenir le théorème, on commence par établir que $\omega_{[\clp\kg/\clp G]}$ admet la section nulle comme support singulier et ensuite on utilise les propriétés de fonctoralité pour le support singulier.
Une fois établi ce premier énoncé, intéressant pour lui-même, On s'intéresse à la façon dont le suppport singulier va interagir avec la stratification par valuation radicielle introduite par Goresky-Kottwitz-McPherson.

Soit $(B,T)$ une paire de Borel de $G$, $W$ le groupe de Weyl, on suppose que la caractéristique est première à l'ordre de $W$, soit $\kc=\kt/W$, $\chi:\kg\ra\kc$ le morphisme de Chevalley et $\clp\kc^{\bullet}$ l'espace d'arcs, génériquement réguliers semisimples.
Goresky-Kottwitz-McPherson construisent alors une stratification:
\[\clp\kc^{\bullet}=\coprod_{w,r}\kc_{w,r},\]
où les $\kc_{w,r}$ sont des schémas affines, localement fermés, irréductibles et lisses \footnote{ici, par \og lisse\fg, on entend une limite projective filtrante de schémas lisses avec des flèches de transition affines lisses, surjectives} de $\clp\kg$, avec $w\in W$ et $r:R\ra\bQ_{+}$ une fonction sur l'ensemble des racines $R$.
De plus, pour  chaque $w\in W$, on a un tore maximal $T_{w}$ sur $k((t))$ qui admet un plongement $T_w\hra G$ unique à  $G$-conjugaison près et un schéma localement fermé $\kt_{w,r}\subset\clp\kt_w$, avec $\kt_w=\Lie(T_w)$ de telle sorte que l'on a un morphisme fini étale $\kt_{w,r}\ra\kc_{w,r}$.
 Soit $\cS^{\bullet}_{\kg,\cL}$ la restriction de $\cS_{\kg,\cL}$ à l'ouvert des génériquement réguliers semisimples et $\kg_{w,r}=\chi^{-1}(\kc_{w,r})$.
On montre alors l'énoncé suivant \ref{incl-loc}:
\bthm\label{int-inc}
Si l'ordre de $W$ est premier à la caractéristique de $k$, on a l'inclusion:
\[C_{\kg}\subset\coprod_{(w,r)}T_{[\kg_{w,r}/\cL G]}^{*}[\cL\kg/\cL G].\]
\ethm

\subsection{Sur la théorie du support singulier}\label{intro-SS}
\subsubsection{Enoncés principaux}
 Jusqu'à récemment, on ne disposait d'une théorie du support singulier pour des faisceaux constructibles que sur $\bC$ avec les travaux de Kashiwara-Schapira \cite{KS}. Inspiré par les travaux de Kashiwara-Schapira et la reformulation de la théorie du support singulier en termes de cycles proches, Beilinson et Saito (\cite{Bei}, \cite{Sai}) ont développé une théorie du support singulier pour des faisceaux étales pour des corps arbitraires.
Il s'agit donc de voir dans quel mesure ce formalisme peut se généraliser à des ind-schémas ou même à des champs tels que $[\cL\kg/\cL G]$ où $\cL G$ agit sur $\cL\kg$ par l'action adjointe.
\medskip

A la base du formalisme de Beilinson-Saito, il y a le cotangent d'un schéma lisse, la notion de locale acyclicité et la notion clé de paires tests.
Chacune de ces notions est délicate à définir en dimension infinie. Tout d'abord pour la lissité, elle va être remplacée par la formelle lissité, mais pour les objets que l'on manipule, la notion qui apparaît naturellement est celle où l'on demande un critère infinitésimal de relevèment localement pour la topologie étale. Il s'agit donc de comparer ces deux notions qui vont coïncider dans une généralité suffisamment vaste (Thm. \ref{etfl} et \ref{gab-fl}). Il faut ensuite une théorie du cotangent et plus généralement une théorie du complexe cotangent (sect. 1.5). On isole alors une classe d'$\infty$-champs pour lesquels on a une théorie du complexe cotangent, dits de Mittag-Leffler (\ref{def-ML}) qui contient par exemple les ind-schémas et les quotients de ind-schémas par des ind-schémas en groupes.
Pour définir la locale acyclicité, on a besoin de foncteurs cohomologiques pour des $\infty$-champs, le formalisme développé dans \cite{BKV} nous donne un bon cadre pour cela (\ref{mor-loc}) avec les propriétés de changement de base et de descente habituelles, ce sont des $\infty$-champs qui \og satisfont le recollement\fg~ (\ref{glu}).

Enfin, se pose la question du choix des paires-tests associées à un $\infty$-champ $\cX$ qui consistent en des paires  $(h,g):\cV\leftarrow\cW\ra\cX$. Pour celles-ci, le choix est déterminant et amène des théories radicalement différentes. Il faut à la fois accommoder le fait que pour des applications à la théorie de Springer affine, on veut des énoncés de fonctorialité du support singulier à la fois pour $f_*$ où $f$ est un morphisme ind-propre et pour $f^!$ dans le cas d'un morphisme schématique.
En particulier, pour les paires tests, on ne peut pas demander qu'elles consistent en des schémas si on veut la fonctorialité pour les morphismes ind-propres et en même temps pour pouvoir se raccorder à la théorie de Beilinson, on se rend compte que l'on ne peut pas demander que les morphismes soient arbitraires, il faut qu'ils soient schématiques.
Une fois que ces problèmes de fondements sont résolus, il s'agit d'établir les énoncés de fonctorialité par image inverse et image directe et de comparer la théorie obtenue à celle de Beilinson-Saito. Cela fait l'objet des énoncés \ref{propfonct}, \ref{fonct-clo} et \ref{fonct1}- \ref{SSex0}.
Il est à noter deux spécificités de la dimension infinie. Tout d'abord, on ne travaille qu'avec l'ensemble des microsupports et on ne cherche pas à montrer qu'il admet un plus petit élément, qui est le support singulier, on se contente d'une \og borne \fg~ et elle est déjà suffisante. De plus, il n'est pas clair qu'en général, le support singulier d'un faisceau sur un ind-schéma existe toujours. Par exemple, si $X$ est un ind-schéma ind-lisse, on peut montrer que le complexe dualisant $\omega_{X}$ admet la section nulle comme support singulier, mais déjà pour un ind-schéma tel que la grassmannienne affine, qui est formellement lisse, mais qui est une colimite de schémas singuliers, la question est ouverte.

La deuxième spécificité est donc de montrer que l'ensemble des microsupports est non-vide. Pour des schémas qcqs qui sont des limites projectives de schémas lisses et des faisceaux constructibles, cela va s'obtenir à l'aide du théorème de Beilinson (\ref{SSex0}), mais ensuite, il faut \og sortir\fg~ du monde des schémas pour aller dans celui des ind-schémas et là on se retrouve confronté à toute une série de subtilités d'algèbre différentielle sur des ind-schémas qui tient au fait que l'on considère des modules qui ne sont pas de pas de type fini sur des anneaux non-noethériens, ce qui fait de l'énoncé \ref{fonct-clo} sans doute le plus délicat. On utilise pour ce faire plusieurs résultats de Drinfeld et Gabber, non-publiés pour nombre d'entre eux.
On montre donc le théorème d'existence suivant:
\bthm
Soit un $\aleph_0$-ind-schéma ind-placide formellement lisse $X$, plongeable (cf.\ref{plong}). Soit $K\in\cD_{c}(X)$ borné (cf. \ref{borne}), alors $K$ admet un support singulier.
\ethm
La condition d'être plongeable est essentiellement que l'on peut trouver un ind-schéma $T$ ind-(placide formellement lisse) et une immersion fermée de présentation finie $X\hra T$ et que la catégorie de tels plongements est filtrante. Une telle condition est vérifiée par tous les ind-schémas de ind-type fini ainsi que les espaces de lacets de $k$-schémas affines lisses. Ici, on entend par borné que $K$ est supporté par un sous-schéma fermé de présentation finie.
\subsubsection{Relations avec d'autres travaux}
Dans la littérature, il y a au moins deux instances de calculs de support singulier sur des ind-schémas et qui peuvent être également vus comme une motivation et une justification de l'intérêt de développer une théorie du support singulier dans un cadre suffisamment général. Le premier travail est celui d'Evens-Mirkovic \cite[Thm. 0.1]{Ev-Mir} qui calcule le support singulier des cellules de Schubert dans la grassmannienne affine. Bien que le calcul ait lieu en dimension finie, l'objectif affirmé est d'obtenir une caractérisation microlocale des faisceaux $G(\cO)$-équivariants sur la grassmannienne affine.

  Le second travail est celui déjà mentionné de \cite{AGKRRV} où la preuve de la conjecture de Laumon repose sur un calcul de support singulier dans un champ de  Hecke \cite[sect. 20.5]{AGKRRV} qui fait intervenir des fibres de Springer affines. Le point particulièrement intéressant est que la conjecture de Laumon est globale sur des champs d'Artin, mais que sa preuve utilise des calculs locaux sur des ind-schémas. Il est cependant à noter que les auteurs se cantonnent à des ind-schémas de ind-type fini, ce qui simplifie considérablement la théorie des fasiceaux cohérents.

\subsection{Une conjecture de Goresky-Kottwitz-McPherson}\label{gkmconj}
Dans la théorie des fibres de Springer affines, on peut dégager deux conjectures principales, toutes deux formulées par Goresky-Kottwitz-McPherson.
La première porte sur la pureté des fibres de Springer affines  (\cite{GKM2}) et la deuxième sur la locale constance de l'homologie des fibres de Springer affines le long de la stratification par valuation radicielles (\cite[p. 1255]{GKM}). Nous n'avons rien à ajouter à la première, voyons comment l'énoncé \ref{thm1} peut éclairer la seconde. On commence par rappeler la conjecture de locale constance. On se place dans la situation où l'ordre de $W$ est premier à la caractéristique de $k$.
On considère la flèche:
\[f_{w,r}: X_{w,r}=\{(g,t)\in\Gr_{G}\times\kt_{w,r}, \ad(g)^{-1}\g\in\clp\kg\}\ra\kt_{w,r}\]
où $\Gr_{G}=\cL G/\clp G$ est la grassmannienne affine. La flèche $f_{w,r}$ est ind-propre et après réduction est localement de présentation finie d'après \cite[Thm. 4.3.3]{BKV}, la conjecture est la suivante \cite{GKM2}:
\bconj[Goresky-Kottwitz-McPherson]\label{loc-gkm}
Pour toute paire $(w,r)$, le faisceau $(f_{w,r})_{!}\omega_{X_{w,r}}$ est localement constant.
\econj
La conjecture n'est connue que dans le cas dit homogène, montré par Yun-Oblombkov comme corollaire de la pureté des fibres \cite[Cor. 5.4.4, 5.4.5]{YunIII}.
Dans le cas de $GL_n$, il est à noter que la stratification radicielle de Goresky-Kottwitz-McPherson correspond à la stratification par les séries de Puiseux de l'espace de modules des germes de courbes planes ou strates équisingulières. Dans ce cas, la paire $(w,r)$ revient à fixer le nombre de composantes irréductibles du germe ainsi que les valuations des racines du polynôme caractéristique (\cite{Tei}).
Elle admet une variante globale pour les fibres de Hitchin, mais qui a surtout étudiée dans le cas de $GL_n$ via le lien avec les jacobiennes compactifiées et le polynôme HOMFLY.
En effet, une version de la locale constance pour la caractéristique d'Euler a été établie par Maulik et Shende (\cite[Thm. 1.1, Cor. 6.7]{Mau}, \cite[Thm. A]{Sh} \cite{Mig}).

Si l'on était en dimension finie, sur $\bC$ et pour peu que la stratification soit de Whitney, il résulte d'un théorème de Kashiwara-Schapira \cite[Prop. 8.4.1]{KS}  que cette conjecture est impliquée par l'énoncé \ref{int-inc}:
\[SS(\cS^{\bullet}_{\kg})\subset\coprod_{(w,r)}T_{[\kg_{w,r}/\cL G]}^{*}[\cL\kg/\cL G].\]
La difficulté reste alors de voir comment déduire de \ref{int-inc} la conjecture \ref{loc-gkm}.
En suivant l'argument en dimension finie de Kashiwara-Schapira, cela se ramène à obtenir une formule pour le support singulier lorsque que l'on restreint à une strate. L'inconvénient est que cette strate n'est pas transverse au support singulier.
Il s'agit ainsi de montrer que la restriction $\cS_{w,r}$ aux strates $i_{w,r}:[\kg_{w,r}/\cL G]\hra[\cL\kg/\cL G]$ admet la section nulle comme support singulier.
Expliquons  quelles sont les difficultés.
 Tout d'abord, on restreint notre faisceau $\cS_{\kg}$ le long d'une strate qui n'est pas transverse au support singulier. Dans cette situation, il n'existe pas de formule générale pour calculer $SS(i_{w,r}^{!}\cS_{\kg})$, essentiellement à cause de problèmes de ramification sauvage. En revanche, si l'on se place sur $\bC$, Kashiwara-Schapira ont une formule pour $SS(i_{w,r}^{!}\cS_{\kg})$ en terme d'une opération $i^{\sharp}$ sur les fermés côniques. On peut alors vérifier dans le cas de $G=SL_2$ que l'opération $i^{\sharp}$ donne le bon résultat. La conjecture sur $\bC$ se ramène essentiellement à montrer la formule de Kashiwara-Schapira en dimension infinie.
\subsection{Organisation de l'article}
La première section consiste à rappeler les fondements $\infty$-catégoriques sur les $\infty$-champs ainsi que la définition dans un contexte suffisamment générale d'une notion de locale acyclicité, ainsi que les propriétés de fonctorialité associées. 
La deuxième section développe les outils de théorie des déformations, notamment une notion adaptée de formelle lissité,  pour avoir une théorie raisonnable du cotangent ainsi qu'une définition fonctorielle des notions de transversalité. Dans la troisième section, on isole une classe d'$\infty$-champs qui admettent une théorie du microsupport et on montre des énoncés d'existence de support singulier pour des schémas placides, des ind-schémas et de champs quotients (Thms. \ref{fonct-clo}, \ref{SSex0}, \ref{ex-ind2}, \ref{SS-indu}). La quatrième section traite des applications, notamment montre le calcul du support singulier des faisceaux de Grothedienck-Springer affine et le lien avec la conjecture de locale constance (Thms. \ref{micro-comp}, \ref{incl-loc}).
Enfin, dans deux appendices, nous établissons un certaines nombres de propriétés fondamentales sur les ind-schémas, notamment leur faisceau cotangent et les propriétés liées à la formelle lissité (Thms. \ref{ind-Tate}, \ref{gab-fl}) obtenues par Drinfeld et Gabber.

\subsection{Remerciements}
Le point de départ de ce travail  a été une suggestion de Michael Finkelberg qui nous a fait part d'un calcul de V. Ginzburg sur le support singulier du faisceau de Grothendieck-Springer usuel et de son lien avec la locale constance de la cohomologie des fibres de Springer. Je le remercie donc chaleureusement pour avoir partagé cette idée.
Tout au long de ce travail, j'ai  beaucoup bénéficié de mes entretiens avec Gérard Laumon et Eric Vasserot et de l'intérêt qu'ils ont pu porter à ce travail, je leur en suis très reconnaissant. Enfin, je souhaite remercier toutes les personnes avec qui j'ai pu également discuter ou correspondre, Alexandre Beilinson, Benjamin Hennion, Pierre Schapira, Vivek Shende, Bernard Teissier, Jean-Baptiste Teyssier et Jean-Loup Waldspurger.

\section{Faisceaux $\ell$-adiques sur des $\infty$-champs}
Dans cette section, on introduit les catégories de faisceaux $\ell$-adiques nécessaires pour établir une théorie du microsupport ainsi que les fonctorialités qui vont avec.
De plus, on développe une notion de morphismes localement acycliques pour une classe suffisamment générale d'$\infty$-champs. La deuxième section développe les outils de théorie des déformations pour avoir une bonne notion de cotangent et isoler une classe d'$\infty$-champs qui vont être admette une théorie du microsupport. On introduit également les propriétés de transversalité ainsi que les opérations de base sur celles-ci.

\subsection{Objets placides}
\subsubsection{Espaces algébriques placides}
\bdefi\label{fpl}
\benumr
\remi
 Soit un morphisme d'espaces algébriques $f:X\ra S$, il admet une présentation placide ou est placidement présenté (resp. fortement pro-lisse) s'il existe une présentation filtrante $X\simeq\varprojlim X_{\al}$, où les morphismes de transition sont affines lisses et $X_{\al}$ des $S$-espaces algébriques de présentation finie (resp. et de plus les $X_{\al}$ sont lisses sur $S$). 
\remi
Pour un morphisme d'espaces algébriques $f:X\ra S$ placide ou fortement pro-lisse, il est dit $\aleph_0$ si $f$ admet une présentation fortement pro-lisse indexée par un système projectif filtrant dénombrable.
\remi
Un morphisme d'espaces algébriques $X\ra Y$ est dit placide s'il existe $Y'\ra Y$ étale surjectif quasi-compact, tel que $X'=X\times_{Y}Y'\simeq\coprod_{\al}X_{\al}$ où $X_{\al}\ra Y$ est placidement présenté.
\remi
Si $S$ est le spectre d'un corps, on dit que $X$ est placide (resp. admet une présentation placide, resp. fortement pro-lisse) si la flèche $X\ra\Spec(k)$ l'est.
\eenum
\edefi
 \brems
\remi\label{fpl21}
D'après \cite[1.1.3]{BKV}, la classe des morphismes fortement pro-lisses est stable par changement de base, composition, morphisme fortement pro-lisses.
\remi\label{fpl22}
D'après \cite[11.1.6]{BKV}, si $f:X\ra S$ est placidement présenté, alors on a une présentation canonique $X\cong\varprojlim_{X\ra Y} Y$, où la limite parcourt les morphismes fortement pro-lisses $X\ra Y$ et les flèches de transition sont lisses affines.
\remi\label{fpl23}
Comme une colimite filtrante d'algèbres plates est plate, on obtient qu'un morphisme fortement pro-lisse est plat.
De plus, comme le complexe cotangent commute aux colimites filtrantes \cite[II.1.2.3.4]{Ill}, et que le complexe cotangent d'un morphisme lisse est quasi-isomorphe à un module projectif de type fini concentré en degré zéro, on obtient que pour un morphisme fortement pro-lisse $f:X\ra S$, $L_{X/S}$ est plat et concentré en degré zéro.
\erems

\bprop\label{plac}
Soit un morphisme d'espaces algébriques qcqs $f:X\ra S$ formellement lisse placidement présenté, alors il est fortement pro-lisse. En particulier, $L_{X/S}$ est localement projectif concentré en degré zéro.
Réciproquement si $X$ est $\aleph_0$-fortement pro-lisse, alors il est formellement lisse.
\eprop

\bpf
Soit une présentation placide $X\simeq\varprojlim X_{\al}$. Fixons $\al$, on a un triangle distingué de complexes cotangents:
\[p_{\al}^{*}L_{X_{\al}/S}\ra L_{X/S}\ra L_{X/X_{\al}},\]
avec $p_{\al}:X\ra X_{\al}$ qui est plate, fortement pro-lisse. Ainsi, $L_{X/X_{\al}}$ est plat et concentré en degré zéro.
De plus, comme $f$ est formellement lisse, on a d'après \cite[Tag. 0D10, Tag. 08RB]{Sta} et \cite[Tag. 0D0W, Tag. 08VE]{Sta}, $H^{-1}(L_{X/S})=0$ et $H^{0}(L_{X/S})=\Omega^{1}_{X/S}$ est localement projectif.
Par exactitude de $p_{\al}^*$ et la suite exacte de cohomologie, on déduit:
\[H^{-1}(p_{\al}^{*}L_{X_{\al}/S})=0~\text{et}~~ H^{0}(p_{\al}^{*}L_{X_{\al}/S})=p_{\al}^*\Omega^{1}_{X_{\al}/S}\]
 est localement projectif. A nouveau par exactitude de $p_{\al}^*$ et comme $\Omega^{1}_{X_{\al}/S}$ est cohérent, on obtient que dans un voisinage ouvert $U_{\al}$ de $p_{\al}(X)$, on a $H^{-1}(L_{U_{\al}/S})=0$ et $\Omega_{U_{\al}/S}$ est projectif de type fini. Ainsi, d'après \cite[Tag. 0D0L]{Sta}, $U_{\al}\ra S$ est lisse (mais pas nécessairement de présentation finie). 
Comme $p_{\al}:X\ra X_{\al}$ se factorise par $U_{\al}$ et que $X\simeq\varprojlim\limits_{\al'\geq\al} X_{\al'}$, il existe $\al'$ tel que $X\ra U_{\al}\subset X_{\al}$ se factorise en $X\ra X_{\al'}\stackrel{p_{\al'\al}}{\rightarrow} U_{\al}\subset X_{\al}$.
Comme $S$ est qcqs et $p_{\al'\al}$ est affine lisse, $U'_{\al}=p_{\al'\al}(X_{\al'})\subset U'_{\al}$ est un ouvert quasi-compact de $X_{\al}$ (on utilise la descente noethérienne \cite[Tag. 07SU]{Sta} pour se ramener au cas noéthérien où c'est clair). Donc quitte à remplacer $U_{\al}$ par $U'_{\al}$,  $X\ra S$ est fortement pro-lisse. Ainsi, $L_{X/S}$ est plat concentré en degré zéro et comme $\Omega^{1}_{X/S}$ est déjà localement projectif, $L_{X/S}$ l'est aussi.

Pour la réciproque, il suffit de montrer la formelle lissité, comme le système projectif est dénombrable, on applique le critère infinitésimal à chaque cran.
\epf
On a le renforcement suivant, dont on donne la preuve en appendice, cf \ref{ml2}:
\bprop\label{placbis}
Soit un morphisme d'espace algébriques qcqs $f:X\ra Y$ formellement lisse placidement présenté, alors il admet une présentation fortement pro-lisse $X\simeq\varprojlim X_{\al}$ et pour tout $\al$, $X\ra X_{\al}$ est formellement lisse.
\eprop
\blem\label{plac2}
Soit $S$ un schéma non-vide. Soit un morphisme de $S$-espaces algébriques $f:X\ra Y$  fortement pro-lisse ou formellement lisse, alors $df:f^{*}\Omega_{Y/S}^{1}\ra\Omega^{1}_{X/S}$ est universellement injectif.
\elem
\bpf
On a un triangle distingué de complexes cotangents:
\[f^{*}L_{Y/S}\ra L_{X/S}\ra L_{X/Y}.\]
Comme $f$ est fortement pro-lisse, $L_{X/Y}$ est quasi-isomorphe à $\Omega^{1}_{X/Y}$ placé en degré zéro et qui est plat sur $X$.
En passant à la suite exacte longue de cohomologie, on obtient ainsi une suite exacte:
\[0\ra f^{*}\Omega_{Y/S}^{1}\ra\Omega^{1}_{X/S}\ra\Omega^{1}_{X/Y}\ra 0\]
et $df$ est universellement injective par platitude de $\Omega^{1}_{X/Y}$.

Si $f$ est formellement lisse, l'argument est le même sauf que l'on utilise le fait que c'est seulement le tronqué $\tau_{\geq-1}L_{X/Y}$ qui est concentré en degré zéro et localement projectif d'après \cite[Tag. 0D10]{Sta}.
\epf
\bexa\label{arc}
Soit un $k[[t]]$-schéma affine de type fini $X$, on considère le foncteur sur les $k$-algèbres $R\mapsto X(R[[t]])$, alors il est représentable par un $k$-schéma qcqs, l'espace d'arcs $\clp X$ et $\clp X$ est fortement pro-lisse si $X$ est lisse sur $k[[t]]$, d'après \cite[3.1.1, 3.1.4.(e)]{BKV}.
\eexa
\subsubsection{Ind-schémas ind-placides}\label{ind-sc-ind-pl}
Soit un corps $k$, pour tout $k$-schéma qcqs $S$, on note $\Sch_{/S}$ la catégorie des $S$-schémas qcqs.
On appelle ind-schéma  sur $S$  un foncteur $T:(\Sch_{/S})^{op}\ra\Ens$ qui est représentable par une colimite filtrante $T\simeq\colim T_{\al}$ de $S$-schémas qcqs avec des flèches de transition qui sont des immersions fermées.
Si $T$ est un ind-schéma et si l'on désigne par $\cC_{/T}$ la catégorie des $k$-schémas qcqs $S$ avec une immersion fermée $S\hra T$, alors d'après \cite[Lem. 1.1.4]{Gai}, la catégorie est filtrante et l'on a une écriture canonique pour $T$:
\[T=\colim_{S\in\cC_{/T}}S.\]

\bdefi\label{ratio}
Un ind-schéma $T$ est \textsl{raisonnable} (resp. \textsl{ind-placide}) s'il s'écrit  $T\simeq\colim T_{\al}$ avec des flèches de transition qui sont des immersions fermées de présentation finie (resp. et les $T_{\al}$ sont des $k$-schémas placides qcqs). 
\edefi
L'exemple de référence d'un ind-schéma ind-placide est le suivant:
\bexa\label{loop}
Soit $k$ un corps et un $k((t))$-schéma $X$ affine de type fini, alors le foncteur $R\mapsto X(R((t)))$ est représentable par un ind-schéma ind-placide $\cL X$ (\cite[Thm. 6.3]{Dr}).
\eexa
De même, on a une notion de ind-espace algébrique et de ind-placidité pour ceux-ci. 
\subsection{$\infty$-champs}
\subsubsection{Généralités}
Soient une $\infty$-catégorie $C$ et $\Spc$, l'$\infty$-catégorie  des $\infty$-groupoïdes. On  considère la catégorie des préfaisceaux  $\PreSh(C)=\Fonct(C^{op},\Spc)$. Si $C$ est équipée d'une topologie de Grothendieck $\cT$, on considère $\Sh(C)$ la sous-$\infty$-catégorie des faisceaux pour la $\cT$-topologie.
Dans le cas où $C=\Aff_k$ est la catégorie des $k$-schémas affines munie de la topologie étale, on appelle catégorie des \textsl{préchamps} (resp. des $\infty$-\textsl{champs}) la catégorie $\PreSh(\Aff_k)$ (resp. $\St_k=\Sh(\Aff_k)$).
L'inclusion $\Sh(C)\hra\PrSh(C)$ admet un adjoint $\PrSh(C)\ra \Sh(C)$, donné par la faisceautisation \cite[6.2.2.7]{Lu1}.

On reprend la définition de \cite[1.2.4]{BKV}.
Soit (P) une classe de morphismes de $\infty$-champs $f:\cX\ra Y$ avec $Y\in\Aff_k$, stable par changement de base.

Un morphisme $f:\cX\ra\cY$ de $\infty$-champs vérifie (P), si pour tout changement de base $Y\ra \cY$, avec $Y$ un schéma affine, le changement de base $\cX\times_{\cY}Y\ra Y$ est dans (P).

On peut donc ainsi parler pour des morphismes de $\infty$-champs $f:\cX\ra\cY$ de morphismes propres/ de présentation finie/ind-fp-propres/ind-représentables, en prenant pour (P) la classe des morphismes de schémas $\cX\ra Y$ qui sont propres/de présentation finie, la classe des morphismes $\cX\ra Y$ ind-fp-propres avec $\cX$ un ind-espace algébrique et la classe des morphismes $\cX\ra Y$ avec $\cX$ un ind-espace algébrique, avec dans tous ces cas $Y$ un schéma affine.

\bexa\label{exquot}
Pour un $\infty$-champ $\cX$ muni d'une action d'un ind-schéma en groupes $H$, on peut former l'$\infty$-champ quotient $[\cX/H]$.
De plus d'après \cite[Lem. 2.3.3]{B}, la flèche $\cX\ra [\cX/H]$ est ind-représentable.
Si de plus $H$ est  ind-(quasi-affine), la flèche $\cX\ra [\cX/H]$ est ind-schématique.
\eexa
\brem
Par ind-représentable, on entend que pour tout $\Spec(A)\ra [\cX/H]$, le changement de base est représentable par un ind-espace algébrique.
\erem

\bdefi\label{indfp}
\benumr
\item
Soit un morphisme $f:T\ra Y$ d'un ind-espace algébrique vers un schéma affine, il est dit \textsl{ind-fp-propre}, s'il admet une présentation comme colimite filtrante $T\simeq\colim T_{\al}$ où les $T_{\al}\ra Y$ font fp-propres et les morphismes de transition sont des immersions fermées.
\item
Soit un morphisme $f:\cX\ra Y$ d'un $\infty$-champ vers un schéma affine, il est dit \textsl{localement ind-fp-propre}, s'il existe $Y'=\Spec(A')\ra Y$ étale surjectif tel que $\cX\times_{Y}Y'\ra Y'$ est ind-fp-propre.
\eenum
\edefi

Souvent dans les applications, il est facile d'obtenir des morphismes localement ind-fp-propres.
\bexa\label{loopeq}
Soit un morphisme d'$\infty$-champs $f:\cX\ra\cY$ localement ind-fp-propre et équivariant pour l'action d'un $\infty$-champ en groupes $H$, alors $f:[\cX/H]\ra [\cY/H]$ est localement ind-fp-propre (\cite[1.2.9]{BKV}).

Une illustration de cet énoncé général est la suivante; soit un $k$-ind-schéma $X$ et un groupe $G$ connexe réductif sur un corps $k$ tel que $\cL G$ agit sur $X$. Alors, on peut former les $\infty$-champs $[X/\clp G]$ et $[X/\cL G]$ et la flèche induite $[X/\clp G]\ra[X/\cL G]$ est localement ind-fp-propre, puisque la projection $X\times_k[\cL G/\clp G]\ra X$ est ind-fp-propre puisque la grassmannienne affine $[\cL G/\clp G]$ l'est.
\eexa
On a en fait le lemme de descente suivant \cite[Lem. 2.3.7]{B}:
\blem\label{desc-ind-fp}
Un morphisme de $\infty$-champs localement ind-fp-propre est ind-fp-propre.
\elem

A tout $\infty$-champ $\cX\cong\colim\limits_{S\ra\cX} S$, on associe son réduit $\cX_{red}\cong\colim\limits_{S\ra\cX}S_{red}$.
On a donc une flèche canonique $\cX_{red}\ra \cX$.

\'{E}tant donné un morphisme d'$\infty$-champs $f:\cX\ra\cY$, on dit que c'est une \textsl{équivalence topologique}, s'il induit un isomorphisme $f_{red}:\cX_{red}\stackrel{\sim}{\rightarrow}\cY_{red}$.
D'après \cite[1.5.1]{BKV}, cette notion est locale pour la topologie étale sur la base, stable par pullbacks, composition et passage au quotient.

Pour tout morphisme de $\infty$-champs $i:\cY\ra X$, on dit que $i$ est une \textsl{immersion topologiquement constructible} si pour tout schéma affine $X\ra \cX$, la composée :
\[(\cY\times_{\cX}X)_{red}\ra\cY\times_{\cX}X\ra X\]
est une immersion localement fermée de présentation finie.
\subsubsection{Espace topologique sous-jacent}
A tout $\infty$-champ $\cX$, on associe son espace topologique $[\cX]$ de la façon suivante:
\begin{itemize}
	\item 
	L'ensemble sous-jacent est l'ensemble des classes d'équivalences de paires $(K,[z])$ où $K/k$ est une extension de corps et $[z]\in\pi_{0}(\cX(K))$ avec $(K,[z])\sim (K',[z'])$ s'il existe $K''\supset K,K'$ tel que $[z]$ et $[z']$ ont la même image dans $\pi_{0}(\cX(K''))$.
	\item
	Un sous-ensemble $U\subset [\cX]$ est ouvert si $U=[\cU]$ où $\cU$ est un sous-préchamp ouvert de $\cX$.
\end{itemize}
Un morphisme $f:\cX\ra\cY$ de $\infty$-champs induit une application continue $[f]:[\cX]\ra[\cY]$ d'espaces topologiques.
D'après \cite[Lem. 2.2.2]{BKV}, un sous-ensemble $U\subset\cX$ est ouvert si pour tout morphisme $f:X\ra\cX$ avec $X$ un schéma affine, $[f]^{-1}(U)$ est ouvert.

\subsection{Catégorie de faisceaux $\ell$-adiques}
\subsubsection{Définition}\label{l-ad}
Soit $\ell$ un premier différent de la caractéristique de $k$ et $\AlgSp^{tf}_k$ (resp. $\AlgSp_k$) la catégorie des $k$-espaces algébriques de type fini (resp. des $k$-espaces algébriques qcqs) que l'on munit de la topologie étale.
\begin{itemize}
	\remi
Soit $\Catst$  l'$\infty$-catégorie des petites $\infty$-catégories, stables et $\bql$-linéaires avec des foncteurs exacts, i.e. qui préservent les colimites finies, comme morphismes. D'après \cite[1.1.4.4, 1.1.4.6]{Lu2}, elle contient toutes les petites colimites filtrantes et toute les petites limites.
\remi
Soit $\PrCat$ la sous-$\infty$-catégorie des  $\infty$-catégories présentables (\cite[5.5.0.1]{Lu1}), stables, $\bql$-linéaires, avec des $\infty$-foncteurs continus, i.e. les catégories sont stables par petites colimites et tous les foncteurs commutent aux petites colimites.
Elle contient toutes les petites limites et colimites d'après \cite[4.2.4.8, 5.5.3.13, 5.5.3.18]{Lu1}.
\end{itemize}
Tous les foncteurs qui apparaissent ici sont des $\infty$-foncteurs et les limites et colimites sont à prendre au sens homotopique.
Etant donnée $\cC\in\Catst$, on peut former sa catégorie des ind-objets $\Ind(\cC)$ (\cite[5.3.5.1]{Lu1}) qui est stable (\cite[1.1.3.6]{Lu2}) et donc présentable (elle admet toutes les colimites finies et les colimites filtrantes, donc toutes les colimites). On a ainsi un foncteur naturel :
\[\Ind:\Catst\rightarrow\PrCat\]
donné par $\cC\mapsto \Ind(\cC)$, qui commute aux petites colimites filtrantes (\cite[5.3.5.10]{Lu1}, \cite[1.9.2]{GRI} et \cite{Roz}).
Pour un $k$-espace algébrique de type fini $Y$, on dispose, d'après  Liu-Zheng (\cite{LZ1}, \cite{LZ2}), d'une $\infty$-catégorie $\cD_{c}(Y):=\cD^{b}_{c}(Y,\bql)$ dont la catégorie homotopique est $D^{b}_{c}(Y,\bql)$.
Dans la suite, on utilise librement le formalisme des six foncteurs pour les $\infty$-catégories (cf. \cite{LZ1},\cite{LZ2}).
Il y a un foncteur naturel:
\begin{equation}
\cD_{c}:(\AlgSp^{tf}_k)^{op}\rightarrow\Catst
\label{tfdef1}
\end{equation}
qui  à chaque morphisme $f:X\rightarrow Y$ associe $f^{!}:\cD_{c}(Y)\rightarrow\cD_{c}(X)$.
On définit alors $\cD:=\Ind\circ\cD_{c}$. En faisant l'extension de Kan à gauche et à droite, on obtient des foncteurs :
\begin{equation}
\cD_{c}:\PreSh(\AlgSp_k)\rightarrow\Catst
\label{kan1}
\end{equation}
\begin{equation}
\cD:\PreSh(\AlgSp_k)\rightarrow\PrCat
\label{kan2}
\end{equation}
On note $\cD_{\bullet}$ pour désigner les foncteurs $\cD$ et $\cD_c$.
De manière explicite si $X\simeq\varprojlim X_{\al}$ où les $X_{\al}$ sont des $k$-espaces algébriques de type fini avec des morphismes de transition affines, on a $\cD_{\bullet}(X)\simeq\colim_{f^!}\cD_c(X_{\al})$ et si $T=\colim T_{\al}$ est un ind-espace algébrique avec $T_{\al}\in\AlgSp_k$, on a $\cD_c(T)=\varprojlim_{f^!}\cD_{\bullet}(T_{\al})$.
D'après \cite[5.2.8]{BKV}, $\cD_{\bullet}$ est un faisceau pour la topologie étale, donc ils se factorisent par $\Sh(\AlgSp_k)$ via le morphisme de faisceautisation.
En particulier, comme le plongement $\iota:\Aff_k\hra\AlgSp_k$ induit une équivalence de $\infty$-catégories par pullback $\iota^{*}:\St_k=\Sh(\Aff_k)\stackrel{\sim}{\rightarrow}\Sh(\AlgSp_k)$, de telle sorte que l'on peut voir $\cD_{\bullet}$ comme un foncteur sur $\St_k$.

\brem\label{kan-tor}
Si l'on fait la même opération avec les coefficients de torsion $\bZ/n$ ($n$ premier à la caractéristique), il résulte de \cite[Tag. 09YU]{Sta} que les foncteurs obtenus sur la catégorie des $k$-espaces algébriques qcqs revient à considérer l'$\infty$-catégorie des faisceaux étales constructibles de $\bZ/n$-modules sur un espace algébrique qcqs et ensuite de prendre la Ind-complétion.
\erem
\subsubsection{Fonctorialité}\label{ss-fonct}
Pour tout morphisme de $\infty$-champs $f:\cX\ra\cY$, on dispose d'un foncteur $f^!$ pour $\cD_{\bullet}$.
De plus, pour toute équivalence topologique $f:\cX\ra\cY$, on a d'après \cite[5.3.6]{BKV} une équivalence de catégories:
\begin{equation}
f^{!}:\cD_{\bullet}(Y)\stackrel{\sim}{\rightarrow}\cD_{\bullet}(\cX).
\label{topeq}
\end{equation}
\bprop\label{Base}
Soit $f:\cX\ra\cY$ un morphisme ind-fp-propre d'$\infty$-champs (resp. une équivalence topologique).
Alors,  $f^!$ admet un adjoint à gauche $f_!$ et	pour tout carré cartésien entre $\infty$-champs:
$$\xymatrix{\ti{\cX}\ar[r]^{\ti{h}}\ar[d]_{\ti{f}}&\cX\ar[d]^{f}\\\ti{\cY}\ar[r]^{h}&\cY}$$
l'application de changement de base:
 \begin{equation}
\ti{f}_!\ti{h}^{!}\ra h^{!}f_!
\label{basemor}
\end{equation} est une équivalence.
\eprop
\bpf
Si $f$ est ind-fp-propre, l'énoncé est obtenu dans \cite[Prop. 5.3.7]{BKV}. Si $f$ est une équivalence topologique, $\ti{f}$ en est une et $f^{!}$ et $\ti{f}^{!}$ sont des équivalences dont les adjoints à gauche $f_!$ et $\ti{f}_!$ sont les inverses.
Il suffit de montrer que \eqref{basemor} est une équivalence, après application de $\ti{f}^{!}$. On a d'une part $\ti{f}^{!}\ti{f}_!\ti{h}^{!}\cong\ti{h}^{!}$ et d'autre part $\ti{f}^{!} h^{!}f_!\cong \ti{h}^{!}f^{!}f_{!}\cong\ti{h}^!$, ce qu'on voulait.
\epf

Pour traiter les hensélisations qui apparaissent dans la définition de morphismes localement acycliques, on a besoin considérer des morphismes pro-étales:
\blem\label{pro-et}
Soit $f:\cX\ra\cY$ un morphisme pro-étale quasi-compact d'$\infty$-champs, alors $f^!$ admet un adjoint continu à droite $f_*$.
De plus, pour tout carré cartésien entre $\infty$-champs:
$$\xymatrix{\ti{\cX}\ar[r]^{\ti{h}}\ar[d]_{\ti{f}}&\cX\ar[d]^{f}\\\ti{\cY}\ar[r]^{h}&\cY}$$
l'application de changement de base $h^{!}f_*\ra \ti{f}_*\ti{h}^!$ est une équivalence.
\elem
\brem
Ici, par pro-étale, on entend qu'après tout changement de base $S\ra\cY$ avec $S$ affine, le morphisme $f_{S}:\cX\times_{\cY}S\ra S$ admet une présentation comme limite projective filtrante de morphismes étales quasi-compacts avec des morphismes de transition affines.
\erem
\bpf
Comme $f^{!}$ commute aux colimites, d'après \cite[5.5.2.9]{Lu1}, $f^{!}$ admet un adjoint à droite $f_{*}$.
On commence par vérifier le changement de base avant d'en déduire la continuité.
Le même dévissage que \cite[5.3.7 Step.1]{BKV}, permet de se ramener à montrer le changement de base dans le cas où $\ti{\cY}$ et $\cY$ sont des schémas affines et alors tous les objets considérés sont des schémas qcqs. En considérant une présentation pro-étale $X\simeq\varprojlim X_{\al}$ avec $X_{\al}\ra Y$ affine étale, on a alors $\cD(X)\simeq\colim_{\pi^!}\cD(X_{\al})$. Comme les morphismes de transition sont affines étales, d'après l'adjonction de Lurie \cite[5.5.3.3]{Lu1}, on a :
\[\colim\limits_{\pi^!}\cD(X_{\al})\simeq\varprojlim_{\pi_*}\cD(X_{\al}).\]
En particulier, l'énoncé de changement de base se ramène au cas où $f$ est affine étale et cela se déduit alors de \cite[5.2.7.(a)]{BKV}.

Montrons maintenant que $f_*$ est continu.
D'après \cite[Prop. 5.5.7.2]{Lu1}, si $\cD(\cX)$ et $\cD(\cY)$ sont compactement engendrées, $f_{*}$ est continu si et seulement si $f^!$ préserve les objets compacts.
Si $Y$ est un espace algébrique qcqs, $\cD(Y)$ est compactement engendrée avec $\cD(Y)=\Ind(\cD_c(Y))$ (\cite[5.2.4.(c)]{BKV}) et $f^!$ préserve $\cD_c$ par construction, donc $f_{*}$ est continu.

Passons au cas général, comme $\cD(\cY)=\varprojlim\limits_{\al_{Y}:Y\rightarrow\cY}\cD(Y)$ où $Y$ est un schéma affine. Commes les $\al_{Y}^{!}$ sont continus, il suffit de voir que $\al_{Y}^{!}f_*$ est continu.
Or, d'après ci-dessus, le morphisme de changement de base:
\[\al_{Y}^{!}f_*\simeq f_{*}\al_{X}^{!}\]
 avec $\al_{X}:X=Y\times_{\cY}\cX\ra\cX$ est un isomorphisme. Le foncteur de droite est clairement continu d'après le cas des espaces algébriques qcqs, ce qui conclut.
\epf
Pour la suite, on a besoin de considérer une classe intermédiaire entre les ind-schémas et les $\infty$-champs, qui contient les objets tels que $[\cL X/\cL G]$ pour $G$ connexe réductif sur $k$ qui agit sur $X\in\Aff_k^{tf}$.
Tout d'abord, pour toute immersion topologiquement constructible de $\infty$-champs $\eta:\cY\hra\cX$, on dispose d'un foncteur $\eta_*$ (\cite[5.4.4]{BKV}).
On s'intéresse à l'existence d'un adjoint à gauche.
\bdefi\label{glu}
Soit $\cX$ un $\infty$-champ, on dit qu'il satisfait le recollement si pour toute immersion topologiquement constructible $\cY\hra\cX$, le foncteur $\eta_{*}$ admet un adjoint à gauche $\eta^*$.
\edefi 

\brem\label{glustk}
Les $\infty$-champs suivants satisfont le recollement:
\benumr
\item
Si $\cX$ satisfait le recollement alors pour toute immersion topologiquement constructible $\eta:\cY\hra\cX$, $\cY$ satisfait le recollement (\cite[Lem. 5.5.5.(a)]{BKV}).
\item
Si $\cX\simeq\colim \cX_{\al}$ où chaque $\cX_{\al}$ admet le recollement et les morphismes de transition sont des immersions ouvertes quasi-compactes (resp. des immersions fermées de présentation finie) alors $\cX$ satisfait le recollement (\cite[Lem. 5.5.5.(b)]{BKV}).
\item
D'après \cite[5.5.6]{BKV}, on obtient que tout espace algébrique placide satisfait le recollement ainsi que tout ind-espace algébrique ind-placide d'après (ii).
\item
D'après \cite[Prop. 5.5.7]{BKV}, les quotients $[X/H]$ d'un ind-espace algébrique ind-placide $X$ par un groupe ind-placide $H$, i.e. un objet en groupe dans la catégorie des ind-espaces algébrique ind-placides satisfont le recollement. Cela s'applique donc à des quotients  $[\cL Y/\cL G]$ pour $G$ affine lisse sur $k$ qui agit sur $Y\in\Aff_k^{tf}$.
\eenum
\erem

\subsection{Morphismes localement acycliques}\label{mor-loc}
Soit $f:\cX\rightarrow S$ un morphisme de $\infty$-champs, avec $\cX$ qui satisfait le recollement et $S$ un $k$-schéma de type fini.
Soient $s$ un point géométrique de $S$ et $t$ un point géométrique de l'hensélisé $S_{(s)}$, $j_{t}:\cX\times_{S}t\rightarrow \cX$ et $i_{s}:\cX_{s}\rightarrow\cX$, qui sont des composées d'immersions fermées de présentation finie et de morphismes pro-étales quasi-compacts comme $S$ est de type fini. Soit $K\in\cD(\cX)$, on écrit:
\[\Psi^{f}_{t\mapsto s}(K)=i_{s}^{*}(j_{t,*}j_t^{*}K),\]
où les foncteurs sont bien définis comme $\cX$ satisfait le recollement et par \ref{pro-et}.

\bdefi\label{def-locac}
On dit que $K\in\cD(\cX)$ est \textsl{localement acyclique} par rapport à $f$ si pour tout point géométrique $s$ de $S$ et tout point géométrique $t\in S_{(s)}$, la flèche canonique:
\begin{equation}
i_{s}^*K\rightarrow \Psi^{f}_{t\mapsto s}(K)
\label{locac}
\end{equation}
est un isomorphisme. 
\edefi
\brem\label{equiv-locac}
Si $X$ est un $k$-schéma qcqs et si $K\in\cD(X,\bZ/n)$ avec $n$ premier à la caractéristique, pour $f:X\ra S$ un  morphisme de $k$-schémas qcqs,  Deligne, dans \cite[Exp. 7.2.12]{SGA4-1/2}, définit une notion de locale acyclicité par rapport à $f$ si pour tout point géométrique $x$ de $X$ et tout point géométrique $t$ de $S_{(f(x))}$, la flèche canonique :
\begin{equation}
R\Gm(X_{(x)},K)\ra R\Gm((X_{(x)})_{t},K)
\label{l-ac}
\end{equation}
est un isomorphisme avec $(X_{(x)})_{t}=X_{(x)}\times_{S_{(f(x))}}t$.
Ainsi, en passant aux fibres en tout point géométrique  $x$ de la fibre $X_{s}$ dans \eqref{locac} et d'après \ref{kan-tor}, on obtient que la définition \ref{def-locac} est équivalente à celle de Deligne.
Les deux définitions sont également équivalentes pour des $\bql$-coefficients sur un $k$-schéma de type fini.
\erem

\blem\label{acprop}
Soient $f:\cT'\rightarrow \cT$ un morphisme  ind-fp-propre (resp. une équivalence topologique) d'$\infty$-champs satisfaisant le recollement, $g:\cT\rightarrow S$ avec $S$ un $k$-schéma de type fini et $K\in\cD(\cT')$, on suppose que	$K$ est localement acyclique par rapport à $g\circ f$,
alors $f_!K$ est localement acyclique par rapport à $g$.
\elem
\bpf
Soient $s$ un point géométrique de $S$ et $t$ un point géométrique de $S_{(s)}$, $i_s:\cT_s\rightarrow \cT$ et $i'_s:\cT'_s\rightarrow \cT$.
Il s'agit de voir que l'application canonique:
\begin{equation}
i_{s}^*f_!K\rightarrow \Psi^{g}_{t\mapsto s}(f_{!}K)
\label{canloc}
\end{equation}
est un isomorphisme.
La flèche $s\ra S$ se décompose en $s\ra s_0\ra\ov{\{s_0\}}\hra S$, où $s_0\in S$ est le point sur lequel s'envoie $s$. On a donc une composée de morphismes pro-étales \footnote{L'inclusion d'un point générique $\eta\hra\overline{\{\eta\}}$ est même pro-Zariski.} pour lesquels $g^!\simeq g^*$ d'après \ref{pro-et} et d'une immersion fermée de présentation finie. Il en est donc de même de $i_{s}$ et d'après \ref{Base} et \cite[Cor. 5.5.4]{BKV}, on a un isomorphisme $i_{s}^*f_!K\cong f_!i^{',*}_{s} K$.
Pareillement, on a le même type de décomposition pour $j_{t}$ et d'après \ref{commut}, on a une flèche canonique $f_!j_{t,*}\ra j_{t,*}f_! $ qui est un isomorphisme \footnote{On a besoin de \ref{commut} seulement pour traiter la partie pro-étale, c'est immédiat pour la partie immersion fermée de présentation finie.}.  On obtient ainsi:
\[\Psi^{g}_{t\mapsto s}(f_{!}K)\cong f_!(\Psi^{g\circ f}_{t\mapsto s}(K)).\]
En particulier, \eqref{canloc} s'identifie à :
\[f_!(i^{',*}_{s}K\rightarrow \Psi^{g\circ f}_{t\mapsto s}(K)),\]
qui est donc un isomorphisme comme $g\circ f$ est localement acyclique.
\epf
\blem\label{commut}
Considérons un carré cartésien:
$$\xymatrix{\cU'\ar[r]^{j'}\ar[d]_{f'}&\cT'\ar[d]^{f}\\\cU\ar[r]^{j}&\cT},$$
où $f$ est ind-fp-propre ou une équivalence topologique et $j$ pro-étale quasi-compact alors, on a une flèche canonique $j_{*}f'_!\ra f'_!j'_*$ qui est un isomorphisme.
\elem
\bpf
Construisons d'abord la flèche; on a un morphisme canonique $f_!j^{*}j_{*}\ra f_!$. Maintenant, si $f$ est ind-fp-propre ou une équivalence topologique, par changement de base, on a  une flèche $j^{*}f_!j_*\stackrel{\sim}{\rightarrow}f_!j^{*}j_{*}\ra f_!$, soit par adjonction une flèche
\begin{equation}
 f_!j_*\ra j_{*}f_!.
\label{morcomm}
\end{equation}
Il s'agit maintenant de voir que c'est un isomorphisme. Si $f$ est une équivalence topologique, comme $f_!$ est une équivalence, $f^{!}$ est aussi un adjoint à gauche de $f_!$, donc \eqref{morcomm} est un isomorphisme. 

Supposons maintenant $f$ ind-fp-propre, comme $\cT\cong\colim_{\Spec(A)\ra\cT}\Spec(A)$ et comme $f_!$ et $j_*$ satisfont le changement de base d'après \ref{Base}, \ref{pro-et} et \cite[5.5.4]{BKV}, on se ramène d'abord à $\cT=\Spec(A)$.
Comme les foncteurs sont continus, donc commutent aux colimites, on est alors ramené au cas où $f$ est fp-propre entre espaces algébriques qcqs et cela se déduit de l'existence d'un adjoint à gauche $f_*$ d'après \cite[5.2.5(c)]{BKV} de telle sorte que $f_!\cong f_*$.
\epf
\blem\label{pulLis}
\benumr
\item
Pour $f:\cX\ra\cY$ étale quasi-compact entre $\infty$-champs satisfaisant le recollement, $g:\cY\ra Z$ où $Z$ lisse de type fini, si $K\in\cD(\cY)$ est $g$-acyclique  alors $f^{!}K$ est $g\circ f$ acyclique et il y a équivalence si $f$ est surjectif.
\item
Soit $f:Y'\rightarrow [Y/H]$  fortement pro-lisse avec $Y,Y'$ des espaces algébriques et $H$ un schéma en groupes placide.
On suppose que $Y$ est également un ind-espace algébrique ind-placide.
 Soit $g:[Y/H]\ra Z$ avec $Z$ lisse de type fini et $K\in\cD([Y/H])$  alors si $K$ est localement acyclique par rapport à $g$, $f^!K$ est localement acyclique par rapport à $g\circ f$.
\item
 Si de plus, $f$ est surjectif et admet des sections localement pour la topologie étale, $f^{!}K$ est $g\circ f$-acyclique si et seulement si $K$ est $g$-acyclique.
\eenum
\elem
\bpf
(i) Cela se résume à vérifier des identités de changement de base pour vérifier que \eqref{locac} commute avec le passage à  $f^!$, et comme $f$ est étale cela se déduit d'après \ref{pro-et} de l'identité $f^{!}\cong f^{*}$ et du changement de base \cite[5.4.5]{BKV}. La conservativité de $f^{!}$ et la descente étale assurent que \eqref{locac} est une équivalence si elle l'est après pullback par $f^{!}$.

On montre (ii) et (iii) simultanément. On commence par le cas où $H=\{1\}$. 
Tout d'abord, d'après (i), si $h:\ti{Y}\ra Y$ est étale quasi-compact surjectif, on a $h^{!}K$ $g\circ h$-acyclique si et seulement si $K$ est $g$-acyclique.
Comme $Y$ est un espace algébrique, d'après \cite[Tag. 04NN]{Sta}, il admet une base de voisinage par des ouverts quasi-compacts. Comme de plus $Y$ est ind-placide, on en déduit que l'on peut écrire $Y=\bigcup U_{\al}$ où les $U_{\al}$ sont des ouverts quasi-compacts placides \footnote{chaque $U_{\al}$ est quasi-compact dans un fermé placide, donc placide lui-même.}.
Après localisation Zariski et recouvrement étale, on se ramène donc au cas où $Y$ est placidement présenté et alors $Y'$ est aussi placidement présenté.
On considère la flèche :
\[i_{s}^*K\rightarrow \Psi^{g}_{t\mapsto s}(K).\]
On applique alors $f^!$ et comme $f$ est fortement pro-lisse, on a par changement de base \cite[5.2.7.(d),5.4.5]{BKV} et \ref{pro-et}:
\[i_s^{*}f^!K\cong f^!(i_{s}^*K)\stackrel{\sim}{\rightarrow} f^{!}\Psi^{g}_{t\mapsto s}(K)\cong\Psi^{g\circ f}_{t\mapsto s}(f^!K).\]
Ainsi si $K$ est localement acyclique par rapport à $g$,  $f^!K$ est localement acyclique par rapport à $g\circ f$.
Si de plus, $f$ est surjectif et localement trivial pour la topologie étale, alors $f^!$ est conservatif d'après \cite[5.3.1.(d)]{BKV}.

Passons au cas général. Comme $Y\ra [Y/H]$ est trivial, localement pour la topologie étale, il existe alors un recouvrement étale quasi-compact $\rho:Y''\ra Y'$ tel qu'on a un diagramme commutatif:
$$\xymatrix{Y''\ar[r]\ar[d]_{\rho}&Y\ar[d]\\Y'\ar[r]^-{f}&[Y/H]}$$
Ainsi, d'après le cas des espaces algébriques, $f^{!}K$ est $g\circ f$-acyclique si et seulement si $\rho^{!}f^{!}K$ est $g\circ f\circ\rho$-acyclique.
On se ramène alors au cas $Y''=Y'$ et $f$ qui se factorise via $\pi:Y\ra [Y/H]$. A nouveau, le cas des espaces algébriques montre qu'il suffit de montrer l'assertion pour $\pi$.Il alors suffit de vérifier que pour une immersion fermée de présentation $i: [X/H]\hra\cY=[Y/H]$ et le diagramme cartésien:
$$\xymatrix{X\ar[d]_-{\ti{\pi}}\ar[r]^{\tilde{i}}&Y\ar[d]_{\pi}\\[X/H]\ar[r]^-{i}&[Y/H]},$$
 on a un isomorphisme:
\[\ti{i}^{*}\pi^{!}\cong \ti{\pi}^!i^{*}.\]
ce qui est l'objet de \cite[Lemme. 5.5.8]{BKV}, une fois que l'on passe aux bar-complexes.
Enfin, l'argument de conservativité reste valable dans le cas général, ce qui conclut.
\epf
Avec les grosses catégories, pour montrer que \eqref{locac} est un isomorphisme, il n'est a priori pas suffisant de vérifier l'isomorphisme sur les fibres. En revanche cela ne se produit pas pour des faisceaux constructibles:
\bprop\label{check-ac}
Soit un ind-schéma ind-placide $V$, $K\in\cD_{c}(V)$ supporté par un sous-schéma fermé de présentation finie, soit $f:V\ra Y$ un morphisme vers un $k$-schéma de type fini, alors $K$ est localement acyclique pour $f$, si $K_{x}$ est $f_{x}$-acyclique pour tout point géométrique $\bar{x}\in\supp(K)$, avec $f_{x}:V_{(x)}\ra V\ra Y$, où $V_{(x)}$ est le localisé strict en $x$.
\eprop
\bpf
En effet, soit un sous-schéma fermé de présentation finie $i:Z\subset V$ sur lequel $Z$ est supporté de telle sorte que $K=\la_{*}\la^{!}K_{0}$. Pour montrer que \eqref{locac} est un isomorphisme, on se ramène immédiatement à l'énoncé pour $K_0$ et $V=Z$ et comme le problème est local pour la topologie étale d'après \ref{pulLis}, on peut supposer $Z$ placide affine. Dans ce cas, un complexe constructible est nul si et seulement si ses fibres sont nulles.
Comme les foncteurs $i_{s}^{*}$ et $\Psi^{f}_{t\mapsto s}$ préservent la constructibilité et commutent à la localisation étale, le cône de la flèche est constructible et par hypothèses ses fibres sont nulles donc il est nul.
\epf
\section{Théorie de la déformation et formelle lissité}
Le but de la section est  de définir le faisceau cotangent et plus généralement le pro-complexe cotangent pour une classe assez générale de préchamps, dont on aura besoin pour définir une notion de transversalité.
Si l'on considère un ind-schéma $X$ sur un corps $k$, comme il est  colimite filtrante de $k$-schémas qcqs, 
on dispose du faisceau pro-quasi-cohérent cotangent $\Omega^{1}_{X/k}$, obtenu par le système projectif $(i_{S,*}\Omega_{S/k}^{1})$ pour $i_S:S\hra X$ qui parcourt la catégorie filtrante des sous-schémas fermés de $X$.

Pour de tels objets, on peut définir le ind-schéma cotangent sans avoir recours à une théorie générale. En revanche, on a également besoin de traiter des champs quotients tels que $[X/H]$ avec $X$ un ind-schéma et $H$ un ind-schéma en groupes. De tels champs ne s'expriment pas comme colimite filtrante de schémas qcqs et il n'est pas clair a priori qu'ils admettent un $\Omega^{1}$ raisonnable. On se  place donc dans le cadre plus général des préchamps et même des préchamps dérivés pour formuler les conditions nécessaires à l'existence d'un pro-complexe cotangent (cf. \ref{defT}). Une fois ceci fait, pour avoir de vrais espaces cotangents qui sont des espaces vectoriels et non des pro-espaces vectoriels, il faudra introduire une notion de Mittag-Leffler (cf. sect. \ref{TMod}).
\subsection{Pro-complexe cotangent}
 
\subsubsection{Pro-catégories}
Soit $\cC$ une $\infty$-catégorie, elle est dite \textsl{accessible} (\cite[5.4.3]{Lu1}), si elle est petite et stable par facteur direct.
Soit $\cC$ une $\infty$-catégorie accessible avec des limites finies, on définit la pro-catégorie $\Pro(\cC)^{op}$ comme la catégorie des foncteurs accessibles $F:\cC\ra\Spc$ qui commutent aux limites finies (\cite[A.8.1.1]{Lu3}).
D'après Lurie \cite[A.8.1.2, A.8.1.3]{Lu3}, on a $\Pro(\cC)=\Ind(\cC^{op})^{op}$ et on a un foncteur pleinement fidèle $\cC\hra\Pro(\cC)$.

De plus, si $\cC$ est stable, $\Pro(\cC)$ est stable (\cite[1.1.1.13, Prop. 1.1.3.6]{Lu2}). Enfin, si $\cC$ est équipée d'une $t$-structure, elle induit une $t$-structure sur $\Pro(\cC)$ en définissant $\Pro(\cC)^{\leq 0}=\Pro(\cC^{\leq 0})$ et $\Pro(\cC)^{\geq 0}=\Pro(\cC^{\geq 0})$ (\cite[6.1.2.(a)]{BKV}).
\subsubsection{Faisceaux quasi-cohérents}
Comme indiqué en début de section, même lorsque l'on traite le cas des schémas, le bon problème de modules à formuler pour construire le complexe cotangent est celui de la géométrie dérivée. Soit $\dAff\supset\Aff$, la catégorie des schémas affines dérivés. Sa catégorie opposée  $\dAff^{op}$ consiste en l'$\infty$-catégorie des anneaux simpliciaux où l'on localise par les équivalences faibles.
On définit également la catégorie des préchamps dérivés $\dPrStk=\PreSh(\dAff)$ et on note $\Cat_{st}$ l'$\infty$-catégorie des $\infty$-catégories stables qui admet toutes les  petites colimites.
Enfin, il est commode de considérer la sous-catégorie pleine $\Ind(\dAff)\subset\dPrStk$ qui consiste en les ind-objets de $\dAff$ (sans condition sur les morphismes de transition).
Par Yoneda \cite[5.1.3.1]{Lu1}, on a un foncteur pleinement fidèle $\eta:\dAff\ra\dPrStk$ qui se factorise en une série de foncteurs pleinement fidèles:
\begin{equation}
\xymatrix{\dAff\ar[r]_-{\eta_1}\ar@/^2pc/[rr]^{\eta}&\Ind(\dAff)\ar[r]_{\eta_{2}}&\dPrStk}.
\label{faith}
\end{equation}

Soit $Z\in\dAff_k$, on définit $\QCoh(Z)$ l'$\infty$-catégorie stable dont la catégorie homotopique est la catégorie dérivée quasi-cohérente sur $Z$ \cite[1.3.5.8]{Lu2} et $\QCoh(Z)^{\heartsuit}$ le coeur de $\QCoh(Z)$ pour la $t$-structure standard (\cite[1.3.5.16, 1.3.5.21]{Lu2}). Pour tout morphisme $f:Z'\ra Z$ dans $\dAff_k$, on a un foncteur $f^{*}:\QCoh(Z)\ra\QCoh(Z')$, on a donc des foncteurs:
\[\QCoh:\dAff_k^{op}\ra\Cat_{st}.\]
En faisant l'extension de Kan à gauche, on obtient un foncteur:
\[\QCoh:\dPrStk^{op}\ra\Cat_{st}.\]
En particulier, si $\cX$ est un préchamp dérivé, on a par définition
\[\QCoh(\cX)=\varprojlim\limits_{f:\Spec(A)\ra\cX,f^{*}}\QCoh(\Spec(A)).\]
D'après \cite[Ch.3, 1.5]{GRI}, elle est également munie d'une $t$-structure telle que :
\[\QCoh(\cX)^{\leq 0}=\{\cF\in\QCoh(\cX), \forall~ S\stackrel{x}{\rightarrow}\cX, x^*\cF\in\QCoh(S)^{\leq 0}\}.\]
En revanche, la catégorie $\QCoh(\cX)^{\geq 0}$ n'admet pas une description aussi agréable et nous n'en ferons pas usage.
D'après \cite[6.2.3.1]{Lu3} et \cite[2.3.3]{GRI}, $\QCoh$ est un faisceau pour la topologie fpqc.
Pour tout $Z\in\dAff_k$, l'$\infty$-catégorie $\QCoh(Z)$ est accessible avec des limites finies, donc on peut considérer $\Pro(\QCoh(Z))$ (resp. $\Pro(\QCoh(Z)^{\leq 0})$) que l'on étend ensuite à nouveau par extension de Kan à gauche en des foncteurs:
\[\Pro(\QCoh), \Pro(\QCoh)^{\leq 0}:\dPrStk^{op}\ra\Cat_{st}.\]
\subsubsection{Foncteur des dérivations}
Pour définir une théorie du pro-complexe cotangent, deux approches sont possibles, celles de Gaitsgory-Rozenblyum \cite[Ch.I]{GRII} et celle de Hennion \cite{He} (dans le cas des ind-champs). Dans la première approche, étant donné un préchamp $\cX$, on commence par définir pour tout schéma affine $x:S\ra\cX$, ce que serait la restriction  $x^*L_{\cX/\bZ}$ et ensuite on s'assure d'une compatibilité par rapport aux différents pullbacks pour obtenir un objet dans $\Pro(\QCoh(\cX))$. Dans la deuxième approche, on définit directement un foncteur des dérivations au niveau des préchamps et on demande qu'il soit pro-représentable. On adopte ce point de vue dans cette section.
 On tire la définition suivante de \cite[Déf. 1.1.15]{Hen}.
\bdefi
Pour toute $\infty$-catégorie $\cC$ avec des colimites finies, on considère le foncteur $\Fact_{\cC}:\cC^{\Delta^{1}}\ra\infty-\Cat$
qui à tout morphisme $\phi:x\ra y$ associe l'$\infty$-catégorie des factorisations $x\ra c\ra y$ de $\phi$.
\edefi
Soit un schéma affine $S=\Spec(A)$, pour tout $A$-module $M$, on considère le schéma $S[M]=\Spec(A\oplus M)$ avec comme produit:
\[(a,m).(a',m')=(aa',am'+a'm),\]
de telle sorte que  l'on  a une factorisation $S\hra S[M]\ra S$, où la première flèche est une extension infinitésimale de carré nul. Cette construction fournit un foncteur $S[-]:\QCoh(S)^{\heartsuit}\ra \Fact_{\Aff}(\Id_{S})$, fonctorielle pour $S\in\Aff$.
 En particulier, on obtient une transformation naturelle:
\[[-]:\QCoh_{\Aff}^{\heartsuit}\ra \Fact_{\Aff}(\Id_{-}).\]
(N.B: On note $\Aff$ en indice pour indiquer que l'on regarde le foncteur $\QCoh$ sur la catégorie $\Aff_k$).
D'après Lurie \cite[25.3.1.1]{Lu3}, il existe une unique extension à $\dAff$ et à $\QCoh^{\leq 0}$ qui commute aux colimites:
\[[-]:\QCoh_{\dAff}^{\leq 0}\ra \Fact_{\dAff}(\Id_{-}).\]
Appliquant le foncteur $\Pro$, on obtient une transformation naturelle:
\[[-]:\Pro(\QCoh_{\dAff}^{\leq 0})\ra\Pro(\Fact_{\dAff}(\Id_{-}))\cong\eta_{1}^{*}(\Fact_{\Ind(\dAff)}(\Id_{-}))\ra \eta^{*}(\Fact_{\dPrStk}(\Id_{-})),\]
où $\eta_1^{*}$ et $\eta^{*}$ sont les morphismes au niveau des catégories de foncteurs, induits par la restriction à $\dAff$ (cf. \eqref{faith}).
On applique alors l'extension de Kan à gauche $\eta_!$ qui est  un adjoint à gauche de $\eta^{*}$ et comme $\eta$ est pleinement fidèle, la counité $\eta_!\eta^{*}\ra\Id$ est une équivalence.
On obtient ainsi le foncteur:
\begin{equation}
[-]:\Pro(\QCoh_{\dPrStk}^{\leq 0})\ra\eta_!\eta^{*}(\Fact_{\dPrStk}(\Id_{-}))\cong\Fact_{\dPrStk}(\Id_{-}).
\label{foncultim}
\end{equation}
On peut maintenant faire la définition suivante:
\bdefi
Soit un morphisme $\cX\ra\cY$, on définit le foncteur des dérivations:
\[\Der_{\cY}(\cX,-)=\Map_{\cX/./\cY}(\cX[-],\cX).\]
\benumr
\remi
On dit que $\cX\ra\cY$ admet un pro-complexe cotangent relatif, s'il existe $L_{\cX/\cY}\in\Pro(\QCoh(\cX))$ tel que pour tout $E\in\Pro(\QCoh^{\leq 0}(\cX))$, on a:
\[\Der_{\cY}(\cX,E)\cong\Map(E,L_{\cX/\cY}).\]
\remi
Si $\cY=\Spec(\bZ)$, on dit que $\cX$ admet un pro-complexe cotangent si $\cX\ra\Spec(\bZ)$ admet un pro-complexe cotangent relatif.
\eenum
\edefi
\brems
\remi
Dans la suite, on utilisera cette notion que pour des préchamps classiques (i.e. non-dérivés), mais il est plus naturel de formuler cette définition dans ce cadre plus général.
\remi
Si l'on dispose d'un morphisme $f:\cX'\ra\cX$ entre préchamps dérivés qui admettent un pro-complexe cotangent, par fonctorialité de la construction, on a un morphisme :
 \begin{equation}
f^{*}L_{\cX/\bZ}\ra L_{\cX'/\bZ},
\label{canL}
\end{equation}
de telle sorte qu'en appliquant $\Map(-,E)$ pour tout $E\in\Pro(\QCoh^{\leq 0}(\cX'))$, $f$ admet un complexe cotangent relatif et 
\begin{equation}
L_{\cX'/\cX}\cong\Cofib(f^{*}L_{\cX/\bZ}\ra L_{\cX'/\bZ}).
\label{Omful}
\end{equation}
\erems
\bexa\label{Tdef-alg}
Soit $X$ un ind-espace algébrique, alors d'après  \cite[Prop. 1.2.19]{Hen}, il admet un pro-complexe cotangent $L_{X/\bZ}$ et si $X\simeq\colim X_{\al}$ alors pour tout morphisme $x\in X(A)$:
\[x^{*}L_{X/\bZ}\simeq\varprojlim x^{*}L_{X_{\al}/\bZ}.\] 
En particulier, on a  $L_{X/\bZ}\in\Pro(\QCoh(X)^{\leq 0})$.
\eexa
Dans la suite, on a également besoin de savoir si les quotients $[X/H]$ d'un ind-schéma par un ind-schéma en groupes admettent un pro-complexe cotangent. Pour ce faire, on a besoin d'un énoncé de descente qui nécessite d'une part des hypothèses supplémentaires sur les ind-schémas et une théorie de la déformation plus \og robuste\fg.

Plus précisément,  si $\cX$ est un préchamp dérivé qui admet un complexe cotangent, alors $L_{\cX/\bZ}$ contrôle les extensions scindées de carré nul. Toutefois, pour avoir une théorie de la déformation raisonnable qui contrôle  toutes les extensions $\cX\hra\cX'$ de carré nul, il faut de plus que $\cX$ soit infinitésimalement cohésif (\cite[Ch. I, sect.6]{GRII}). De manière informelle, si $S\hra S'$ est une extension de carré nul de schémas affines, cette condition décrit $\cX(S')$ en termes de $\QCoh(S)$.
On a la définition suivante tirée de \cite[Ch. I, 7.1.2]{GRII}:
\bdefi\label{defT}
Soit $\cX\ra\cY$ un morphisme de préchamps dérivés, on dit qu'il admet une théorie de la déformation relative si:
\benumr
\remi
$\cX$ est convergent, i.e pour tout $S\in\dAff$, $\Map(S,\cX)\cong\varprojlim_{n\in\NN}\Map(\tau^{\leq n}S,\cX)$ où pour tout $n\in\NN$, $\tau^{\leq n}S$ est le $n$-tronqué de $S$.
\remi
$\cX\ra\cY$ admet un pro-complexe cotangent relatif.
\remi
$\cX$ est infinitésimalement cohésif.
\eenum
Si $\cY=\Spec(\bZ)$, on dit que $\cX$ admet une théorie de la déformation.
\edefi
\brems
\remi
Tout préchamp $n$-tronqué pour un certain $n\in\NN$ est convergent, en particulier dans la suite, on ne considère que des préchamps non-dérivés, la première condition est automatiquement remplie.
\remi\label{rem-defT1}
Si $X$ est un ind-schéma, alors il admet une théorie de la déformation \cite[Ch. II. Prop. 1.3.2]{GRII} et plus généralement  d'après \eqref{Omful}, tout morphisme de ind-schémas $X\ra Y$ admet une théorie de la déformation relative.
Plus généralement, si $X$ est un ind-espace algébrique, alors d'après \cite[17.3.1.7]{Lu3}, il est infinitésimalement cohésif et il admet un complexe cotangent, comme un espace algébrique admet un complexe cotangent et qu'une colimite filtrante d'objets qui admettent un complexe cotangent admet un pro-complexe cotangent (pour la même raison que \cite[Ch. II. Prop. 1.3.2]{GRII}).
\remi\label{rem-defT2}
Soit $f:\cY'\ra\cY$ un morphisme de préchamps (non-dérivés), alors d'après \cite[Lem. 4.2.5, 6.1.4]{GRII} $f$ admet une théorie de la déformation relative, si pour tout schéma affine $S\ra\cY$, $\cY'\times_{\cY}S$ admet une théorie de la déformation relative.
\erems
A partir de maintenant, tous les objets qui apparaissent sont non-dérivés.
Pour obtenir un énoncé de descente du complexe cotangent, on a besoin d'une notion de formelle lissité, ce qui fait l'objet du paragraphe suivant.
\subsection{Formelle lissité}\label{sect-fl}
On commence par donner une définition plus générale que la notion usuelle de formelle lissité, mais plus flexible pour les applications. 
De plus, on verra qu'elle coïncide avec la notion usuelle dans une généralité suffisante (cf.Thm. \ref{etfl}).

\bdefi
\benumr
\item
Soit un morphisme de $\infty$-champs $f:\cX\ra\cY$, il est \textsl{formellement lisse} s'il vérifie le critère infinitésimal localement pour la topologie fppf, i.e. pour toute $k$-algèbre $R$ et tout idéal $I\subset R$ avec $I^2=0$, il existe $\Spec(R')\ra\Spec(R)$ fidèlement plat tel qu'on a un diagramme commutatif:
$$
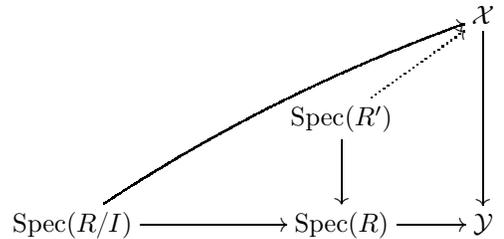
\xymatrix{&&&\cX\ar[dd]\\&&\Spec(R')\ar@{.>}[ur]\ar[d]\\\Spec(R/I)\ar@/^/[uurrr]\ar[rr]&&\Spec(R)\ar[r]&\cY}.$$
\item
On dit que $f$ est formellement lisse standard si l'on peut prendre $R=R'$.
\item
Un $\infty$-champ sur $k$ est \textsl{formellement lisse}  si la flèche $\cX\ra\Spec(k)$ l'est.
\eenum
\edefi
\brem
A vrai dire, dans tout ce travail, les objets considérés vérifient le critère infinitésimal après un recouvrement étale $\Spec(R')\ra\Spec(R)$.
Néanmoins, à la lumière du théorème suivant il est plus commode de considérer la définition plus générale ci-dessus.
\erem
\bthm[Gabber]\label{etfl}
Soit un morphisme formellement lisse  $X\ra\Spec(A)$ de ind-espaces algébriques ind-placides, alors il est formellement lisse standard.
\ethm
\brem
On montre cet énoncé dans l'appendice (cf. Thm. \ref{gab-fl}). L'énoncé est vrai également si l'on remplace placide par seulement raisonnable, mais est beaucoup plus long à démontrer et nous n'en aurons pas l'usage dans ce travail.
\erem
Le théorème précédent nous permet d'obtenir le critère de lissité formelle suivant, essentiellement dû à Gaitsgory-Rozenblyum:
\bthm\label{fl-gr}
Soit un $k$-ind-schéma ind-placide $X$ (cf \ref{ind-sc-ind-pl}), alors il est formellement lisse si et seulement si pour tout morphisme $x:S\ra X$ avec $S\in\Aff_k$, on a:
\benumr
\item
$H_{1}(x^*L_{X/k})=0$.
\item
$H_{0}(x^{*}L_{X/k})$ est pro-projectif.
\eenum
On peut donc résumer cet énoncé en disant que $\tau_{\geq -1}L_{X/k}$ est quasi-isomorphe à $\Omega^{1}_{X/k}$, localement pro-projectif et placé en degré zéro.
\ethm
\brem
On rappelle que $\Omega^{1}_{X/k}$ est le pro-faisceau $(i_{S,*}\Omega^{1}_{S/k})$ où $S$ parcourt la catégorie filtrantes des sous-schémas fermés de $X$ et $i_{S}:S\hra X$.
\erem
\bpf
On utilise \ref{etfl} avec $A=k$ pour obtenir que $X$ est formellement lisse standard et on utilise ensuite la caractérisation de \cite[9.4.2]{GR}.
\epf
\bprop\label{flquot}
Soit un $k$-ind-schéma $X$ avec une action d'un ind-schéma en groupes $H$ ind-placide et formellement lisse, alors la projection $\pi:X\ra [X/H]$
est formellement lisse standard. Si de plus, $X$ est formellement lisse, alors $[X/H]$ est formellement lisse.
\eprop
\bpf
On considère un morphisme $\Spec(R)\ra[X/H]$. On forme $E_{R}:=X\times_{[X/H]}\Spec(R)\ra\Spec(R)$ qui d'après \ref{exquot} est un $H$-torseur et un ind-espace algébrique.
\'{E}tale localement, $E$ est isomorphe à $H\times_k\Spec(R)$, en particulier comme $H$ est ind-placide, formellement lisse, il en est de même pour $E$ d'après \ref{etfl} et \cite[Tag. 02KZ]{Sta}, ce qui conclut.
Enfin, si l'on a une paire $(R,I)$ avec $I^2=0$ et une section $x\in[X/H](R/I)$, on la relève après recouvrement étale $\Spec(B)\ra\Spec(R/I)$ en $x'\in X(B)$ et il résulte de  \cite[Exp. III, Cor. 6.8]{SGA1} qu'il existe un recouvrement étale $\Spec(B')\ra\Spec(R)$ qui induit modulo $I$ le précédent, donc par formelle lissité de $X$, $x'\in X(B)$ se relève en $\ti{x}'\in X(B')$, ce qui conclut.
\epf
\brem\label{desc-fl}
Le même argument que ci-dessus donne que, tout morphisme formellement lisse $\cX\ra\cY$ entre $\infty$-champs $H$-équivariant pour un préchamp en groupes $H$, induit un morphisme formellement lisse $[\cX/H]\ra[\cY/H]$.
\erem
La formelle lissité est particulièrement utile pour obtenir l'énoncé de descente suivant \cite[Ch. I, Prop. 7.4.2]{GRII}:
\bprop\label{desc-GR}
Soit $f:\cY\ra\cX$ un morphisme de préchamps, on suppose que :
\benumr
\item
$\cX$ satisfait la descente étale.
\item
$f$ admet des sections localement pour la topologie étale.
\item
$\cY$ est formellement lisse standard et admet une théorie de la déformation.
\item
$f$ admet une théorie de la déformation relative.
\eenum
Alors, $\cX$ admet une théorie de la déformation.
\eprop
\brem
Pour un morphisme de préchamps, on entend que $f$ admet des sections localement pour la topologie étale si c'est le cas après pullback à tout schéma affine $S\ra\cX$.
\erem
L'énoncé précédent nous permet alors de construire des complexes cotangents équivariants:
\bprop\label{cot-eq}
Soit $X$ un $k$-ind-schéma ind-placide formellement lisse avec une action d'un ind-schéma en groupes $H$, alors le quotient $\pi:X\ra[X/H]$ admet une théorie de la déformation et on a une suite fibrante:
\[\pi^{*}L_{[X/H]/k}\ra L_{X/k}\ra L_{X/[X/H]}.\]
\eprop
\brem
Si l'on a la généralisation de \ref{etfl} à un ind-schéma raisonnable, on peut remplacer dans la proposition ind-placide par raisonnable.
\erem
\bpf
Comme $X$ est un ind-schéma ind-placide formellement lisse, d'après \ref{etfl}, il est formellement lisse standard et d'après \ref{rem-defT1} il admet une théorie de la déformation.
Le morphisme $X\ra[X/H]$ admet des sections localement pour la topologie étale et $[X/H]$ satisfait la descente étale. Enfin, d'après \eqref{exquot}, la flèche $X\ra[X/H]$ est ind-représentable, donc d'après \eqref{rem-defT1} et \eqref{rem-defT2}, il admet une théorie de la déformation relative. On déduit donc de \ref{desc-GR} que $[X/H]$ admet une théorie de la déformation et il résulte de \eqref{Omful} que l'on a une suite fibrante:
\[\pi^{*}L_{[X/H]/k}\ra L_{X/k}\ra L_{X/[X/H]}.\]
\epf

\subsubsection{Suites fibrantes pour des immersions fermées}
\bprop\label{exa-clo}
On considère une immersion fermée de présentation finie $i:Z\ra X$ entre ind-schémas ind-placides formellement lisses, soit $\cI$ le pro-faisceau d'idéaux qui définit $i$, alors on a une suite exacte:
\[0\ra C_{Z/X}=i^{*}\cI/\cI^{2}\ra i^{*}\Omega^{1}_{X/k}\ra\Omega^{1}_{Z/k}\ra 0\]
de pro-faisceaux localement pro-projectifs et $i^{*}\cI/\cI^{2}$ est pro-projectif de type fini.
\eprop
\bpf
On considère $x:\Spec(A)\ra X$.  Etant donné que $\Omega^{1}_{X/k}$ correspond au pro-faisceau $(i_{S,*}\Omega^{1}_{S/k})$ où $S$ parcourt la catégorie filtrantes des sous-schémas fermés de présentation finie de $X$, comme $X$ ind-placide et que 
\begin{equation}
C_{Z/X}=i^{*}\cI/\cI^{2}=(\ti{i}_{S,*}C_{Z_{S}/S}).
\label{cone-def}
\end{equation}
avec $\ti{i}_{S}: Z\times_{S}X\hra S\hra Z$, de présentation finie comme $Z\hra X$ l'est, de telle sorte que $C_{Z/X}$ est un système projectif filtrant de faisceaux de présentation finie.
De plus, pour une factorisation $S\hra S'\hra X$, on obtient un carré cartésien:
$$\xymatrix{Z_{S}\ar[d]\ar[r]^{\ti{i}_{SS'}}&Z_{S'}\ar[d]\\S\ar[r]^{i_{SS'}}&S'}$$
où toutes les flèches sont des immersions fermées de présentation finie et d'après \cite[Tag. 0473]{Sta}, on a une surjection:
\begin{equation}
\tilde{i}_{SS'}^{*}C_{Z_{S'}/S'}\twoheadrightarrow C_{Z_{S}/S},
\label{cone-surj}
\end{equation}
 donc $\cI/\cI^{2}$ est un pro-ensemble de Mittag-Leffler.
On a immédiatement par réduction au cas des schémas et passage au pro-système une suite:
\begin{equation}
x^{*}\cI/\cI^{2}\ra x^{*}i^{*}\Omega^{1}_{X/k}\ra x^*\Omega^{1}_{Z/k}\ra 0.
\label{Om-exa}
\end{equation}
D'autre part, on  a une suite fibrante au niveau des pro-complexes cotangents:
\[x^{*}L_{Z/X}\ra x^{*}i^{*}L_{X/k}\ra x^{*}L_{Z/k}.\]
Ainsi, en passant à l'homologie et comme $X$ et $Z$ sont formellement lisses, il résulte de \ref{fl-gr} et \eqref{Om-exa} que $H_{0}(x^{*}L_{Z/X})=0$ et une suite exacte:
\begin{equation}
0\ra H_{1}(x^{*}L_{Z/X})\ra i^{*}\Omega^{1}_{X/k}\ra\Omega^{1}_{Z/k}\ra 0.
\label{cotan-clo}
\end{equation}
De plus, par définition de la $t$-structure sur les pro-quasi-cohérents,  $H_{1}(x^{*}L_{Z/X})$ se calcule à partir du pro-système $H_{1}(x^{*}L_{Z_{S}/S})$ pour $S$ fermé dans $X$ et donc par réduction au cas des schémas, d'après \cite[Tag. 08RG]{Sta}, on a :
\[H_{1}(x^{*}L_{Z/X})=x^{*}\cI/\cI^2.\]
Enfin, comme les deux derniers termes de \eqref{cotan-clo} sont pro-projectifs, il en est de même du premier d'après \ref{proj-rem2} et $\cI/\cI^{2}$ est bien localement pro-projectif de type fini d'après \ref{flat-ML}.
\epf
La preuve précédente et plus spécifiquement \eqref{cone-def} et \eqref{cone-surj} permettent d'obtenir le corollaire suivant:
\bcor\label{cone-surj2}
Soit une immersion fermée de présentation finie $i:Z\ra X$ entre ind-schémas ind-placides, soit $h:V\ra X$ un morphisme de ind-schémas raisonnables, alors on a une surjection:
\[\ti{h}^{*}C_{Z/X}\thra C_{Z_{V}/V}\]
avec $\ti{h}:Z_{V}=Z\times_{X}V\ra Z$.
\ecor
Dans la suite, on veut remplacer les faisceaux pro-quasi-cohérents par des quasi-cohérents topologiques.
\subsection{Champ cotangent}
\subsubsection{Fibrés généralisés}
On commence par donner une construction générale, dont le cotangent est un cas particulier.
Soit $\cX$ un préchamp et $\Omega\in\Pro(\QCoh(\cX)^{\heartsuit})$, i.e. pour tout $z\in\cX(A)$, on dispose de $\Omega_{A,z}\in\Pro(\QCoh(\Spec(A))^{\heartsuit})$  et pour tout $f:\Spec(B)\ra\Spec(A)$, on a un isomorphisme $\phi_{B}:\Omega_{B,z\circ f}\cong f^{*,\heartsuit}\Omega_{A,z}$.
On construit alors le champ $\ab(\Omega)\ra\cX$ défini de la façon suivante:
pour tout $\phi\in\cX(A)$, on pose:
\[\ab(\Omega)_{\vert A}=\Sym((\Omega_{A,z})^{\vee}):=\colim\limits_{\Omega_{A,z}\ra\cF} \Sym(\cF^{\vee})\]
où $\cF$ parcourt la catégorie cofiltrante des faisceaux quasi-cohérents qui reçoivent un morphisme de $\Omega_{A,z}$.
Pour toute factorisation $\Spec(A)\stackrel{h}{\rightarrow}\Spec(C)\stackrel{z_{C}}{\rightarrow}\cX$, on vérifie alors aisément que l'on dispose de morphismes de transition $\ab(\Omega)_{\vert C}\times_{\Spec(C)}\Spec(A)\cong \ab(\cF)_{\vert A}$ et on définit ensuite :
\begin{equation}
\ab(\Omega)=\colim_{\Spec(A)\ra \cX}\ab(\Omega)_{\vert A}\ra \cX.
\label{cstr-fib}
\end{equation}
On applique maintenant cette construction dans le cas du cotangent à un ind-schéma.
Soit $X$ un ind-schéma sur $k$, d'après \ref{Tdef-alg}, il admet un pro-complexe cotangent $L_{X/k}\in\Pro(\QCoh(X)^{\leq 0})$.
Considérons $z\in X(A)$, on a alors toujours d'après \ref{Tdef-alg}:
\[H_0(z^{*}L_{X/k})\cong (H_{0}(z_{V}^{*}L_{V/k}))_{V},\]
où $V$ parcourt les sous-schémas fermés $V\subset Z$ tels qu'on a une factorisation $z:\Spec(A)\stackrel{z_{V}}{\rightarrow}V\ra X$.
Comme $L_{V/k}$ est à composantes plates et dans $\QCoh(V)^{\leq 0}$, on a:
\[H_{0}(z_{V}^{*}L_{V/k})\cong H_{0}(z_{V}^{*,\heartsuit}L_{V/k})\cong z_{V}^{*,\heartsuit}H_{0}(L_{V/k})\cong z_{V}^{*,\heartsuit}\Omega^{1}_{V/k}.\]
Ainsi, pour tout morphisme $f:\Spec(B)\ra\Spec(A)$,on obtient que 
\[f^{*,\heartsuit}H_0(z^{*}L_{X/k})\cong H_0(f^{*}(z^{*}L_{X/k}))\]
On définit alors $\Omega^{1}_{X/k}\in\Pro(\QCoh(X)^{\heartsuit})$ tel que pour tout $z\in X(A)$, on a:
\[z^{*,\heartsuit}\Omega^{1}_{X/k}=H_{0}(z^{*}L_{X/k}).\]
La construction précédente fournit canoniquement un morphisme de ind-schémas ind-affine:
\[T^{*}X=\ab(\Omega^{1}_{X/k})\ra X.\]
Si $X$ est un $\aleph_0$-ind-schéma sur $k$, alors $\Omega_{X/k}^{1}$ est un pro-ensemble strictement de Mittag-Leffler, i.e. pour tout morphisme $x:\Spec(A)\ra X$, $x^{*,\heartsuit}\Omega_{X/k}^{1}$ est un pro-ensemble strictement de Mittag-Leffler \ref{pro-ml}. En effet, comme c'est un ind-schéma, quitte à choisir un présentation $X\simeq\colim X_{\al}$, on obtient que $x^{*,\heartsuit}\Omega_{X/k}^{1}$ est équivalent à un  système projectif avec des flèches de transition surjectives et comme le système projectif est dénombrable, c'est suffisant pour être strictement de Mittag-Leffler d'après \ref{mldef1}. Ainsi, on déduit que $\Omega^{1}_{X/k}\in\QCoh(X)^{\heartsuit}$.

\bdefi\label{def-ML}
On dit qu'un $\infty$-champ $\cX$ est de Mittag-Leffler s'il admet une théorie de la déformation au sens de \cite[7.1.2]{GRII} et si $\Omega^{1}_{\cX/k}=H_{0}(L_{\cX/k})$ est un pro-ensemble de Mittag-Leffler (cf. \ref{pro-ml}). Dans ce cas, $\Omega^{1}_{\cX/k}$ provient d'un objet de $\QCoh(\cX)^{\heartsuit}$.
\edefi
La discussion précédente montre qu'un $\aleph_0$-ind-schéma est de Mittag-Leffler. Dans la section suivante, on voit aussi que c'est le cas pour certains champs quotients $[X/H]$.

\subsubsection{Pro-complexe cotangent équivariant}
Soit un $\aleph_0$-ind-schéma $X$ sur $k$ ind-placide, formellement lisse avec une action d'un $\aleph_0$-ind-schéma en groupes $H$ ind-placide, formellement lisse.
Il résulte alors de la section précédente que les pro-complexes cotangents $L_{X/k}$ et $L_{H/k}$ sont en fait des objets de $\QCoh$.
On considère l'$\infty$-champ quotient $[X/H]$, d'après \ref{cot-eq}, il admet un complexe cotangent et l'on a une suite fibrante:
\[\pi^{*}L_{[X/H]/k}\ra L_{X/k}\ra L_{X/[X/H]}.\]
Calculons $L_{X/[X/H]}$. Comme $\QCoh$ est un faisceau fpqc, par descente, on a $L_{X/[X/H]}\cong e^{*}L_{H/k}\otimes_{k}\cO_{X}$ où $e$ est la section unité de $H\ra\Spec(k)$. 
De plus, comme $H$ est formellement lisse et ind-placide, d'après \ref{fl-gr}, $\tau_{\geq -1}L_{H/k}\cong\Omega^{1}_{H/k}$.
En appliquant la suite exacte longue, on obtient:
\[0\ra \pi^{*}\Omega^{1}_{[X/H]/k}\ra \Omega^{1}_{X/k}\stackrel{da}{\rightarrow}\Lie(H)^{\vee}\otimes_{k}\cO_{X}\ra 0.\]
Ainsi, comme $\pi^{*}\Omega^{1}_{[X/H]/k}$ est le noyau d'un morphisme de faisceaux quasi-cohérents, il est aussi quasi-cohérent et donc par descente $[X/H]$ est un $\infty$-champ de Mittag-Leffler. 
La construction \eqref{cstr-fib} appliquée à $\Omega_{[X/H]/k}^{1}$ donne un champ cotangent:
\[T^{*}[X/H]\ra [X/H].\]
Décrivons la flèche $da$, pour tout $x\in X$, on a un morphisme $a_x:H\ra X$ donné par $h\mapsto h.x$ qui induit un morphisme:
\[da_x:T_{x}^{*}X\ra \Lie(H)^{\vee}.\]
qui est la fibre en $x$ de l'application $\Omega^{1}_{X/k}\ra \Lie(H)^{\vee}\otimes_k\cO_X$.
On obtient ainsi une application moment:
\[\mu_{H}: T^{*}X\ra \Lie(H)^{\vee}\]
donnée par $(x,\xi)\mapsto da_{x}(\xi)$ ainsi que l'identification :
\begin{equation}
T^{*}[X/H]\times_{[X/H]}X\stackrel{\sim}{\rightarrow}\mu^{-1}_{H}(0).
\label{cotanH}
\end{equation}
qui induit un isomorphisme :
\begin{equation}
[\mu_{H}^{-1}(0)/H]\stackrel{\sim}{\rightarrow}T^{*}[X/H].
\label{cotanH2}
\end{equation}
\subsection{Transversalité}
On reprend les définitions de \cite{Bei}, étendues dans le cas des $\infty$-champs de Mittag-Leffler. Pour le moment il n'est pas nécessaire d'introduire la moindre propriété de lissité. 

\bdefi
Soit un fermé cônique $C\subset T^*\cX$, $B(C)$ l'image de $C$ dans $\cX$; des morphismes de $\infty$-champs de Mittag-Leffler $h:\cU\rightarrow \cX$ et $f:\cX\rightarrow \cZ$.
\begin{enumerate}
	\item 
Soit $u$ (resp. $x$) un point géométrique de $\cU$ (resp. $\cX$), On dit que $h$ (resp. $f$) est $C$-\textsl{transverse} en $u$ (resp. $x$) si pour tout point géométrique $u\in \cU$ (resp $x\in \cX)$, $(\Ker(dh_u)\cap C_{h(u)})\backslash\{0\}=\emptyset$ (resp. $(df_x)^{-1}(C_x)\backslash\{0\}=\emptyset$).
\item
On dit que $h$ (resp. $f$) est $C$-transverse si elle l'est en tout point géométrique de $\cU$ (resp. $\cX$).
\end{enumerate}
\edefi
Soit un morphisme $h:\cU\rightarrow \cX$  entre $\infty$-champs de Mittag-Leffler, $dh:T^{*}\cX\times_{\cX} \cU\ra T^{*}\cU$ la flèche induite.
Soient deux fermés côniques $C\subset T^{*}\cX$, $C'\subset T^{*}\cU$, on définit alors les sous-ensembles côniques:
\[h^{\circ}C=[dh](C)\subset T^{*}\cU, h_{\circ}C'=[p]([dh]^{-1}(C'))\subset T^{*}\cX,\]
avec $p:T^{*}\cX\ra \cX$. Pour rappel, les crochets désignent l'application au niveau des espaces topologiques.
En général, ces sous-ensembles n'ont pas de raison d'être fermés. 
Si de plus, $B(C')$ est propre sur $\cX$, on obtient ainsi un fermé cônique  de $T^{*}\cX$ avec $B(h_{\circ}C')=h(B(C'))$.
On commence par établir des propriétés de changement de base.

\blem\label{Cbc} 
\benumr
	\item 
On considère un carré cartésien entre $\infty$-champs de Mittag-Leffler:
$$\xymatrix{\cX_{V}\ar[d]_{g_{\cV}}\ar[r]^{h_{\cV}}&\cX\ar[d]^{g}\\\cV\ar[r]^{h}&\cZ}.$$
Soit $C\in T^{*}\cX$ un fermé cônique, alors on a $h^{\circ}g_{\circ}C=g_{\cV,\circ}h_{\cV}^{\circ}C$.
\item
Soit un morphisme $g:\cX\rightarrow \cZ$  entre $\infty$-champs de Mittag-Leffler, $C\subset T^*\cX$ un fermé cônique, alors pour tout morphisme $\infty$-champs de Mittag-Leffler $h:\cV\rightarrow \cZ$ $g_{\circ}C$-transverse, $h_{\cV}:\cX_{\cV}=\cX\times_{\cZ}\cV\rightarrow \cX$ est $C$-transverse.
\eenum
\elem
\bpf
(i) Le calcul de $g_{\cV,\circ}h_{\cV}^{\circ}C$ se fait à l'aide du diagramme suivant, où le carré central est cartésien:
$$\xymatrix{T^{*}\cZ\times_{Z}\cX_{\cV}\ar[r]\ar[d]&T^{*}\cZ\times_{\cZ}\cX\ar[d]\ar[r]^-{dg}& T^{*}\cX\\T^{*}\cZ\times_{\cZ}\cV\ar[r]\ar[d]^{dh}&T^{*}\cZ\\T^{*}\cV}$$
par image inverse et image directe. 
Par ailleurs, on a également un autre diagramme commutatif avec un carré central cartésien qui permet de calculer $h^{\circ}g_{\circ}C$:
\begin{equation}
\xymatrix{T^{*}\cZ\times_{\cZ}\cX_{\cV}\ar[r]^{dg}\ar[d]_{dh}&T^{*}\cX\times_{\cX}\cX_{\cV}\ar[d]^{dh_{\cV}}\ar[r]^-{p}& T^{*}\cX\\T^{*}\cV\times_{\cV}\cX_{\cV}\ar[d]\ar[r]^-{dg_{\cV}}&T^{*}\cX_{\cV}\\T^{*}\cV}.
\label{cart1}
\end{equation}
L'énoncé vient alors du fait que l'on a l'égalité ensembliste $dh(dg^{-1}(p^{-1}(C)))=dg_{\cV}^{-1}(dh_{\cV}(p^{-1}(C)))$.

(ii) Il résulte du diagramme \eqref{cart1} que l'on a une application surjective induite par $dg$:
\[\Ker(dh)\twoheadrightarrow\Ker(dh_{\cV})\]
et donc une surjection:
\[\Ker(dh)\cap dg^{-1}(p^{-1}(C))\twoheadrightarrow\Ker(dh_{\cV})\cap p^{-1}(C).\]
et comme $h$ est $g_{\circ}C$-transverse, on en déduit que $dh_{\cV}$ est $C$-transverse.
\epf
\blem\label{Cbc2}
Soient $p:\cX\rightarrow \cY$, $f:\cY\rightarrow \cZ$ entre $\infty$-champs de Mittag-Leffler, $C\in T^{*}X$ un fermé cônique, alors il a équivalence entre:
\benumr
	\item 
	$f\circ p$ est $C$-transverse.
	\item
	$f$ est $p_{\circ}C$-transverse.
\eenum
\elem
\bpf
 L'argument est le même que \cite[Lem. 3.8]{Sai}, une fois que tous les objets sont bien définis. La propreté dans loc.cit. n'est là que pour assurer que $p_{\circ}C$ est fermé et la lissité de $\cX,\cY, \cZ$ est inutile (seulement là parce que dans loc.cit., on ne considère la $C$-transversalité que pour des schémas lisses).
\epf
\bprop\label{C1}
Soient $Y$ un $k$-schéma lisse de type fini, $C\subset T^*Y$ un fermé cônique, $Z\in\Sch_k$ placide. 
Soit un morphisme $h:Z\rightarrow Y$ $C$-transverse, alors il existe $Z_{\al}$ de type fini tel que $h$ se factorise étale localement sur $Z$, en $Z\stackrel{p_{\al}}{\rightarrow} Z_{\al}\stackrel{h_{\al}}{\rightarrow} Y$ où $h_{\al}$ est $C$-transverse  et $p_{\al}$ fortement pro-lisse. Si de plus, $Z$ est formellement lisse, on peut supposer $Z_{\al}$ lisse.
	\eprop
\bpf
(i) Comme $Z$ est placide quasi-compact, soit $Z'\ra Z$ étale tel que $Z'\in\Sch_k$ admet une présentation placide.
La composée $Z'\ra Z\ra Y$ reste $C$-transverse et il suffit de montrer qu'elle admet une factorisation telle que dans l'énoncé. On suppose donc dans la suite que $Z$ admet une présentation placide.
On l'écrit alors $Z\simeq\varprojlim Z_{\al}$ avec $Z_{\al}$ de type fini et les flèches de transition affines, lisses.
Comme $Y$ est affine de type fini, il existe $\al$ tel que $h$ se factorise en $Z\rightarrow Z_\al\rightarrow Y$. Soient $z\in Z$, $z_{\al}=p(z)$ et $y=h(z)$, alors la flèche $dh_z:T^{*}_{y}Y\rightarrow T_{z}^{*}Z$ se factorise par $T^{*}_{z_{\al}}Z_{\al}$. Pour conclure, il suffit de remarquer que d'après \ref{plac2}, on a un morphisme universellement injectif de $\cO_{Z}$-modules:
\[p_{\al}^*\Omega^1_{Z_{\al}/k}\rightarrow\Omega^1_{Z/k}.\]
Comme la condition de $C$-transversalité est ouverte entre $k$-schémas de type fini d'après \cite[Lem. 1.2]{Bei} \footnote{ dans loc.cit. les schémas sont lisses, mais la même preuve vaut pour des $k$-schémas de type fini.}, on peut donc, quitte à remplacer $Z_{\al}$ par un voisinage ouvert de $p_{\al}(Z)$, supposer que $Z_{\al}\rightarrow Y$ est bien $C$-transverse.
Si de plus, $Z$ est formellement lisse, alors d'après \ref{plac}, on peut supposer que tous les $Z_{\al}$ sont lisses.
\epf
Enfin, il résulte de \ref{plac2} que pour un morphisme de schémas qcqs $h:X\ra Y$ fortement pro-lisse ou formellement lisse, on a pour tout fermé cônique $C\subset T^{*}Y$:
\begin{equation}
h^{\circ}C=C\times_{Y}X.
\label{Clis}
\end{equation}
Ainsi, si de plus $B(C)$ est fermé, il en est de même de $B(h^{\circ}C)$. 
\blem\label{open}
Soit $X$ un ind-schéma de Mittag-Leffler, $C\subset T^{*}X$ un fermé cônique et un morphisme de ind-schémas $f:X\ra Y$ avec $Y$ un $k$-schéma de type fini, alors la propriété d'être $C$-transverse pour $f$ est ouverte.
\elem
\bpf
C'est la même preuve que \cite[Lem. 1.2]{Bei}, le lieu de $C$-transversalité corrrespond au complémentaire de l'image du projectivisé $\bP((df)^{-1}(C))\subset\bP(T^{*}Y_{X})$ dans $X$ qui est bien fermé puisque comme $Y$ est de type fini, on a $\bP(T^{*}Y_{X})\ra X$ qui est un morphisme projectif et que $df^{-1}(C)$ reste un fermé cônique de $T^{*}Y_{X}$.
\epf
\subsubsection{Stabilité pour les fermés côniques}
Dans cette section, on cherche à déterminer un contexte plus large que \eqref{Clis} qui nous assure que $h^{\circ}C$ reste un fermé cônique de base fermée si $C$ l'est.
Typiquement, on va voir que c'est le cas pour un morphisme $C$-transverse avec quelques hypothèses supplémentaires sur les schémas.
On commence par introduire la définition suivante:
\bdefi\label{defquot}
Soit $\cX\in\St_k$, on dit que c'est un champ quotient s'il existe un isomorphisme $\cX\simeq[X/H]$ où $X$ est un $\aleph_0$-ind-schéma raisonnable (cf. \ref{ratio}) formellement lisse et $H$ un $\aleph_0$-ind-schéma en groupes raisonnable,  ind-affine et formellement lisse.
On appelle $X$ un atlas.
\edefi
\brem
La condition d'être un $\aleph_0$-ind-schéma assure automatiquement d'être Mittag-Leffler et étant donné que dans la suite seuls des sytèmes dénombrables apparaissent, cette hypothèse simplificatrice est commode.
On impose que $H$ soit ind-affine pour que la flèche $X\ra[X/H]$ soit ind-schématique.
\erem

\blem\label{modI}
On considère une algèbre  graduée positivement $R=\bigoplus\limits_{n\geq0} R_n$, soit $P=(P_{\al})$ une $R$-pro-algèbre graduée et l'idéal $I=\bigoplus\limits_{n\ssup 0} R_{n}$. On suppose que $P/I$ est une $R/I$-pro-algèbre entière (i.e. pour tout $\al$, $P_{\al}/I$ est une $R/I$-algèbre entière), alors $P$ est une $R$-pro-algèbre entière.

De plus, si $P$ est de Mittag-Leffler et pour tout $\al$ et $\beta\geq\al$, $P_{\beta}\ra P_{\al}$ a un noyau de type fini et $P_{\beta}/I\rightarrow P_{\al}/I$ est à noyau nilpotent, $P_{\beta}\rightarrow P_{\al}$ est aussi à noyau nilpotent.
\elem

\bpf
On commence par considérer le cas où $P$ est une algèbre. Soit $x\in P$, on définit le $R$-module gradué $M=R[x]=\bigoplus_{n\in\NN} M_{n}$ avec la graduation $M_{n}=\bigoplus\limits_{i+j=n} R_ix^j$. Comme $P/I$ est entière, on en déduit que pour tout $\al$, $M/IM$ est un $R/I$-module de type fini, donc $M$ est un $R$-module de type fini d'après \cite[Lem. 3.2]{Sai} et $x$ est entier sur $R$, ce qu'on voulait.

Passons au cas général, en appliquant le cas précédent aux triplets $(A,R,P_{\al})$ on obtient que $P_{\al}$ est une $R$-algèbre entière. Supposons de plus que $P$ est  Mittag-Leffler, fixons $\al$ et $\beta\geq\al$, soit $K_{\beta\al}:=\Ker(P_{\beta}\rightarrow P_{\al})$.
Par hypothèse, $K_{\beta\al}/I$ est nilpotent d'ordre un certain $r$, comme $K_{\beta\al}$ est le noyau d'un morphisme d'algèbres graduées, il est également gradué, ainsi que $K_{\beta\al}^r$ qui est aussi de type fini. Soient $x_1,\dots, x_m$ ses générateurs, sans restreindre la généralité, on peut les supposer homogènes et alors le $R$-module $N$ engendré par les $x_i$ est un sous-module gradué de type fini avec $N/I=0$ donc $N=0$ et $K_{\beta\al}^r=0$.
\epf
On utilise le lemme précédent sous la forme du corollaire suivant:
\bcor\label{cormod}
Soient $T$, $T'$ des ind-schémas ind-affines munis chacun de présentations $\bG_m$-équivariantes avec $T'$ raisonnable, soit $h:T'\rightarrow T$ un morphisme de ind-schémas $\bG_m$-équivariant, on suppose que la flèche au niveau des points fixes $h^{\bG_m}:T'\times_{T}T^{\bG_m}\rightarrow T^{\bG_m}$ est  ind-entière avec des flèches de transition nilpotentes , alors il en est de même de $h$. De plus, $h$ est un morphisme fermé.
\ecor
\brem
Comme $T$ admet une présentation équivariante, on a immédiatement par réduction au cas des schémas que $T^{\bG_m}$ est un sous-ind-schéma fermé.
\erem
\bpf
Comme $T$ admet une présentation $\bG_m$-équivariante et est ind-affine, il suffit de vérifier l'assertion après restriction à tout sous-schéma fermé $\bG_m$-équivariant. On se ramène au cas où $T$ est un schéma affine, et comme $T'$ est raisonnable, on applique alors \ref{modI}, pour obtenir que $T'\simeq \colim T'_{\al}$ où les $h_{\al}:T'_{\al}\rightarrow T$ sont entières, donc fermées et les flèches  de transition sont des immersions fermées nilpotentes.
On obtient donc qu'au niveau des espaces topologiques, on a pour tout $\al$ et tout fermé $F\subset T$, $h(F)=h_{\al}(F_{\al})$, qui est bien fermé, ce qui conclut.
\epf
 Si l'on ajoute des hypothèses de transversalité, on a la proposition suivante pour l'image inverse:
\bprop\label{Cind}
Soit $\cX$ est un champ quotient, $C\subset T^{*}\cX$ un fermé cônique tel que $B(C)$ est fermé dans $\cX$.
Soit un morphisme  $C$-transverse entre $\infty$-champs $h:\cU\rightarrow \cX$ , avec $h$ est $\aleph_0$-ind-schématique raisonnable et $\cU$ formellement lisse,  alors $dh_{\cU}:C_{\cU}:=C\times_{\cX}\cU\rightarrow T^{*}\cU$ est d'image fermée. Ainsi, $h^{\circ}C$ est fermé et on a $B(h^{\circ}C)=h^{-1}(B(C))$.
\eprop
\bpf
Comme $\cX$ est un champ quotient, on l'écrit $\cX\simeq[X/H]$ comme dans \ref{defquot} et on forme $U=\cU\times_{\cX}X$.
Alors comme $h$ est $\aleph_0$-ind-schématique raisonnable, $U$ est un $\aleph_0$-ind-schéma raisonnable et on a $[U/H]\cong\cU$ d'après \ref{loopeq}. Comme $\cU$ est formellement lisse ainsi que $H$, $U$ est formellement lisse.
On commence par montrer l'assertion dans le cas $H=\{1\}$.
On écrit $U=\colim U_{\al}$, il suffit de vérifier la propriété après restriction à chaque $U_{\al}$ et ensuite, on peut supposer que $U_{\al}=\Spec(A)$ est affine.
Le morphisme $T^*U\rightarrow U$ est ind-affine et  $T^{*}U_{\vert U_{\al}}:=T^{*}U\times_{U}U_{\al}=\colim Z_{\al}$ où $Z_{\al}$ admet une action de $\bG_m$ ainsi que $C_{U_{\al}}$. Comme $h$ est ind-schématique raisonnable, $C_{U_{\al}}$ est aussi un ind-schéma raisonnable. Enfin, on a $(T^{*}U_{\vert U_{\al}})^{\bG_m}=U_{\al}$ et la flèche $\sigma:U_{\al}=(T^{*}U_{\vert U_{\al}})^{\bG_m}\rightarrow T^{*}U_{\vert U_{\al}}$ s'identifie à la section nulle.
D'après \ref{cormod}, il suffit donc de montrer que $dh$ est entière après restriction à la section nulle et la condition de $C$-transversalité donne que la flèche s'identifie à l'inclusion du fermé:
\[(\sigma^{*}C_{U_{\al}})_{red}=h^{-1}(B(C))\subset U_{\al},\]
ce qui conclut.

Passons au cas général, on a le carré cartésien suivant:
$$\xymatrix{U\ar[d]_{\ti{\pi}}\ar[r]^{\ti{h}}&X\ar[d]^{\pi}\\\cU\ar[r]^{h}&\cX}$$
On a $\cU\simeq [U/H]$ et d'après \eqref{cotanH} et \eqref{cotanH2}, on a également que $[\mu^{-1}(0)/H]\simeq T^{*}\cU$ avec $\mu^{-1}(0)=T^{*}\cU\times_{\cU}U$ qui est fermé dans $T^{*}U$.
En particulier, comme $T^{*}\cU$ est muni de la topologie quotient, $h^{\circ}C$ est fermé si et seulement si $h^{\circ}C\times_{\cU}U$ l'est et si et seulement si son image dans $T^{*}U$ l'est aussi, laquelle s'identifie  à  $\ti{\pi}^{\circ}(h^{\circ}C)$.
Montrons donc que $\ti{\pi}^{\circ}(h^{\circ}C)$ est fermé.
On a $\ti{\pi}^{\circ}(h^{\circ}C)=\ti{h}^{\circ}(\pi^{\circ}C)$. De plus, d'après \ref{flquot}, $\pi$ est formellement lisse et il résulte à nouveau des descriptions \eqref{cotanH} et \eqref{cotanH2}, que $\pi_{\circ}(\pi^{\circ}C)=C$ et que $\pi^{\circ}C$ est fermé avec $B(\pi^{\circ}C)=\pi^{-1}(B(C))$. En particulier, d'après \ref{Cbc}(2), $\ti{h}$ est $\pi^{\circ}C$ transverse et le cas des ind-schémas donne que $\ti{h}^{\circ}(\pi^{\circ}C)$ est fermé comme souhaité et l'égalité sur les bases suit.
\epf
\section{Microsupports}
Dans cette section, on introduit la classe des champs microsupportables et on montre des énoncés de fonctorialité par pullback et poussé-en-avant pour des morphismes ind-propres formellement lisses. L'énoncé particulièrement est \ref{fonct-clo} puisqu'il permet de passer des schémas aux ind-schémas. On établit ensuite des énoncés d'existence de supports singuliers pour des schémas placides \ref{SSex0} et on montre que la théorie développée dans ce travail généralise les travaux de Beilinson-Saito \ref{Beisupp} et \ref{fonct1}.

\subsection{Champs microsupportables}\label{gen}
On commence par introduire une classe de $\infty$-champs suffisamment large pour laquelle on peut parler de microsupports.
\bdefi\label{defmicro}
Soit $\cX\in\St_k$, on dit qu'il est microsupportable si:
\benumr
\item
	$\cX\simeq\coprod_{\al}\cX_{\al}$ où chaque $\cX_{\al}$ est un champ quotient au sens de \ref{defquot}.
	\item
	Pour tout $\al$, il existe une présentation $\cX_{\al}=[X_{\al}/H_{\al}]$ telle que dans \ref{defquot} où l'on suppose de plus que $X_{\al}$ et $H_{\al}$ sont ind-placides.
	\item
	$\cX$ est formellement lisse.
	\eenum
		\edefi
\brems
\remi\label{defmic1-1}
La condition (i) assure que chaque composante connexe est de Mittag-Leffler et que l'on peut définir une notion de $C$-transversalité avec des énoncés de fonctorialité par pullbacks pour les fermés côniques.
Un quotient de deux objets ind-placides satisfait le recollement de faisceaux et une colimite de tels objets avec des morphismes de transition qui sont des immersions ouvertes quasi-compactes satisfait aussi le recollement d'après \ref{glustk}. On dispose donc d'une notion de locale acyclicité.
La dernière condition de lissité est nécessaire pour pouvoir utiliser les énoncés sur le support singulier de Beilinson \cite{Bei}.
\remi\label{defmic1-2}
Soit un morphisme d'$\infty$-champs $h:\cY\ra\cX$ quasi-compact schématique, placide et formellement lisse avec $\cX$ microsupportable, alors $\cY$ l'est aussi. On a clairement que $\cY$ est formellement lisse. En passant aux composantes, on peut supposer que $\cX$ est de plus un champ quotient et si $\cX\simeq[X/H]$, on a $\cY\simeq [Y/H]$ avec $Y\simeq \cY\times_{\cX}X$ qui est bien ind-placide et $\aleph_0$, comme $X$ l'est et $h$ placide. 
\remi\label{defmic1-3}
Si $\cY$ est un champ quotient microsupportable avec une présentation $\cY\simeq[Y/H]$ telle que dans \ref{defquot}, d'après \ref{flquot}, on a immédiatement que $Y$ est formellement lisse.
\remi
Dans la suite, on parle de quotients microsupportables pour \og champs quotients microsupportables\fg.
\erems
Le lemme suivant caractérise les ind-schémas qui vont avoir une théorie du microsupport.
\blem\label{eq-indsch}
Soit un $\aleph_0$-ind-schéma $X$ raisonnable, alors il est microsupportable si et seulement si il est ind-placide formellement lisse.
\elem
\bpf
S'il est ind-placide formellement lisse alors on peut prendre $H=\{1\}$ dans la description comme quotient. Réciproquement, s'il est microsupportable alors il est déjà formellement lisse et on considère une décomposition $X=\coprod_{a\in A}[\ti{X}_{a}/H_{a}]$. Comme $X$ est un $\aleph_0$-ind-schéma, cela force $A$ à être dénombrable et pour tout $a$, $[\ti{X}_{a}/H_{a}]$ est un ind-schéma. Il suffit donc de montrer le cas où $X=[\ti{X}/H]$. Il reste à voir que $X$ est ind-placide; comme il est raisonnable, on se ramène alors au cas où $X$ est un schéma qcqs et alors il existe $X'\ra X$ étale quasi-compact tel que $\ti{X}\times_{X}X'\simeq X'\times H$ et comme le membre de gauche est ind-placide, $X'$ l'est aussi et donc également $X$, comme souhaité.

\epf
\bdefi\label{lqc}
\'{E}tant donné un morphisme d'$\infty$-champs $f:\cX\ra\cY$, on dit qu'il est localement quasi-compact (lqc) si:
\benumr
\remi
$f$ est $\aleph_0$-ind-schématique raisonnable.
\remi
$f$ est schématique, i.e. pour tout schéma affine $\Spec(A)\ra\cY$ le changement de base $\cX_{A}$ est représentable par un schéma (mais non-nécessairement quasi-compact).
\eenum
\edefi
Les morphismes lqc sont stables par composition et changement de base. Un morphisme quasi-compact schématique d'$\infty$-champs quotients est évidemment lqc. Un exemple plus intéressant et typique de notre situation est le suivant:
\blem\label{lqc-disc}
Soit $H$ un schéma en groupes qui agit sur un $\infty$-champ quotient $\cX$, on suppose que:
\benumr
\item
On a une décomposition $H\simeq H^{0}\times\La$ où $\La$ est un groupe discret dénombrable.
\item
$H^{0}$ est un schéma en groupes affine placidement présenté.
\eenum
 alors $\cX\ra[\cX/H]$ est lqc.
\elem
\bpf
On peut voir $H$ comme un $\aleph_0$-ind-schéma ind-affine raisonnable sur $k$, de telle sorte que d'après \ref{exquot}, le morphisme $\cX\ra[\cX/H]$ est $\aleph_0$-ind-schématique raisonnable.
De plus, pour tout schéma affine $\Spec(A)\ra\cX$, si $E_{A}=\cX\times_{[\cX/H]}\Spec(A)$ est étale-localement isomorphe à $\Spec(A)\times H$
et comme il est représentable par un schéma, c'est un faisceau fppf avec une diagonale représentable, en particulier, on obtient que $E_{A}$ est un espace algébrique et un ind-schéma donc c'est un schéma (la question est Zariski-locale, donc on peut supposer que $E_{A}$ est quasi-compact et un atlas de $E_{A}$ est un schéma dont l'image se factorise par un sous-schéma fermé de $E_{A}$).
\epf

 \subsection{Microsupports et premières fonctorialités}
\bdefi\label{pairetest}
Soient $\cX$ un $\infty$-champ microsupportable et $C\subset T^{*}\cX$ un fermé cônique tel que $B(C)$ est fermé.
\begin{enumerate}
	\item  Une \textsl{paire test} $(h,f)$ est une correspondance  $\cX\stackrel{h}{\leftarrow}\cU\stackrel{f}{\rightarrow} Y$ où $\cU$ microsupportable, $h$ lqc et $Y$ lisse de type fini.
On dit que la paire $(h,f)$ est $C$-\textsl{transverse} si $h$ est $C$-transverse et $f$ est $h^{\circ}C$-transverse.
\item
Soit $K\in \cD(X)$, une paire test $(h,f)$ est $K$-\textsl{acyclique} si $h^!K$ est localement acyclique relativement à $f$.
On dit que $K$ est $C$-\textsl{microsupporté} si pour toute paire $C$-transverse $(h,f)$, $(h,f)$ est $K$-\textsl{acyclique}.
\item
On définit alors $\cC(K)$ l'ensemble des microsupports de $K$. 
S'il existe, on note $SS(K)$, le support singulier de $K$, le minimum pour l'inclusion de $\cC(K)$.
On considère aussi l'ensemble $\ov{\cC}(K)$ défini de la même manière que $\cC(K)$ mais où l'on n'impose pas que $B(C)$ soit fermé. 
\end{enumerate}
\edefi
\brems
\remi\label{rem-ptest1}
Il résulte de \ref{Cind}, que si $h$ est lqc et $C$-transverse pour un fermé cônique $C$ de base fermé, alors $h^{\circ}C$ reste un fermé cônique de base fermé.
\remi\label{rem-ptest2}
 On remarque que si $X$ est lisse de type fini, d'après \cite[sect. 1.1]{Bei}, $B(C)$ est toujours fermé.
\erems

\blem\label{fonct}
Soit un $\infty$-champ microsupportable $\cX$ et $K\in\cD(\cX)$. Soit $C\in\cC(K)$ et $h:\cU\rightarrow \cX$ $C$-transverse et lqc, avec $\cU$  microsupportable, alors $h^{\circ}C\in\cC(h^!K)$.
\elem
\bpf
Soit $\cV$ un $\infty$-champ microsupportable et $r:\cV\rightarrow \cU$ $h^{\circ}C$-transverse et lqc, alors la composée $\cV\rightarrow \cU\rightarrow \cX$ est $C$-transverse.
En particulier, une paire test $(g,r):W\leftarrow \cV\rightarrow \cU$ $h^{\circ}C$-transverse fournit par composition une paire test $(g,r\circ h):W\leftarrow \cV\rightarrow \cX$ $C$-transverse, donc la paire $(g,r\circ h)$ est $K$-acyclique et d'après \ref{rem-ptest1} $h^{\circ}C\in\cC(h^!K)$.
\epf

\bprop\label{propfonct}
Soit un morphisme d'$\infty$-champs microsupportables $h:\cX\ra\cY$  ind-fp-propre et formellement lisse. Soient $K\in\cD(\cX)$ et $C\in\cC(K)$,
 alors $\overline{h_{\circ}C}\in\ov{\cC}(h_!K)$.
\eprop

\bpf
Soit une paire test $(g,r):W\leftarrow \cV\rightarrow \cY$ $h_{\circ}C$-transverse avec $r$ lqc, $\cV$ microsupportable et $W$ lisse de type fini, il suffit de voir que $r^{!}(h_!K)$ est $g$-acyclique. Ensuite, en prenant l'adhérence, on obtient immédiatement que $\overline{h_{\circ}C}\in\ov{\cC}(h_!K)$. En écrivant $\cV\simeq\coprod\cV_{\al}$ avec $\cV_{\al}$ un champ quotient microsupportable d'après \ref{pulLis}(i), on se ramène au cas où $\cV$ est un quotient microsupportable, ce que l'on suppose dans la suite.
De même, on peut supposer que $\cY$ est un quotient microsupportable.
On forme le carré cartésien:
\begin{equation}
\xymatrix{\cX_{\cV}\ar[d]_-{\ti{h}}\ar[r]^{\ti{r}}&\cX\ar[d]^{h}\\\cV\ar[r]^{r}&\cY}
\label{eq-pr}
\end{equation}
\blem\label{stabmic}
Le champ $\cX_{\cV}$ est  microsupportable.
\elem
\bpf
Comme par changement de base $\cX_{\cV}\ra\cV$ est formellement lisse et que $\cV$ est formellement lisse, $\cX_{\cV}$ l'est.
On écrit alors $\cV\simeq [V/H_1]$ et $\cX\simeq [X/H_2]$ comme dans \ref{defmicro}.
On forme le diagramme suivant où tous les carrés sont cartésiens:
$$\xymatrix{X_{V}\ar[d]\ar[r]&X_{\cV}\ar[r]\ar[d]&X\ar[d]\\\cX_{V}\ar[r]\ar[d]&\cX_{\cV}\ar[d]\ar[r]&\cX=[X/H_2]\ar[d]\\V\ar[r]&\cV\simeq[V/H_1]\ar[r]^-{h}&\cY}$$
La flèche $X_{\cV}\ra X$ est lqc, donc $\aleph_0$-ind-schématique et $X_{V}\ra X_{\cV}$ est $\aleph_0$-ind-schématique d'après \ref{exquot}.
Ainsi, $X_{V}$ est un $\aleph_0$-ind-schéma et l'on a $[X_{V}/(H_1\times_k H_2)]\simeq \cX_{\cV}$ avec $H_1\times_k H_2$ ind-affine, ind-placide et formellement lisse. Il reste à voir que chaque $X_{V}$ est ind-placide.
Or la flèche $X_{V}\ra V$ est la composée d'un $H_2$-torseur et d'un morphisme ind-fp-propre, qui sont tous les deux des morphismes ind-placides, ce qui conclut.
\epf
 Comme $r$ est $h_{\circ}C$-transverse, d'après \ref{Cbc}.(ii), $\ti{r}$ est $C$-transverse et par changement de base, il est lqc. Comme $g$ est $r^{\circ}h_{\circ}C$-transverse, d'après \ref{Cbc}.(ii), $g$ est $\ti{h}_{\circ}\ti{r}^{\circ}C$-transverse et donc d'après \ref{Cbc2}, on obtient que la paire $(g\circ\ti{h},\ti{r})$ est $C$-transverse.
Ainsi, $\ti{r}^{!}K$ est $g\circ\ti{h}$-acyclique et donc d'après \ref{acprop}, $\ti{h}_{!}\ti{r}^{!}K$ est $g$-acyclique et on conclut par le changement de base propre \ref{Base}.
\epf

\blem\label{equivquot}
Soit un $\aleph_0$-ind-schéma $X$ microsupportable avec une action d'un schéma en groupes qcqs fortement pro-lisse $H$. Soient $K\in\cD(X)$ $H$-équivariant et $C\in\cC(K)$ tel que $C\subset\mu_{H}^{-1}(0)$ et est $H$-équivariant, alors $p_{\circ}C\in\cC(\ov{K})$ où $\ov{K}$ est l'unique faisceau tel que $K\cong p^{!}\ov{K}$ avec $p:X\ra[X/H]$.
\elem
\bpf
Tout d'abord, d'après \eqref{cotanH2} et comme $C\subset\mu_{H}^{-1}(0)$ et est $H$-équivariant, on a $p_{\circ}C=[C/H]\subset T^{*}[X/H]$. En particulier, c'est un fermé cônique de base fermée. Montrons que $p_{\circ}C\in\cC(\ov{K})$.
Soit une paire test $(g,r):W\leftarrow \cV\ra [X/H]$ $p_{\circ}C$-transverse avec $\cV$ microsupportable, $W$ lisse de type fini et $r$ lqc, alors d'après \eqref{cotanH2} et \ref{Cbc2}, la paire  $W\stackrel{g\circ\ti{p}}{\leftarrow} X_{\cV}=X\times_{[X/H]}\cV\stackrel{\ti{r}}{\rightarrow} X$ est $C$-transverse et d'après \ref{defmic1-2}, $X_{\cV}$ est microsupportable car $\cV$ l'est. En particulier, $\ti{r}^!K$ est $g\circ\ti{p}$-acyclique.
De plus, la flèche $X_{\cV}\ra\cV$ est un $H$-torseur localement pour la topologie étale, donc d'après \ref{pulLis}, $r^{!}\ov{K}$ est $g$-acyclique si $\ti{p}^{!}(r^{!}\ov{K})\cong \ti{r}^!K$ est $g\circ\ti{p}$-acyclique, ce qui conclut.
\epf
\subsection{Fonctorialité pour une immersion fermée}
L'énoncé principal est \ref{fonct-clo}. Il nous permet d'avoir un énoncé de fonctorialité pour les microsupports pour une immersion fermée et en particulier de pouvoir passer des schémas aux ind-schémas.

Pour un ind-schéma $T$ et un point $x\in T$, on rappelle que l'on note le localisé en $x$:
\[T_{x}=\colim_{S\subset T}\Spec(\cO_{S,x}),\]
où $S$ parcourt les sous-schémas fermés de $T$. Le morphisme $T_{x}\ra T$ est affine pro-Zariski.
\bprop\label{critinflis}
Soient $Z, X$ des  $k$-$\aleph_0$-ind-schémas ind-placides formellement lisses,  $i:Z\hra X$ une immersion fermée de présentation finie. Soit un fermé cônique $C\subset T^{*}Z$ et un morphisme de $\aleph_0$-ind-schémas ind-placides formellement lisses $h:V\ra Z$. On suppose que $h$ est $i_{\circ}C$-transverse, alors pour tout $x\in\ti{h}^{-1}(B(C))$,  $Z_{V,x}=(V\times_{X}Z)_{x}$ est un $k$-ind-schéma formellement lisse  avec $\ti{h}:Z_{V}\ra Z$.
\eprop
\brem
Lorsque l'on travaille avec des ind-schémas, il n'est pas clair que la formelle lissité de $Z_{V,x}$ implique la formelle lissité dans un voisinage.
\erem
\bpf
Formons le carré cartésien:
$$\xymatrix{Z_{V}\ar[d]_{\ti{i}}\ar[r]^{\ti{h}}&Z\ar[d]^{i}\\V\ar[r]^{h}&X}.$$
Par changement de base $Z_{V}\hra V$ est une immersion fermée de présentation finie et $Z_V$ est un $\aleph_0$-ind-schéma ind-placide, comme $V$ en est un.
D'après \ref{exa-clo}, on a un diagramme commutatif avec des lignes exactes et où la ligne du haut est constitué de faisceaux localement pro-projectifs:
$$\xymatrix{0\ar[r]&\ti{h}^{*}C_{Z/X}\ar[r]\ar[d]&(h\circ\ti{i})^*\Omega^{1}_{X/k}\ar[r]^-{di}\ar[d]_{dh}&\ti{h}^*\Omega^{1}_{Z/k}\ar[r]\ar[d]&0\\&
C_{Z_{V}/V}\ar[r]^{\ti{i}^{*}}&\ti{i}^{*}\Omega^{1}_{V/k}\ar[r]^{d\ti{i}}&\Omega^{1}_{Z_{V}/k}\ar[r]& 0}$$
 D'après \ref{exa-clo}, $C_{Z/X}$ est pro-projectif de type fini .
On a aussi une surjection d'après \ref{cone-surj2} :
\begin{equation}
\ti{h}^{*} C_{Z/X}\dhra C_{Z_{V}/V}.
\label{c-surj3}
\end{equation}
Soit alors $x\in \ti{h}^{-1}(B(C))$ et $\phi: S\hra Z_{V}$ un sous-schéma fermé de présentation finie, comme $h$ est $i_{\circ}C$-transverse, on a donc en particulier que :
\[\Ker dh_{\ti{i}(x)}\cap\Ker(di_{\ti{h}(x)})\subset\Ker dh_{\ti{i}(x)}\cap(i_{\circ}C)_{\tilde{h}(x)}=\{0\}.\]
Ainsi, sur le corps résiduel, on a une injection :
\[0\ra\ti{h}^{*} C_{Z/X}\otimes_{\cO_{S}} k(x)\ra \ti{i}^{*}\Omega^{1}_{V/k}\otimes_{\cO_{S}} k(x).\]
et d'après \eqref{c-surj3}, on a aussi $\ti{h}^{*} C_{Z/X}\otimes_{\cO_{S}} k(x)=C_{Z_{V}/V}\otimes_{\cO_{Z_{V}}} k(x)$.
Par construction \ref{cone-def}, $C_{Z_{V}/V}$ est pro-(de présentation finie), donc d'après \ref{Nak1}(i), on obtient donc un isomorphisme de pro-$\cO_{S,x}$-modules:
\[\phi^{*}\ti{h}^{*} C_{Z/X,x}\cong \phi^{*}C_{Z_{V}/V,x}.\]
De plus, d'après \eqref{cone-surj}, $C_{Z/X}$ est de Mittag-Leffler et $\aleph_{0}$, donc d'après \ref{aleph-tate}, il est de Tate et le dual d'un module projectif. Comme $V$ est un $\aleph_0$-ind-schéma ind-placide formellement lisse, d'après \ref{ind-Tate}, $\Omega^{1}_{V/k}$ est de Tate.
On déduit alors de \ref{Nak1}(ii) une suite exacte:
\begin{equation}
0\ra\phi^{*}C_{Z_{V}/V,x}\ra\phi^{*}\ti{i}^{*}\Omega^{1}_{V/k,x}\ra\phi^{*}\Omega^{1}_{Z_{V}/k,x}\ra 0
\label{exac}
\end{equation}
  et que $\phi^{*}\Omega^{1}_{Z_{V}/k,x}$ est topologiquement plat sur $\cO_{S,x}$.
D'après \ref{fond-tate}, comme $Z_{V}\hra V$ est une immersion fermée de présentation finie,  $\Omega^{1}_{Z_{V}/k}$ est pro-(2-presque projectif).
On déduit donc de \ref{pro-pl3}, que  $\phi^{*}\Omega^{1}_{Z_{V}/k,x}$ est pro-projectif. Comme il est par construction $\aleph_0$ et Mittag-Leffler avec des morphismes de transition qui ont des noyaux de type fini \footnote{Cela vient du fait que le ind-schéma est raisonnable.}, il est donc de Tate par \ref{car-mlt}.

De même, il résulte de \eqref{exac}, du fait que $V$ est formellement lisse  et de \ref{fl-gr} que $H_{1}(\phi^*L_{Z_{V,x}/k})=0$, de telle sorte que $Z_{V,x}$ est formellement lisse à nouveau par \ref{fl-gr}.
\epf
\bdefi\label{borne}
Dans la suite, pour un ind-schéma raisonnable $T$ et $K\in\cD_{c}(T)$, il est dit borné s'il est supporté par un sous-schéma fermé de présentation finie $Z\hra T$.
\edefi
\bthm\label{fonct-clo}
On considère une immersion fermée de présentation finie de $k$-$\aleph_0$-ind-schémas ind-placides formellement lisses $i:Z\hra X$.
Soit $K\in\cD_c(Z)$, borné, $C\in\cC(K)$ tel que $\supp(K)\subset B(C)$, alors $i_{\circ}C\in\cC(i_*K)$.
\ethm
\bpf
Soit $(h,f):Y\leftarrow V\ra X$ une paire $i_{\circ}C$-transverse avec $V=\coprod V_{\al}$ avec $V_{\al}$ des ind-schémas microsupportables, $Y$ lisse de type fini et $h$ lqc. En passant aux composantes, on peut supposer que $V$ est un ind-schéma microsupportable.

Soit $\tilde{h}:Z_{V}=V\times_{X}Z\ra Z$, comme $h$ est lqc donc schématique, $Z_{V}$ est un $k$-schéma.
Comme $h$ est $i_{\circ}C$-transverse,  d'après \ref{critinflis}, pour tout $x\in \tilde{h}^{-1}(\supp(K))$, $Z_{V,x}$ est formellement lisse  avec les notations de \ref{critinflis} et donc $f\circ \ti{i}\circ\phi_{x}$ est $\ti{h}^!K_{x}$-acyclique avec $\phi_x: Z_{V,x}\ra Z_{V}$. Comme $K$ est à support schématique, d'après \ref{check-ac}, on en déduit que $f\circ \ti{i}$ est $\ti{h}^!K$-acyclique \footnote{la condition de locale acyclicité sur les localisations Zariski induit la locale acyclicité sur les localisations étales par \ref{pulLis}.}.
Ainsi, on obtient que $f$ est $\ti{i}_{*}\ti{h}^!K$-acyclique d'après \ref{acprop}. Il ne reste plus qu'à appliquer le changement de base propre \ref{Base} pour conclure.
\epf

\subsection{Le cas des schémas placides}
Dans le cas des schémas placides formellement lisses, on va avoir des énoncés de fonctorialité plus forts ainsi que des énoncés d'existence de support singulier pour des faisceaux constructibles.
\subsubsection{Rappels sur le cas de type fini}
On commence par rappeler les résultats établis par Beilinson dans \cite{Bei}. Soit $\La$ un anneau fini artinien de caractéristique $\ell$ différente de la caractéristique de $k$, un $k$-schéma lisse de type fini $X$  et $K\in D^{b}_{c}(X,\La)$, Beilinson considère une notion de  $K$-acyclicité qui porte sur $h^*K$ au lieu de $h^!K$, définit l'ensemble $C_{\La}^{Bei}(K)$ des microsupports correspondants et montre qu'il admet un minimum $SS_{\La}^{Bei}(K)$ (\cite[Thm. 1.3]{Bei})\footnote{Stricto sensu, dans loc. cit. n'est considéré que le cas de $\bZ/\ell^n$, mais tous les arguments valent pour un anneau tel que $\La$. C'est écrit dans ce degré de généralité dans \cite{Sai}.}.

De manière explicite, $\cC_{\La}^{Bei}(K)$ consiste en les fermés côniques tels que pour toute paire $C$-transverse $(h,f)$ de $k$-schémas lisses de type fini, $f$ est $h^{*}K$-acyclique. On dit que $K$ est $(C,*)$-\textsl{microsupporté}. Quitte à passer aux composantes connexes, il est commode de supposer partout \og lisse\fg~au lieu de lisse de type fini. Pour la suite, on a besoin de considérer le cas des $\bql$-coefficients  ainsi que des complexes non-bornés.
On a l'énoncé suivant obtenu par Bartlett dans \cite[sect. 1.5, Thm.]{Bart}.
\bthm\label{Beiql}
Soit un $k$-schéma lisse, pour tout $K\in\cD_c(X)$, l'ensemble $C_{\bql}^{Bei}(K)$ admet un minimum $SS_{\bql}^{Bei}(K)$.
De plus, il existe toujours un modèle $K_{0}\in D^{b}_{c}(X,\bzl)$ tel que $K\simeq K_{0}\otimes_{\bzl}\bql$ avec $SS(K)=SS(K_0)$.
\ethm

Il s'agit maintenant de comparer les ensembles $\cC(K)$ et $\cC^{Bei}(K)$. De manière concrète sur un  $k$-schéma lisse, on va voir que  la propriété d'être un micro-support par rapport aux paires transverses constituées de schémas lisses implique la propriété d'être un microsupport pour toutes les paires transverses constituées de schémas placides formellement lisses.

\bprop\label{Beisupp}
Soient un $k$-schéma  $X$ lisse, $K\in \cD(X)$ et $C\subset T^{*}X$ un fermé cônique, alors $K$ est $C$-microsupporté pour la catégorie des $k$-schémas lisses si et seulement si il est $(C,*)$-microsupporté pour la même catégorie.
\eprop

\bpf
L'énoncé se ramène immédiatement à $X$ lisse de type fini et pour la catégorie des $k$-schémas lisses de type fini, ce que l'on suppose dans la suite.
Supposons que $K$ est $(C,*)$-microsupporté.
Par passage à la limite et comme les foncteurs en jeu commutent aux limites filtrantes, on se ramène immédiatement à $K\in\cD_c(X)$.
On considère une paire $C$-transverse $Y\stackrel{f}{\leftarrow}U\stackrel{h}{\rightarrow}X$ avec $U$ et $Y$ lisses de type fini. Comme $U$ est lisse, on a $h^{!}\bql=\bql[c]$ où $c$ est une fonction localement constante, donc pour montrer que $f$ est $h^!K$-acyclique, il suffit de montrer que:
\begin{equation}
h^{!}K\simeq h^{*}K\otimes^{L}_{\bql} h^{!}\bql.
\label{eqQl}
\end{equation}
Soit $K'\in\cD_{c}^{b}(X,\bzl)$ tel que $K'\otimes_{\bzl}\bql\simeq K$, d'après \ref{Beiql} $K'$ est $(C,*)$-microsupporté.
Maintenant, par réduction à un anneau artinien fini $\La$ et d'après Saito \cite[Prop. 8.13]{Sai}, on obtient que $h^{!}K'\simeq h^{*}K'\otimes^{L}_{\bzl} h^{!}\bzl$ et en localisant on obtient \eqref{eqQl}.

Réciproquement, si $K$ est $C$-microsupporté pour la catégorie des $k$-schémas de type fini. Soit $Z\stackrel{f}{\leftarrow}U\stackrel{h}{\rightarrow}X$ une paire test $C$-transverse de $k$-schémas de type fini lisses. On suppose que $f$ est $h^!K$-acyclique, il s'agit de voir que $f$ est aussi $h^*K$-acyclique.
La même réduction que ci-dessus nous ramène à $K\in\cD^{b}_c(X)$ et on peut supposer qu'il existe $K'\in\cD_{ctf}^{b}(X,\bzl)$ tel $K'\otimes_{\bzl}\bql\simeq K$ et $f$ est $h^!K'$-acyclique. On se ramène alors au cas où $K'\in\cD_{ctf}^{b}(X,\La)$. 
Par dualité, d'après Lu-Zheng et Gabber \cite[Cor. 5.13, Cor. 6.6]{LuZ}, $f$ est $h^*\bD(K')$-acyclique où $\bD$ est le dual de Verdier. On applique alors à nouveau \cite[Prop. 8.13]{Sai} pour obtenir que :
\[h^{!}\bD(K')\simeq h^{*}\bD(K')[c],\]
donc $f$ est $h^{!}\bD(K')$-acyclique et à nouveau par dualité, $f$ est $h^*K'$-acyclique, ce qu'on voulait.
\epf
\bprop\label{C2}
Soit un $k$-schéma $X$ lisse de type fini, $C\subset T^*X$ un fermé cônique, $Y\stackrel{f}{\leftarrow}U\stackrel{h}{\rightarrow} X$ une paire $C$-transverse avec $Y$ lisse de type fini et $U$ placidement présenté et formellement lisse, alors il existe $U_{\al}$ lisse de type fini tel qu'on ait un diagramme:
$$\xymatrix{U\ar[dd]_{f}\ar[drr]^{h}\ar[dr]\\&U_{\al}\ar[dl]^{f_\al}\ar[r]_{h_\al}&X\\Y}$$
avec $(h_\al,f_{\al})$  $C$-transverse et $U\rightarrow U_{\al}$ fortement pro-lisse.
\eprop
\bpf
 D'après \ref{C1}, il existe $\al$ tel qu'on ait un diagramme :
$$\xymatrix{U\ar[dd]_{f}\ar[drr]^{h}\ar[dr]_{p_{\al}}\\&U_{\al}\ar[dl]^{f_\al}\ar[r]_{h_\al}&X\\Y},$$
avec $h_{\al}$ qui est $C$-transverse et $U_{\al}$ lisse de type fini.
Il s'agit de voir que $f_\al$ est $h_\al^{\circ}C$-transverse aux points de $p_{\al}(U)$. Tout d'abord, d'après \eqref{Clis}, on a $h^{\circ}C=h_{\al}^{\circ}C\times_{U_{\al}}U$.
Soient $u\in U$, $u_{\al}=p_{\al}(x)$ et $y=f(x)$, la flèche 
$T^{*}_{y}Y\rightarrow T_u^{*}U$ se factorise par $T^{*}_{u_\al}U_{\al}$ et comme $f$ est $h^{\circ} C$-transverse en $u$, $f_{\al}$ l'est aussi en $u_{\al}$. Enfin, quitte à rétrécir $U_{\al}$ et comme la condition de $C$-transversalité est ouverte d'après \cite[Lem. 1.2]{Bei}, on peut donc supposer que $f_{\al}$ est $h_{\al}^{\circ}C$-transverse.
\epf
\bprop\label{fonct1}
Soit $X$ lisse et $K\in\cD(X)$, alors on a $\cC(K)=\cC^{Bei}(K)$.
\eprop
\bpf
Il suffit de montrer l'énoncé dans le cas lisse de type fini.
Tout d'abord, par \ref{Beisupp}, on a une inclusion $\cC(K)\subset\cC^{Bei}(K)$. Montrons la réciproque, soit $C\in\cC^{Bei}(K)$, on considère une paire $(h,f)$ $C$-transverse où $U$ microsupportable et $h$ lqc et $Y$ lisse de type fini, il s'agit de voir que $f$ est $h^!K$-localement acyclique. Comme $h$ est lqc, il est schématique et $U$ est donc à la fois un ind-schéma microsupportable et un schéma, en particulier il est ind-placide formellement lisse et comme l'énoncé est Zariski local, il suffit de montrer l'énoncé après restriction à tout ouvert quasi-compact $V\subset U$. Et donc on est ramené au cas où $U$ est placide quasi-compact et formellement lisse.
Avec \ref{pulLis}(i), on se réduit au cas où $U$ est formellement lisse avec une présentation placide.
En reprenant les notations de \ref{C2} et à nouveau par \ref{Beisupp}, on a déjà que $f_{\al}$ est $h_{\al}^!K$-acyclique, donc d'après \ref{pulLis} et comme $U\rightarrow U_{\al}$ est fortement pro-lisse, $f$ est $h^!K$-acyclique.
\epf

\bthm\label{SSex0}
Soit un $k$-schéma $X$ formellement lisse placide  alors pour tout $K\in\cD_c(X)$, $SS(K)$ existe et l'on a $B(SS(K))=\supp(K)$.
\ethm
\bpf
On commence par considérer le cas $X$ formellement lisse, placidement présenté.  D'après \ref{plac}, on peut considérer une présentation fortement pro-lisse $X\simeq\varprojlim X_i$ avec des flèches de transition affines, lisses. Comme $K\in\cD_c(X)$, on a $\cD_c(X)\simeq\colim\cD_c(X_i)$ et il existe $i$ tel que $K$ se descend en $K_i\in\cD_c(X_i)$.

Dans ce cas, d'après \ref{Beiql} et \ref{fonct1}, $SS(K_i)$ existe et $B(SS(K_i))=\supp(K_i)$ et  d'après \ref{Clis} et \ref{fonct}, on  $p_i^{\circ}SS(K_i)=SS(K_i)\times_{X_i}X\in C(K)$. Montrons que :
\begin{equation}
SS(K)=p_i^{\circ}(SS(K_i))=p_i^{*}(SS(K_i)).
\label{SS-eg}
\end{equation}
Soit $C\in\cC(K)$, d'après \cite[Thm 1.5]{Bei}, $SS(K_i)$ est l'adhérence des paires $(x_i,df_{x_i})$ où $x_i$ est un point fermé de $X_i$ et $f$ une fonction dans un voisinage Zariski $U_i$ de $x_i$ (éventuellement dont on étend les scalaires si $k$ est fini) qui n'est pas localement acyclique relativement à $K_i$ \footnote{Dans loc.cit., on demande seulement que $f$ ne soit pas localement acyclique relativement à $K_i$ en $x_i$, mais comme $X_i$ est un $k$-schéma de type fini, c'est une notion locale, donc on peut toujours supposer que la non-acyclicité est vraie sur tout $U_i$.}. Soit alors  une paire $(x_i, df_{x_i})\in SS(K_i)$ et $x\in X$ tel que $p(x_i)=x$ et $f$ définie sur un voisinage ouvert $U_i$. On note $\tilde{f}:U=p_{i}^{-1}(U_i)\stackrel{p_i}{\rightarrow} U_i\stackrel{f}{\rightarrow}\ab^1$ la composée.
Ainsi, $(x,d\tilde{f}_{x})\in p_i^{\circ}SS(K_i)$ et d'après \ref{pulLis}, $\tilde{f}$ n'est pas localement acyclique relativement à $K$ en $x$.
Si par l'absurde, $(x,d\tilde{f}_{x})\notin C$, alors $\tilde{f}$ est $C$-transverse en $x$ et d'après \ref{open}, quitte à rétrécir $U$, on peut supposer que $\tilde{f}$ est $C$-transverse sur $U$ et on obtiendrait alors que $K$ est $\ti{f}$-acyclique, une contradiction.
Ainsi, $C$ contient toutes les paires  $(x,d\tilde{f}_{x})$ avec $\ti{f}$ non localement acycliquement en $x$ relativement à $K$ et donc leur adhérence, qui s'identifie à $p_i^{\circ}SS(K_i)$ par platitude topologique comme $p_i$ est fortement pro-lisse, donc plate.
On obtient donc que $SS(K)=p_i^{\circ}(SS(K_i))$ comme souhaité.

Passons maintenant au cas où $X$ est seulement placide.
On a d'abord d'après \eqref{SS-eg} et \cite[Thm.1.4]{Bei} après réduction à un cran fini, que pour tout morphisme lisse $f:T\ra S$ entre schémas formellement lisses placidement présentés et tout $K\in\cD_c(S)$, on a:
\[f^{\circ}SS(K)=SS(f^!K),\]
et le résultat s'obtient par descente étale.
Enfin, l'identité sur  $B(SS(K))$ se déduit de l'identité analogue à un cran fini.
\epf
Pour tout champ placide $\cX$ formellement lisse, on dispose d'après \cite[3.2.4]{B} d'une catégorie $\Loc^{ren}(\cX)$ de systèmes locaux renormalisés.

Si $\cX$ est de plus un schéma placidement présenté, on a :
\[\Loc^{ren}(\cX)\cong\colim_{X\ra Y} \Loc^{ren}(Y),\]
 où la colimite est faite sur tous les morphismes $X\ra Y$ fortement pro-lisses avec $Y$ un $k$-schéma de type fini et $\Loc^{ren}(Y)=\{\omega_{Y}\otimes_{\bql}\cL,\cL\in \Loc_{\bql}(Y)\}$ où $\Loc_{\bql}(Y)$ est la catégorie habituelle des systèmes locaux sur $Y$.
\bcor\label{SSex-bis}
\benumr
\remi
Soit $X$ un $k$-schéma placide formellement lisse, alors pour tout système local $\cL\in\Loc^{ren}(X)$, $SS(\cL)=\{0\}_X$.
\remi
Soit un morphisme $f:X\ra Y$ fortement pro-lisse entre schémas placides formellement lisses, alors  pour tout $K\in\cD_c(Y)$, on a :
\[f^{\circ}SS(K)=SS(f^{!}K).\]
\eenum
\ecor
\bpf
(i) On se ramène immédiatement par descente au cas où $X$ est formellement placidement présenté et on a alors $p^{!}\cL_0=\cL$ avec $p:X\ra Y$, $Y$ lisse de type fini et $\cL_{0}\in\Loc^{ren}(Y)$. L'énoncé se déduit alors de \eqref{SS-eg} et du fait que $SS(\cL_{0})=\{0\}_{Y}$.

(ii) Comme $f$ est fortement pro-lisse, on a clairement l'inclusion $SS(f^{!}K)\subset f^{\circ}SS(K)$.
Comme le résultat est local pour la topologie étale, on peut se ramener à $Y$ affine, placidement présenté et formellement lisse et comme $K\in\cD_c(Y)$, on peut supposer que $Y$ est lisse de type fini et cela se déduit de \eqref{SS-eg}.
\epf

\subsection{Théorème d'existence pour des ind-schémas}
\bprop\label{ex-ind1}
Soit un $\aleph_0$-ind-schéma ind-(placide formellement lisse) $Y$, soit $K\in\cD_{c}(Y)$ à support borné, alors $SS(K)$ existe.
\eprop
\bpf
Soit $K\in\cD_{c}(T)$ à support borné, il est donc supporté par un sous-schéma fermé placide formellement lisse $Z$ et $i:Z\hra Y$ est de présentation finie.
D'après \ref{SSex0}, $SS(i^{!}K)$ existe et donc d'après \ref{fonct-clo}, on en déduit que :
\[i_{\circ}SS(i^{!}K)\in\cC(K=i_{*}i^{!}K).\]
De plus, pour tout $C\in \cC(K)$ et une paire $(x,df_{x})\in\cC$, $K$ est localement acyclique en $x$ pour $f_x$ si et seulement si $(f_{x})_{\vert Z}$ l'est pour $i^{!}K$. Ainsi, $i_{\circ}SS(i^{!}K)\subset C$, comme souhaité.
\epf
Pour les applications, ce cas n'est pas suffisant étant donné que beaucoup de ind-schémas formellement lisses ne sont pas ind-formellement lisses, à commencer par la grassmannienne affine, par exemple. 
On a donc besoin d'un énoncé plus général.
\begin{defi}\label{plong}
Soit un ind-schéma $X$ microsupportable. Il est dit plongeable s'il admet une  immersion fermée de présentation finie $i:X\hra Y$ dans un $\aleph_0$-ind schéma ind-(placide formellement lisse) et que la catégorie de tels plongements est filtrante.
\end{defi}
\brem
On vérifie aisément, par réduction au cas des $k$-schémas de type fini, que si $X$ est un $\aleph_0$-ind-schéma de ind-type fini, alors il est plongeable.

Il en résulte immédiatement qu'il en est de même pour tout ind-schéma raisonnable $X$ qui s'écrit $X\simeq\colim X_{\al}$ avec $X_{\al}\simeq Y\times\ab^{\NN}$ et $Y$ un $k$-schéma de type fini. Et par suite, on peut même supposer que les $X_{\al}$ sont de cette forme localement pour la topologie étale.

Il résulte alors de \cite[Thm. 6.3]{Dr} que tout espace de lacets $\cL X$ pour $X$ affine lisse est plongeable.
\erem

L'énoncé suivant est l'énoncé d'existence souhaité:

\bthm\label{ex-ind2}
Soit un ind-schéma $X$ microsupportable et plongeable. Soit $K\in\cD_{c}(X)$, borné (cf. \ref{borne}) alors il existe un fermé cônique $C\in\cC(K)$ tel que pour tout plongement $i:X\hra Y$ dans $Y$ ind-(placide formellement lisse), on a:
\[(di)^{-1}(C)=SS(i_*K).\]
\ethm
La preuve est analogue à \cite[E. 6.4]{AGKRRV} une fois que l'on a établi \ref{fonct-clo}:
\bpf
Soit un plongement $i:X\hra Y$ avec $Y$ $\aleph_0$-ind-(placide formellement lisse).
Il suffit de vérifier que pour tout $x\in X$, l'intersection:
\[T^{*}_{x}Y\cap SS(i_{*}K)\]
est invariante par translation par les éléments de :
\[Ker(T_x^{*}Y\ra T_x^{*}X).\]
Si l'on dispose de deux plongements $i_{1}:X\ra Y_1$ et $i_{2}:X\ra Y_2$ par hypothèse, on peut compléter le diagramme en:
$$\xymatrix{Y\ar[r]^{i_1}\ar[d]_{i_2}&Y_{1}\ar[d]^{i'_{2}}\\Y_{2}\ar[r]^{i_{1}'}&Y_{12}}$$
avec $Y_{12}$ ind-(placide formellement lisse).
Ainsi, pour un plongement $(Y_1,i_1)$ et une paire $(\xi_1,\xi_1')\in T^{*}_{x} Y_1$ qui se projettent sur le même vecteur de $T^{*}_{x}X$, on peut donc trouver un diagramme commutatif et une paire  $(\xi_{12},\xi_{12}')$ qui s'envoie sur $(\xi_1,\xi_1')$ par $T^{*}_{x}Y_{12}\ra T^{*}_{y}Y_1$ et sur un même élément $\xi_2\in T^{*}Y_2$ par $T^{*}_{x}Y_{12}\ra T^{*}_{y}Y_2$. 
En utilisant alors \ref{fonct-clo} pour $Y_1$ et $Y_2$ et \ref{ex-ind1}, on trouve:
\[  \xi_{1}\in SS((i_{1})_{*}K)\Leftrightarrow\xi_{12}\in SS((i'_{2}\circ i_{1})_{*}K)=SS((i'_{1}\circ i_{2})_{*}K)\Leftrightarrow\xi_{2}\in SS((i_{2})_{*}K) \\
  \]
\[\Leftrightarrow\xi'_{12}\in SS((i'_{1}\circ i_2)_{*}K)= SS((i'_{2}\circ i_1)_{*}K)\Leftrightarrow\xi'_1\in SS((i_{1})_{*}K).\]
comme souhaité.
\epf

\subsection{Situation équivariante}\label{cal-SS}
On considère un $k$-ind-schéma $X$ microsupportable avec une action de $\cL G$ pour un groupe $G$ connexe réductif sur $k$.
Soit $H\subset \cL G$ un sous-groupe fortement pro-lisse tel que $\cL G/H$ est ind-fp-propre. On peut prendre en particulier $H=\clp G$ où $H$ un sous-groupe d'Iwahori.
D'après \ref{loopeq} et \ref{desc-ind-fp}, on dispose d'une flèche de $\infty$-champs $\pi: [X/H]\ra [X/\cL G]$, ind-fp-propre.
En particulier, d'après \ref{Base}, on dispose d'un foncteur $\pi_!$.
Partant d'un faisceau $K\in\cD(X)$ $H$-équivariant, on définit l'induction de $K$ comme :
\[\Ind_{H}^{\cL G}(K)=\pi_{!}K_0,\]
où $K_0$ est l'unique faisceau de $\cD([X/H])$ tel que $p^{!}K_0=K$.
L'énoncé suivant permet d'obtenir une \og borne\fg~ pour les microsupports de $\Ind(K)$ en fonction de ceux de $K$.
\bprop\label{SS-indu}
Soit un $k$-ind-schéma $X$ microsupportable avec une action de $\cL G$, $H\subset \cL G$ un sous-groupe fortement pro-lisse.
Soit $K\in\cD(X)$, $H$-équivariant et $C\in\cC(K)$ tel que $C=p^{\circ}C_0$ pour un fermé cônique $C_0\subset T^{*}[X/H]$, alors $[(\ov{LG.C}\cap\mu_{\cL G}^{-1}(0))/\cL G]\in\ov{\cC}(\Ind_{H}^{\cL G}(K))$.
\eprop
\bpf
Comme $C=p^{\circ}C_0$, on déduit de \eqref{cotanH2} que $C$ est $H$-invariant et $C\subset\mu_{H}^{-1}(0)$.
D'après \ref{equivquot}, comme $p_{\circ}p^{\circ}C_0=C_0$, $C_0\in\cC(K_0)$ et d'après \ref{propfonct}, on a 
$\ov{\pi_{\circ}C_0}\in\ov{\cC}(\Ind_{H}^{\cL G}(K))$. Il s'agit d'avoir une description plus explicite de $\ov{\pi_{\circ}C_0}=\ov{\ti{\pi}_{\circ}C}$ avec $\ti{\pi}:X\ra [X/\cL G]$. Calculons $\ti{\pi}^{\circ}\ti{\pi}_{\circ}C$, puisque d'après \eqref{cotanH2}, on a $\ti{\pi}_{\circ}C=[\ti{\pi}^{\circ}\ti{\pi}_{\circ}C/\cL G]$.
On forme le carré cartésien suivant:
$$\xymatrix{X\times_k\cL G\ar[d]_{q}\ar[r]^-{a}&X\ar[d]^{\ti{\pi}}\\X\ar[r]^-{\ti{\pi}}&[X/\cL G]}$$
où $a$ est l'action et $q$ la projection. D'après \ref{Cbc}(i), on a $\ti{\pi}^{\circ}\ti{\pi}_{\circ}C=q_{\circ}a^{\circ}C$.
En particulier, on a déjà d'après \eqref{cotanH} que :
\begin{equation}
q_{\circ}a^{\circ}C\subset\mu^{-1}_{\cL G}(0).
\label{moment}
\end{equation}
En notant également $q$ la projection $T^{*}X\oplus T^{*}\cL G\ra T^{*}X$, on a immédiatement que :
\[q_{\circ}a^{\circ}C\subset q(a^{\circ}C)\]
ainsi que $a^{\circ}C=a^{*}C=C\times_{X}(X\times\cL G)$. Pour toute paire $(x,g)\in X\times_k\cL G$, on a 
\[a^{*}C\cap T_{(x,g)}^{*}(X\times_k\cL G)=da_{(g,x)}(C\cap T^{*}_{g.x}X)\]
 et la composée $q\circ da_{(g,x)}: T^{*}_{g.x}X\ra T_{x}^{*}X$ s'identifie à la différentielle de l'action de $g$ sur $X$ d'où l'identité $q(a^{\circ}C)=\cL G. C$, soit en combinant avec \eqref{moment}:
\begin{equation}
q_{\circ}a^{\circ}C\subset\cL G. C\cap\mu_{\cL G}^{-1}(0)\subset\ov{\cL G.C}\cap \mu_{\cL G}^{-1}(0)
\label{moment2}
\end{equation}
Il suffit maintenant de voir que $\ti{\pi}^{\circ}(\ov{\ti{\pi}_{\circ}C})=\ov{\ti{\pi}^{\circ}(\ti{\pi}_{\circ}C})$, ce qui est clair à nouveau par \eqref{cotanH2} et de la description de la topologie sur le champ quotient. Comme le dernier terme de droite de \eqref{moment2} est fermé, on obtient l'inclusion:
\[\ov{\pi_{\circ}C_0}\subset[(\ov{LG.C}\cap\mu_{\cL G}^{-1}(0))/\cL G],\]
ce qu'on voulait.
\epf
\section{Microsupports pour le faisceau de Grothendieck-Springer affine}
\subsection{Une fibration et sa stratification}
\subsubsection{Stratification par les valuations radicielles}
Soit un groupe $G$ connexe réductif sur un corps $k$ algébriquement clos. Soit $(B,T)$ une paire de Borel, $W$ le groupe de Weyl, $\kg=\Lie(G)$, $\kt=\Lie(T)$ et $\kc:=\kt/W$. Si $\car(k)\ssup 2$ ou $G$ sans composante de type $(C_n)$, on dispose d'après \cite[Thm. 4.1.10]{BC} d'un isomorphisme de Chevalley $k[\kg]^{G}\cong k[\kt]^{W}$ qui induit un diagramme commutatif:
\begin{equation}
\xymatrix{\kt\ar[dr]_{\pi}\ar[r]&\kg\ar[d]^{\chi}\\&\kc},
\label{Ccom}
\end{equation}
où la flèche horizontale est l'inclusion canonique.
Soit le discriminant $\kD=\prod\limits_{\al\in R}d\al\in k[\kc]$ où $R$ désigne l'ensemble des racines, on note $\kc^{rs}$ le complémentaire de ce diviseur et $\kg^{rs}=\chi^{-1}(\kc^{rs})$. 

On suppose  désormais l'ordre de $W$ premier à la caractéristique.
D'après \cite[2.5.12]{Let}, il existe un isomorphisme $G$-équivariant $\kappa:\kg\stackrel{\sim}{\rightarrow}\kg^{\vee}$ que l'on fixe dans la suite.
L'isomorphisme $G$-équivariant $\kappa:\kg\stackrel{\sim}{\rightarrow}\kg^{\vee}$ induit un accouplement parfait :
\[\kappa^{aff}:\cL\kg\times\cL\kg\ra\ab^1,\]
donné par $(x,y)\mapsto\Res\kappa(x,y)$.

Soient $m=\vert W\vert$, $\cO=k[[t]]$ et $F=k((t))$, $\co'=k[[t^{\frac{1}{m}}]]$, $F'=\Frac(\co')$, $\kt'=\Res_{\cO'/\co}(\kt\times_k\cO')$ et $\kc'=\Res_{\co'/\co}(\kc\times_k\cO')$.
On définit alors $T_w=\Res_{\cO'/\co}(\kt\times_k\cO')^{w\mu_{m}}$ où $\mu_{m}$ agit sur $\cO'$ par l'action de Galois. D'après \cite[Lem. 15.3.1, 15.4.1]{GKM}, c'est un schéma en groupes lisse sur $\co$, qui est un tore sur $F$ et $\kt_w:=\Lie(T_w)=\Lie(\kt')^{w\mu_m}$.
D'après \cite[4.2, 4.3]{GKM2}, $T_{w}$ se déploie en $T$ sur $F'$ et on a un plongement $T_w\hra G$ unique à $G$-conjugaison près. De plus, tout tore maximal de $G$ est $G$-conjugué à un $T_w$, pour $w\in W$.
En particulier, tout élément $\g\in\kg^{rs}(F)$ est $G(F)$-conjugué à un élément de $\kt_{w}(F)$.
En considérant la composée $\kt_w\ra\kg\ra\kc$, on obtient une flèche canonique, indépendante des choix, définie sur $F$, $\pi_w:\kt_w\ra\kc.$
\blem\label{w-fin}
La flèche $\pi_w$ est finie.
\elem
\bpf
Par descente, il suffit de montrer que le morphisme est fini après extension de scalaires $F'/F$. Sur $F'$, $T_w$ se déploie en $T$ et la flèche $\pi_{w}\times_{F} F'$ s'identifie à la flèche $\pi:\kt\ra\kc$, qui est également finie, ce qu'on voulait.
\epf
Pour toute fonction $r:R\ra\frac{1}{m}\NN$ sur l'ensemble des racines, on peut considérer le sous-schéma localement fermé de présentation finie $\kt'_{r}$ où l'on fixe l'ensemble des valuations radicielles. En intersectant avec $\clp\kt_w$, on obtient ainsi la strate $\kt_{w,r}$.
De manière concrète, $\clp\kt_{w}$ classifie les séries $\sum\g_{i}t^{\frac{i}{m}}$ telles que $w^{-1}(\g_{i})=\xi^{i}\g_{i}$ et $\kt_{w,r}\subset\clp\kt_{w}$ est donnée par les équations:
\begin{equation}
d\al(\g_{i})=0~~ \text{si}~~0\leq i\sinf r(\al)~~ \text{et}~ d\al(\g_{i})\neq0~~ \text{si}~~ r(\al)=i.
\label{eqsplit}
\end{equation}
Ainsi, $\kt_{w,r}$ est donné par un nombre fini d'équations et d'inéquations linéaires, donc $\kt_{w,r}$ est fortement pro-lisse.
Pour $w\in W$ et une fonction $r$, on considère $u(r):R\ra\frac{1}{m}\NN$ donnée par $u(r)(\al)=r(u^{-1}.\al)$. 
On a donc une action de $W$ sur les paires $(w,r)$  donnée par $u.(w,r)=(uwu^{-1},u(r))$.
On note $\kc_{w,r}\subset\clp\kc$ l'image de $\kt_{w,r}$. D'après \cite[Thm. 8.2.2]{GKM}, c'est un localement fermé de présentation finie, irréductible et fortement pro-lisse et  d'après \cite[3.3.4]{BKV}, on a une flèche finie étale $\kt_{w,r}\ra\kc_{w,r}$ de groupe $W_{w,r}$, où $W_{w,r}\subset W$ est le stabilisateur de $(w,r)$ pour l'action de $W$ sur les paires $(w,r)$.
La flèche $\kg^{rs}\ra\kc^{rs}$ étant surjective lisse \cite[4.2.6]{BC}, la flèche induite $\cL\kg^{rs}\ra\cL\kc^{rs}$ est formellement lisse, et comme $\kc_{w,r}$ est fortement pro-lisse, $\kC_{w,r}=\kc_{w,r}\times_{\cL\kc^{rs}}\cL\kg^{rs}$ est formellement lisse.

\blem\label{orthdec}
Pour tout $s\in\kg^{rs}(F)$, on a une décomposition orthogonale sur $F$:
\[\kg=\kg_{s}\oplus[s,\kg]=\Lie(C_{G}(s))\oplus[s,\kg],\]
avec $\kg_{s}=\Ker([s,.])$ et $C_{G}(s)$ le centralisateur de $s$.
\elem
\bpf
Tout d'abord, comme $s$ est semisimple, $C_{G}(s)$ est $F$-lisse d'après \cite[III.9.1, Prop.]{Bor} et on a $\Lie(C_{G}(s))=\kg_{s}$ d'après \cite[1.10]{Hum}.
En considérant $\kappa$ comme une forme bilinéaire symétrique, comme elle est $G$-équivariante d'après \cite[10.7]{Hum2}, on a:
\[\forall~(x,y,z)\in\kg^3, \kappa(x,[y,z])=\kappa([x,y],z).\]
On a ainsi l'inclusion $\kappa(\kg_{s})\subset(\kg/[s,\kg])^{*}$ et l'égalité par un argument de dimension.
Enfin, pour vérifier que la somme est directe, il suffit de le faire après une extension finie $L/F$ qui déploie $s$, de telle sorte que le centralisateur de $s$ s'identifie au tore maximal $T$ et  la décomposition s'identifie alors à $\kg=\kt\oplus(\kn\oplus\kn^{-}).$
\epf
\subsubsection{La fibration de Grothendieck-Springer affine}
On applique maintenant les résultats de \ref{cal-SS} dans le cas de l'action adjointe de $\cL G$ sur $\cL\kg$.
Pour $\cL\in\Loc^{ren}([\clp\kg/\clp G])$, on pose alors $\cS_{\kg,\cL}=\Ind_{\clp G}^{\cL G}(i_{*}\cL)\in\cD([\cL \kg/\cL G])$ avec $i:\clp\kg\hra\cL\kg$.
\subsection{Un calcul de microsupport}

Le morphisme de Chevalley induit une flèche $\cL\chi:\cL\kg\ra\cL\kc$ et on considère le sous-ind-schéma fermé des éléments compacts $\kC=\cL\chi^{-1}(\clp\kc)$.  Comme $\pi_{w}$ est finie d'après \ref{w-fin}, on obtient par critère valuatif l'inclusion:
\begin{equation}
(\kC\cap\cL\kt_w)(k)\subset(\clp\kt_w)(k).
\label{C-incl}
\end{equation}
Enfin, on introduit les ouverts génériquement semisimples $\kC^{\bullet}=\chi^{-1}(\clp\kc^{\bullet})=\cL\chi^{-1}(\clp\kc-\clp\kD)$  et $\wkC^{\bullet}=p^{-1}(\kC^{\bullet})$. Ce ne sont pas des ind-schémas puisque $\clp\kc^{\bullet}$ n'est pas quasi-compact. En revanche, lorsque l'on borne la valuation du discriminant, on obtient des ouverts quasi-compacts.
On s'intéresse aux micro-supports de $\cS_{\kg}$.
Le théorème fondamental est le suivant:
\bthm\label{micro-comp}
 Posons $C_{\kg}:=(\kC^{\bullet}\times\kC)\cap\{(\g,\xi)\in\cL\kg\times\cL\kg, [\g,\xi]=0\}$ alors  $[C_{\kg}/\cL G]$ est un microsupport pour $\cC(\cS_{\kg})$.
\ethm
\bpf
On décompose la preuve en plusieurs étapes, on commence par le cas sphérique (i).

\textbf{\'{E}tape 1}: Détermination d'une borne supérieure pour le support singulier de $i_{*}\cL$ pour  $\cL\in\Loc^{ren}([\clp\kg/\clp G])$.

On note de la même manière le pullback de $\cL$ à $\clp\kg$.
Le schéma $\clp\kg$ est placide formellement lisse, de telle sorte que d'après \ref{SSex-bis} $SS(\cL)=\{0\}_{\clp\kg}$. On applique alors \ref{fonct-clo}
à l'immersion de présentation finie $i:\clp\kg\hra\cL\kg$ de telle sorte que :
\[i_{\circ}\{0\}_{\clp\kg}=T^{*}_{\clp\kg}\cL\kg\in\cC(i_{*}\omega_{\clp G}).\]
Or, on a l'égalité $(\clp\kg)^{\perp}=\clp\kg$, d'où l'on déduit:
\begin{equation}
\clp\kg\times\clp\kg\in\cC(i_{*}\omega_{\clp G}).
\label{etap1}
\end{equation}
\medskip

\textbf{\'{E}tape 2}: Descente à $[\cL\kg/\clp G]$.

Il s'agit de calculer l'application moment
\[\mu_{\clp G}:T^{*}\cL\kg\ra (\clp\kg)^{*}\]
 donnée par $(x,\xi)\mapsto da_{x}(\xi)$ où $a_x:\clp G\ra\cL\kg$ est donnée par $g\mapsto\ad(g).x$. On a une identification à l'aide de $\kappa^{aff}$, $(\clp\kg)^{*}\cong t^{-1}\kg(k[t^{-1}])$
et $\mu$ est donnée par $(x,\xi)\mapsto \theta([x,\xi])$ avec $\theta:\cL\kg\ra[\cL\kg/\clp\kg]\cong t^{-1}\kg(k[t^{-1}])$. En particulier, on obtient que :
\[\mu_{\clp G}^{-1}(0)=\{(g,\xi)\in(\cL\kg)^{2}, [g,\xi]\in\clp\kg\}.\]
Ainsi, $C=\clp\kg\times\clp\kg\subset\mu_{\clp G}^{-1}(0)$ et est $\clp G$-équivariant de telle sorte que $C=p^{\circ}C_0$ et d'après \ref{equivquot} et \eqref{etap1}, $C_{0}\in\cC(K_0)$ où $K_0$ est l'unique faisceau de $\cD([\cL\kg/\clp G])$ tel que $p^{!}K_0\cong K$.
\medskip

\textbf{\'{E}tape 3}: Passage à l'induction.

L'application moment $\mu_{\cL G}$ pour l'action adjointe de $\cL G$ sur $\cL\kg$ est donnée par le crochet.
Ainsi, on obtient $\mu_{\cL G}^{-1}(0)=\{(\g,\xi)\in\cL\kg\times\cL\kg, [\g,\xi]=0\}$.
De plus, on a :
\[\cL G.(\clp\kg\times\clp\kg)\subset\kC\times\kC,\]
où $\cL G$ agit par l'action diagonale. Comme $\kC\times\kC$ est fermé, on a aussi $\ov{\cL G.(\clp\kg\times\clp\kg)}\subset\kC\times\kC$.
On applique alors \ref{SS-indu} et on obtient:
\[[C_{\kg}/\cL G]\subset\ov{\cC}(\cS_{\kg}),\]
avec $C_{\kg}:=(\kC\times\kC)\cap\{(\g,\xi)\in\cL\kg\times\cL\kg, [\g,\xi]=0\}$. Il reste à voir qu'en fait $[C_{\kg}/\cL G]\subset\cC(\cS_{\kg})$, i.e. que $B(C_{\kg})$ est fermé dans $\cL\kg$, ce qui est clair puisque :
\begin{equation}
B(C_{\kg})=\kC.
\label{baseS}
\end{equation}
\end{proof}
\subsection{Interaction avec la stratification par valuations radicielles}
\bthm\label{incl-loc}
On a l'inclusion:
\[[C^{\bullet}_{\kg}/\cL G]\subset\coprod_{(w,r)}T_{[\kg_{w,r}/\cL G]}^{*}[\cL\kg/\cL G].\]
\ethm
\bpf
On a une flèche canonique:
\[\coprod_{w\in W}\cL G\times\cL\kt_{w}^{rs}\ra\cL\kg^{rs}\]
qui se factorise en $\pi_1:\coprod_{w\in W}\cL G\times\cL\kt_{w}^{rs}\ra\coprod_{w\in W}\cL \kg^{rs}\times_{\cL\kc^{rs}}\cL\kt_{w}^{rs}$ et $\pi_{2}:\coprod_{w\in W}\cL \kg^{rs}\times_{\cL\kc^{rs}}\cL\kt_{w}^{rs}\ra\cL\kg^{rs}$.
La première flèche induit alors d'après \cite[Preuve de 4.1.8]{BKV} et \cite[Cor. 2.6.5]{B}, un isomorphisme:
\[\coprod_{w\in W}[\cL\kt_{w}^{rs}/\cL T_w]\simeq \coprod_{w\in W}[\cL \kg^{rs}\times_{\cL\kc^{rs}}\cL\kt_{w}^{rs}/\cL G]\]
et la flèche $\pi_2$ est un morphisme fini étale, par changement de base et \cite[Prop. 3.1.11]{BKV}.
En particulier, en passant aux strates, on a un isomorphisme:
\begin{equation}
[\kt_{w,r}/(\cL T_w\rtimes W_w)]\cong[\kg_{w,r}/\cL G],
\label{r-quot}
\end{equation}
avec $W_w\subset W$ le stabilisateur de la strate $\kt_{w,r}$.
Enfin, la flèche $\cL\kg^{rs}\ra\cL\kg$ est formellement étale de telle sorte que pour un point $t\in\kt_{w,r}(K)$, on a un isomorphisme d'espaces vectoriels:
\[T^{*}_{t}[\kt_{w,r}]\cong T^{*}_{t}[\kg_{w,r}/\cL G].\]
et la flèche :
$T^{*}_{t}[\cL\kg/G]\ra T^{*}_{t}[\kg_{w,r}/\cL G]$
s'identifie à $T^{*}_{t}\cL\kt_w\ra T^{*}_{t}\kt_{w,r}$.
On déduit de \ref{orthdec}, une décomposition orthogonale $G$-équivariante:
\[\cL\kg=\cL\kg_{t}\oplus[t,\cL\kg]\]
où $\kg_{t}$ est le centralisateur de $t$ qui s'identifie à $\cL\kt_{w}$ comme $t\in\kt_{w,r}$.
On a donc en particulier que $(\cL\kt_{w})^{*}=([t,\cL\kg])^{\perp}=\cL\kt_{w}$ et la forme bilinéaire $\kappa$ reste non-dégénérée quand on la restreint à $\cL\kt_w$.
Maintenant pour $(t,\xi)\in C_{\kg}^{\bullet}(K)$, on a donc $\xi\in T^{*}\cL\kt_w=\cL\kt_w$ et $\xi\in\kC$. Or, on a $\chi(\xi)=\pi_{w}(\xi)\in\kc(K[[t]])$ avec $\pi_{w}:\kt_{w}\ra\kc$ qui est fini daprès \ref{w-fin}, donc par critère valuatif $\xi\in\clp\kt_{w}$.

Maintenant, on a la série d'inclusions $\kt_{w,r}\ra\clp\kt_{w}\ra\cL\kt_w$ de telle sorte que $\clp\kt_w=(T_{t}\clp\kt_{w})^{\perp}\subset (T_{t}\kt_{w,r})^{\perp}$,
d'où l'on déduit que :
\[(t,\xi)\in T_{[\kg_{w,r}/\cL G]}^{*}[\cL\kg/\cL G],\]
comme souhaité.
\epf
\brem
En dimension finie sur $\bC$ et si l'on sait que la stratification est de Whitney, il résulte d'un théorème de Kashiwara-Schapira \cite[Prop. 8.4.1]{KS} que l'énoncé \ref{incl-loc} implique la locale constance. Dans ce contexte, un tel énoncé n'est pour l'heure pas disponible. En fait, la principale difficulté vient du fait que l'on manque d'informations sur la stratification sur les $\kc_{w,r}$, la nature des adhérences, est-elle de Whitney, etc...
\erem
\section{Appendice A : Classes de modules topologiques}
\subsection{Catégorie des pro-modules}
\subsubsection{Définitions}
Soit $A$ un anneau et $M=(M_{\al})\in\Pro(\Mod^{\heartsuit}_{A})$, pour tout morphisme $A\ra A'$, on a un foncteur de changement de base de la catégorie des $A$-pro-modules dans celle des $A'$-pro-modules que l'on note $M\mapsto M\widehat{\otimes}_{A}A'$.
Tout $A$-module peut être vu comme un $A$-pro-module. On peut obtenir une autre classe de pro-modules via la construction suivante.
Un $A$-module est un ind-objet dans la catégorie des $A$-modules de présentation finie de telle sorte que le foncteur $M\mapsto M^{\vee}=\Hom_{A}(M,A)$ s'étend en un foncteur:
\begin{equation}
(-)^{\vee}:\Mod^{\heartsuit}_{A}\ra \Pro(\Mod^{\heartsuit}_{A}).
\label{eqdual}
\end{equation}
La restriction de ce foncteur à la sous-catégorie des $A$-modules plats est pleinement fidèle, en effet,
par le théorème de Govorov-Lazard (\cite{Gov}, \cite{Laz}), tout $A$-module plat est une colimite filtrante de $A$-modules projectifs de type fini de telle sorte que $M\cong \Hom_{A}(M^{\vee},A)$.

Un $A$-pro-module est dit $\aleph_0$ ou dénombrablement engendré s'il peut-être représenté par un système projectif dénombrable.
\subsubsection{Pro-modules et modules topologiques}\label{TMod}
Soit $A$ un anneau et $M=(M_{\al})\in\Pro(\Mod^{\heartsuit}_{A})$, on a foncteur canonique vers la catégorie des $A$-modules topologiques:
\[F:\Pro(\Mod^{\heartsuit}_{A})\ra\TopMod_{A}\]
donnée par $(V_{\al})\mapsto\varprojlim V_{\al}$ qui admet un adjoint à gauche :
\[G:\TopMod_{A}\ra\Pro(\Mod^{\heartsuit}_{A}).\]
qui envoie un $A$-module topologique $M$ sur le système projectif de ses quotients discrets.
En général, les morphismes d'adjonction $G\circ F\ra \Id$ et $\Id\ra F\circ G$ ne sont pas des isomorphismes. Déterminons quand peut-on remplacer un pro-module par un $A$-module topologique usuel.
\bdefi\label{pro-ml}
Soit $A$ un anneau commutatif et $M=(M_{\al})\in\Pro(\Mod^{\heartsuit}_{A})$.
\benumr
	\remi 
	On dit que $M$ est un pro-ensemble de Mittag-Leffler (resp. strictement de Mittag-Leffler), s'il est équivalent à un système projectif filtrant $(M'_{\al})$ où les flèches de transition sont surjectives  (resp. et de surcroît les flèches $\varprojlim M'_{\al}\ra M'_{\beta}$ sont surjectives pour tout $\beta$).
	\remi
Soit un $A$-module plat $M$, il est dit de Mittag-Leffler (resp. strictement de Mittag-Leffler) si $M^{\vee}$ est un pro-ensemble de Mittag-Leffler (resp. strictement de Mittag-Leffler).
\eenum
\edefi
\brems
\remi\label{mldef1}
D'après \cite[Tag. 0597]{Sta}, si le système projectif est dénombrable, alors il y a équivalence entre être de Mittag-Leffler et strictement de Mittag-Leffler.
\remi\label{mldef2}
Si $N\in\TopMod_{A}$, alors par construction $G(N)$ est un pro-ensemble strictement de Mittag-Leffler. En particulier, pour $M\in\Pro(\Mod^{\heartsuit}_{A})$, la flèche $G\circ F(M)\ra M$ est un isomorphisme si et seulement si $M$ est un pro-ensemble strictement de Mittag-Leffler.
\erems

Etudions maintenant pour quelle classe de $A$-modules topologiques $M$, la flèche $M\ra F\circ G(M)$ est
un isomorphisme.
Soit $A$-module topologique $M$, on dit qu'il est \textsl{linéairement topologisé} s'il admet une base de voisinages de 0 formée par des sous-modules $(M_{\al})$.
Un $A$-module topologique $M$ linéairement topologisé est dit \textsl{complet} et \textsl{Hausdorff }si :
\[M\cong\varprojlim M/M_{\al}.\]

\brem\label{adj-lin}
Par construction, si $M$ est un $A$-module topologique, linéairement topologisé complet et Hausdorff, alors la flèche $M\ra F\circ G(M)$ est
un isomorphisme.
\erem
On peut  donc considérer ces modules topologiques comme des $A$-pro-modules.
Il y a une classe particulièrement utile de modules topologiques, celles des modules de Tate.
\subsubsection{Modules de Tate}
Soit un $A$-module $M$, on peut considérer son dual $M^{\vee}=\Hom_{A}(M,A)$ qui admet une structure de $A$-module topologique pour laquelle les ouverts sont les complémentaires orthogonaux des sous-modules de présentation finie. On rappelle la définition suivante due à Drinfeld \cite[3.2.1]{Dr}:
\bdefi\label{Dr-Tate}
Soit un anneau commutatif $A$, un $A$-module de Tate est  un facteur direct d'un $A$-module topologique de la forme $P\oplus Q^{\vee}$ pour $P$ et $Q$ des $A$-modules projectifs discrets.
\edefi
Si $M$ est un $A$-module topologique, un sous-module $L\subset M$ est un \textsl{réseau} s'il est ouvert et si pour tout sous-module ouvert $U\subset L$, $L/U$ est de type fini.
On a la propriété suivante pour les $A$-modules de Tate, \cite[Thm.3.2]{Dr}:
\bthm\label{Tate-prop}
Soit $M$ un $A$-module de Tate alors on a :
\benumr
	\remi
	$M$ admet une base de voisinages de zéro formée par des réseaux.
		\remi 
	$M$ est complet et Hausdorff.
	\remi
	Le foncteur qui à tout $A$-module $N$ associe le group $\Hom_c(M,N)$ des morphismes continus est exact.
\eenum
\ethm
En particulier, il résulte de \ref{adj-lin} que l'on peut également voir les $A$-modules de Tate comme des $A$-pro-modules.
Si $M$ est un $A$-module de Tate, on dispose de son dual $M^{\vee}=\Hom_{c}(M,A)$ où la base de la topologie est formée par les complémentaires orthogonaux des sous-modules bornés ouverts $L\subset M$.
Alors $M^{\vee}$ est également un $A$-module de Tate et l'on a $M\cong M^{\vee\vee}$ (\cite[3.2.4]{Dr}).

\subsection{Pro-modules pro-projectifs}
\bprop\label{def-pro-proj}
Soit un anneau commutatif $A$ et $M=(M_{\al})$ un $A$-pro-module, alors les assertions suivantes sont équivalentes:
\benumr
\remi
Le foncteur qui associe à tout $A$-module $N$ le groupe:
\[\Hom_c(M,N)=\colim_{\al} Hom_{A}(M_{\al},N)\]
est exact.
\remi
Pour tout $\al$, il existe $\beta\geq\al$ tel que le morphisme de transition $M_{\beta}\ra M_{\al}$ admet une factorisation $M_{\beta}\ra P\ra M_{\al}$ avec $P$ projectif (ou de manière équivalente $P$ libre).
\remi
Le pro-module $M$ est isomorphe à un système projectif filtrant de $A$-modules projectifs.
\eenum
Dans ce cas, un pro-module $M$ qui vérifie une de ces assertions est dit pro-projectif.
\eprop
\bpf
Les implications $(ii)\Longrightarrow (iii)$ et $(iii)\Longrightarrow (i)$ sont claires.
Pour $(i)\Longrightarrow (ii)$, fixons $\al$, on considère  une surjection $A^{(I)}\ra N_{\al}$ et on applique le foncteur $\Hom_c(M,)$ qui est exact et on obtient le résultat voulu.
\epf
\brems
\remi\label{proj-rem1}
Par le théorème de Govorov-Lazard (\cite{Gov},\cite{Laz}), tout $A$-module plat est une colimite filtrante de $A$-modules projectifs de type fini, de telle sorte que le pro-module dual d'un $A$-module plat est pro-projectif de type fini.
\remi\label{proj-rem2}
On a immédiatement à l'aide de la caractérisation (iii) et de l'énoncé analogue pour les $A$-modules que si l'on a une suite exacte:
\[0\ra M\ra P\ra N\ra0\]
de $A$-pro-modules avec $P$ et $N$ pro-projectifs, alors $M$ est pro-projectif.
\erems
Dans la suite, on a besoin d'extraire une certaine classe de pro-modules pro-projectifs. On commence par une classe de pro-modules pro-projectifs de type fini.
\blem\label{flat-ML}
Soit un $A$-pro-module $M$ alors les assertions suivantes sont équivalentes:
\benumr
\remi
$M$ est le dual d'un $A$-module plat de Mittag-Leffler.
\remi
$M$ est un pro-objet dans la catégorie des $A$-modules projectifs de type fini et est un pro-ensemble de Mittag-Leffler.
\remi
$M$ est pro-projectif et peut être représenté comme un système projectif filtrant de $A$-modules de type fini $M_{\al}$ avec des morphismes de transition $M_{\beta}\ra M_{\al}$ surjectifs.
\eenum
\elem
\bpf
$(ii)\Longrightarrow (i)$  et $(ii)\Longrightarrow (iii) $ sont clairs. $(i)\Longrightarrow (ii)$  se déduit de \ref{proj-rem1}.
Montrons $(iii)\Longrightarrow (ii)$, on utilise la caractérisation (ii) de \ref{def-pro-proj}, en particulier comme $M$ est pro-projectif, tout morphisme de transition $M_{\beta}\ra M_{\al}$, peut se factoriser en $M_{\beta}\ra P\ra M_{\al}$ avec $P$ libre, que l'on peut supposer de type fini comme $M_{\beta}$ est de type fini. 
\epf

On va maintenant expliciter le lien entre être de Tate et être pro-projectif. On commence par le cas pro-projectif de type fini, qui est plus simple.
\bprop\label{aleph-tate}
Soit un $A$-pro-module $M$ pro-projectif de type fini, $\aleph_0$ et est un pro-ensemble de Mittag-Leffler, alors c'est le dual d'un module de projectif.
\eprop
\bpf
Comme $M$ est pro-projectif de type fini, il est dual d'un $A$-module plat $N$, comme de plus, c'est un pro-ensemble de Mittag-Leffler, on en déduit que $N$ est plat de Mittag-Leffler, comme $M$ est $\aleph_0$, $N$ l'est aussi et par Raynaud-Gruson \cite[p.73-74]{RG}, cela se déduit du fait qu'un $A$-module plat de Mittag-Leffler dénombrablement engendré est projectif.
\epf
Dans la pratique, cela est trop restrictif de considérer des pro-projectifs de type fini. Typiquement, les pro-modules que l'on veut caractériser sont les $\Omega^{1}_{X}$ pour $X$ un $\aleph_0$-ind-schéma raisonnable formellement lisse. On étudie cela dans la section \ref{mod-mlt}.
\subsection{Modules 2-presque projectifs}
On tire les définitions suivantes de \cite[sect. 4]{Dr}:
\bdefi\label{tate}
Soit $A$ un anneau.
Un $A$-module  est $2$-presque projectif  s'il s'écrit comme facteur direct d'un $A$-module de la forme $P\oplus M$, avec $P$ un $A$-module projectif et $M$ un $A$-module de présentation finie qui sont dits  modules $2$-presque projectifs  élémentaires.
\edefi
\brems
\remi\label{2-quot1}
Le $2$ de 2-presque projectif se rapporte au fait qu'un tel module s'écrit dans la catégorie dérivée $P\oplus M^{\bullet}$ où $M^{\bullet}$ 
est un complexe de $A$-modules de type fini concentré en degrés négatifs et dont les deux premiers crans, sont des modules projectifs de type fini.
\remi\label{2-quot2}
Si $P$ est 2-presque projectif et $N\subset P$ est un sous-module de type fini, alors $M=P/N$ est 2-presque projectif.
En effet, si de plus $P=A^{(I)}$ est libre, on choisit $l$ tel que $N\subset A^{l}\subset A^{(I)}$, de telle sorte que $M\dhra A^{(I)}/A^{(l)}=K$ a un noyau de présentation finie et est scindé comme $K$ est projectif.
Si $P$ est projectif, si $P$ est facteur direct de $A^{(I)}$, $P/N$ est facteur direct de $A^{(I)}/N$.
Enfin, si $P$ est 2-presque projectif, il existe une surjection $\phi:P'\dhra P$ avec $P'$ projectif dont le noyau est de type fini et tel que $\phi^{-1}(N)$ est de type fini, on se ramène alors au cas projectif.
\remi\label{2-quot3}
La notion de 2-presque projectivité est locale pour la topologie fpqc \cite[Thm. 4.2.(i)]{Dr}.
\remi\label{2-quot4}
Pour tout $A$-module 2-presque projectif $M$, il existe un recouvrement Nisnevich $\Spec(A')\ra\Spec(A)$ tel que $M\otimes_{A}A'$ devienne élémentaire \cite[Thm. 4.2.(ii)]{Dr}.
\erems

Le lemme suivant est établi par Drinfeld \cite[Prop. 9.11]{Dr2}. On en donne la preuve pour la commodité du lecteur:
\blem\label{2nis}
Soit $A$ un anneau commutatif, $M$ un $A$-module 2-presque projectif et plat, alors $M$ est projectif.
\elem
\bpf
Soit $M$ un tel module. On utilise l'observation suivante de Govorov-Lazard (\cite{Gov},\cite{Laz}): Tout morphisme de $A$-modules $F\ra M$ avec $F$ de présentation finie et $M$ plat se factorise par un module libre de type fini (\cite[Tag. 058E]{Sta}).
Soit alors $g:P\oplus F\dhra M$ une surjection scindée avec $P$ libre et $F$ de présentation finie. Alors $g_{\vert F}$ s'écrit comme une composée $F\ra A^{b}\ra M$ et donc $g$ peut s'écrire comme une composée $P\oplus F\ra P\oplus A^{n}\ra M$ et la flèche $P\oplus A^{n}\ra M$ est une surjection scindée comme $g$ l'est.
\epf
\subsection{Pro-platitude}
Dans la catégorie $\Pro(\Mod^{\heartsuit}_{A})$, on a une structure tensorielle naturelle:
\[\hat{\otimes}:\Pro(\Mod^{\heartsuit}_{A})\times\Pro(\Mod^{\heartsuit}_{A})\rightarrow\Pro(\Mod^{\heartsuit}_{A}),\]
qui envoie $((M_{i}), (N_{j}))$ sur $(M_{i}\otimes_{A} N_{j})$.
La notation se justifie par le fait que si $M, N \in\TopMod_{A}$ sont des $A$-modules topologiques, complets, Hausdorff avec une topologie linéarisable, alors:
\[M\hat{\otimes}_{A}N\simeq F(G(M)\hat{\otimes} G(N)).\]

\begin{defi}
Soit $A$ un anneau,	un $A$-pro-module $M$ est dit topologiquement plat si le foncteur:
\[N\mapsto N\hat{\otimes}_{A}M\]
est exact dans la catégorie des $A$-pro-modules.
\end{defi}
\brems
\remi
Comme toute suite exacte de pro-modules est une système projectif cofiltrant de suites exactes et que les limites cofiltrantes  sont exactes, on en déduit que $M\in\Pro(\Mod^{\heartsuit}_{A})$ est topologiquement plat si le foncteur $M\hat{\otimes}_{A}-$ est exact dans la catégorie des $A$-modules.
\remi
On a  immédiatement que si $M=(M_{i})_{i\in\cI}$ est pro-plat, i.e. un système projectif filtrant de $A$-modules plat, alors $M$ est topologiquement plat.
L'étude de la réciproque est une question délicate.
\remi\label{pro-plat}
Si $M$ est un $A$-module de Tate, alors il est topologiquement plat puisque d'après \ref{adj-lin}, \ref{Tate-prop} et \ref{def-pro-proj}, vu comme pro-module, il est pro-projectif.
\erems

\begin{prop}\label{pro-pl1}
Soit $M:=(M_{i})\in\Pro(\Mod^{\heartsuit}_{A})$.
Les assertions suivantes sont équivalentes:
\benumr
	\item 
	$M$ est topologiquement plat.
	\item
	Pour tout $A$-module $N$, $(\Tor_{1}^{A}(N,M_{i}))$ est essentiellement nul.
	\item
	Pour tout $i$, il existe $j\geq i$, tel que le morphisme $M_{j}\rightarrow M_{i}$ induise 0 sur $\Tor^{1}$ pour tout $A$-module de présentation finie $N$.
	\item
	Pour tout idéal $I$ de type fini  de $A$, $(\Tor_{1}^{A}(A/I,M_{i}))$ est essentiellement nul.
	\eenum
	\end{prop}
	\bpf
	L'équivalence entre les deux premières assertions est évidente, en considérant une suite exacte de $A$-modules:
$$\xymatrix{0\ar[r]&K\ar[r]&P\ar[r]&N\ar[r]&0}$$
avec $P$ projectif.
Pour la troisième assertion, il suffit de prendre $N=\oplus N_{\la}$ où $N_{\la}$ parcourt tous les $A$-modules de présentation finie à isomorphismes près.
Pour la dernière, cela utilise le même dévissage que \cite[Tag. 00HD]{Sta}.
	\epf
	Le lemme suivant,  tiré de Gabber-Ramero \cite[2.4.17]{GR} permet d'obtenir une caractérisation de la platitude topologique en termes de systèmes projectifs de $A$-modules plats:
\blem\label{factgr}
Soit $A$ un anneau, $M$ un $A$-module et $C:=\coker(\phi:A^{n}\rightarrow A^{m})$ un $A$-module de présentation finie quelconque. Soit $C':=\coker(\phi^{*}:A^{m}\rightarrow A^{n})$ le conoyau de la transposée de $\phi$, alors il y a un isomorphisme naturel:
\[\Tor_{1}^{A}(C',M)\simeq\Hom_{A}(C,M)/\Ima(\Hom_{A}(C,A)\otimes_{A}M).
\]
\elem
On en déduit alors:
\bprop\label{pro-pl2}
Soit $M=(M_{i})\in\Pro(\Mod^{\heartsuit}_{A})$, alors $M$ est topologiquement plat si et seulement si pour tout $i$, il existe $j\geq i$ tel que pour tout morphisme $N\ra M_j$ avec $N$ de présentation finie, la composée,
$N\ra M_i$ se factorise par un $A$-module libre de type fini.

En particulier, si $M\in\Pro(\Mod^{\heartsuit,pf}_{A})$, $M$ est équivalent à un système projectif filtrant de $A$-modules libres de type fini, donc en particulier pro-plat. Ici, $\Mod^{\heartsuit, pf}_{A}$ désigne la catégorie des $A$-modules de présentation finie.
\eprop
\begin{proof}
Soit $N$ un $A$-module de présentation finie et un morphisme:
\[\theta:N\rightarrow M_{j}.\]
D'après le lemme \ref{factgr} et comme la flèche $M_{j}\rightarrow M_{i}$ induit 0 sur le $\Tor_{1}$ pour tout $A$ module de présentation finie d'après \ref{pro-pl1}, en l'appliquant à $N'$, le composé $\theta':N\rightarrow M_{i}$ est l'image d'un certain élément:
\[\sum\limits_{i=1}^{p}\phi_{i}\otimes m_{i}\in \Hom_{A}(N,A)\otimes_{A}M_{i}.\]
En considérant $L=A^{p}$, le morphisme $\theta'$ se factorise ainsi en:
\[N\stackrel{u}{\rightarrow} A^{p}\stackrel{v}{\rightarrow} M_{i},\]
avec $u(x):=(\phi_{1}(x),\dots,\phi_{p}(x))$ et $v(y_{1},\dots,y_{p})=\sum\limits_{i=1}^{p} y_{i}m_{i}$.
\end{proof}

Ce critère va nous permettre de déterminer un autre contexte où un pro-module topologiquement plat est bien pro-plat.
\bprop\label{pro-pl3}
Soit $M:=(M_{i})\in\Pro(\Mod^{\heartsuit}_{A})$, on suppose $M$ topologiquement plat et les $M_i$ 2-presque projectifs, alors $M$ est pro-projectif, donc pro-plat.
\eprop
\bpf
Fixons $i$, en appliquant \ref{pro-pl1}, il existe $j$ tel que le morphisme $M_j\ra M_i$ induise l'application nulle sur $\Tor_{1}$ pour tout $A$-module de présentation finie.
 Comme chaque $M_j$ est 2-presque projectif, c'est un facteur direct d'un $A$-module de la forme $P\oplus N_{j}$ avec $P$ projectif et $N_j$ de présentation finie.
On a alors une factorisation:
$$\xymatrix{M_{j}\ar@/^1pc/[rr]^{\Id}\ar[r]&P\oplus N_j\ar[r]&M_{j}\ar[r]&M_i}$$
et il résulte de \ref{pro-pl2} que la flèche $N_j\ra M_j\ra M_i$ se factorise par un module libre de type fini $Q$. On en déduit alors que le morphisme $M_j\ra M_i$ se factorise par $P\oplus Q$, ce qui conclut.
\epf

On a besoin d'une version de Nakayama pour les pro-modules:
\blem\label{Nak1}
Soit $A$ un anneau local d'idéal maximal $\km$, $k$ son corps résiduel, $M, N \in \Pro(\Mod_{A})$ et un morphisme de $A$-pro-modules $f:N\ra M$.

\benumr
\remi
Supposons $M, N$ dans $\Pro(\Mod^{\heartsuit,pf}_{A})$, $f$ surjective et $f\otimes_{A}k$ injective, alors $f$ est un isomorphisme.
\remi
Supposons $M$ de Tate et  $N$ le dual d'un module projectif et  $f\otimes_{A}k$ injective, alors $f$ est injective et $M/f(N)$ est topologiquement plat.
\eenum

\elem
\bpf
(i) D'après \cite[Tag. 0519]{Sta}, si $K=\Ker f$, alors $K\in \Pro(\Mod^{\heartsuit,tf}_{A})$ \footnote{Ici $\Mod^{\heartsuit,tf}_{A}$ est la catégorie des $A$-modules de type fini.} et par hypothèse $K\otimes_{A}k=0$. En particulier si $K=(K_{\al})$ est une présentation comme système projectif de $A$-modules de type fini, pour tout $\al$, il existe $\beta$ tel que $f_{\beta\al}:K_{\beta}\ra K_{\al}$ soit l'application nulle modulo $\km$. En appliquant Nakayama à $\Ima f_{\beta\al}$ qui est de type fini comme $K_{\beta}$ l'est, on en déduit que $f_{\beta\al}=0$ et $K=0$.

(ii) Comme $M$ est de Tate c'est un facteur direct d'un $A$-module de la forme $P\oplus Q^{\vee}$ pour $P, Q$ projectifs. En particulier, comme $N$ est le dual d'un module projectif, la flèche induite:
\[N\stackrel{f}{\rightarrow} M\ra P\oplus Q^{\vee} \ra P\]
 a une image de type fini (cf. \cite[Lem. 3.1]{Dr}).
Il existe donc un $A$-module projectif $H$ et une injection $i: H^{\vee}\ra M$ et tel que $f$ se factorise par $H^{\vee}\in \Pro(\Mod_{A}^{\heartsuit,pf})$. On en déduit alors d'après (i) que $N\cong H^{\vee}$ et $f$ injective.
Pour la platitude, il suffit de vérifier d'après \ref{pro-pl1}(iv) que pour tout idéal $I$ de type fini, $\widehat{Tor}^{1}_{A}(M/f(N), A/I)=0$ et comme $M$ est de Tate donc topologiquement plat, cela revient à montrer l'injectivité de :
\[\bar{f}:N/IM\ra M/IM.\]
Or $M/IM$ reste Tate, $N/IM$ reste le dual d'un $A/I$-module projectif et $\bar{f}$ est injectif modulo $\km$, il suffit donc d'appliquer la première partie de l'argument pour obtenir l'injectivité, d'où $M/f(N)$ topologiquement plat.
\epf

\section{Appendice B : Autour de certains résultats sur les ind-schémas}\label{app-drin}
L'objet de cet appendice est de donner la preuve de plusieurs résultats généraux sur les ind-schémas, notamment \ref{ind-Tate} et \ref{gab-fl}.

\subsection{Retour sur les Modules de Mittag-Leffler}
Pour le moment, on a seulement défini les modules plats de Mittag-Leffler, il est commode de faire une définition plus générale qui ne fait pas intervenir la platitude.
La définition suivante est tirée de Raynaud-Gruson \cite[2.1.3, p. 69]{RG}.
\bdefi\label{ML-def}
Soit un anneau commutatif $A$.
Soit un $A$-module $M$, alors il est de Mittag-Leffler si pour tout $A$-module de présentation finie $F$, tout morphisme $f:F\ra M$ se factorise en $F\stackrel{u}{\rightarrow}G\ra M$ avec $G$ de présentation finie tel que pour toute autre factorisation 
$F\stackrel{u'}{\rightarrow} G'\ra M$ avec $G'$ de présentation finie, on ait un morphisme $\phi:G'\ra G$ tel que $u=\phi u'$.
\edefi
\brems
\remi\label{def-eq1}
La définition ci-dessus admet la reformulation suivante; si $\Mod_{A}$ est la catégorie des $A$-modules, alors $\Mod_{A}=\Ind(\Mod_{A}^{pf})$ où $\Mod_{A}^{pf}$ est la catégorie des $A$-modules de présentation finie.
On considère alors $M\in\Mod_{A}$ que l'on écrit comme une colimite filtrante $M\simeq\varinjlim M_{i}$, alors d'après \cite[2.1.4, p. 69]{RG}, $M$ est de Mittag-Leffler si et seulement si pour tout $A$-module $N$, $(\Hom_{A}(M_{i},N))$ est un pro-ensemble de Mittag-Leffler.
\remi\label{def-eq2}
De plus, si $M$ est un $A$-module plat, alors de Mittag-Leffler si et seulement si $(\Hom_{A}(M_{i},A))$ est un pro-ensemble de Mittag-Leffler, ce qui fait le lien avec la définition \ref{pro-ml}. Un module projectif est plat de Mittag-Leffler \cite[Tag. 059Z]{Sta}.
\remi\label{def-eq3}
La condition de Mittag-Leffler satisfait la descente fidèlement plate \cite[Tag. 05A5]{Sta}.
\erems
Cette notion se faisceautise de telle sorte que si $X$ est un schéma et $\cF$ un faisceau quasi-cohérent sur $X$, on peut parler de faisceau de Mittag-Leffler. Un exemple typique qui fait apparaître des faisceaux de Mittag-Leffler est le suivant:
\bprop\label{ml2}
Soit un morphisme $f:X\ra Y$ entre $k$-espaces algébriques qcqs avec une présentation placide, alors:
\benumr
\item
 $\Omega^{1}_{X/Y}$ est une extension d'un faisceau quasi-cohérent plat de Mittag-Leffler et d'un faisceau cohérent de présentation finie.
\item
Si de plus, $f$ est fortement pro-lisse $\Omega^{1}_{X/Y}$ est plat de Mittag-Leffler. 
\item
Si $f$ est formellement lisse, alors $f$ admet une présentation fortement pro-lisse $X\simeq\varprojlim X_{\al}$ telle que pour tout $\al$, $X\ra X_{\al}$ est formellement lisse.
\eenum
\eprop
\bpf
On montre (i) et (ii) simultanément. On considère une présentation placide $X\simeq\varprojlim X_{\al}$ de $f$. Fixons un $\al$ de telle sorte que $f_{\al}:X\ra X_{\al}$ est fortement pro-lisse et $L_{X/X_{\al}}$ est plat concentré en degré zéro.
On a un triangle distingué:
\[f_{\al}^{*}L_{X_{\al}/Y}\ra L_{X/Y}\ra L_{X/X_{\al}},\]
soit en  prenant la suite exacte longue de cohomologie:
\[0\ra f^{*}_{\al}\Omega^{1}_{X_{\al}/Y}\ra \Omega^{1}_{X/Y}\ra \Omega^{1}_{X/X_{\al}}\ra 0\]
Le faisceau $f^{*}_{\al}\Omega^{1}_{X_{\al}/Y}$ est de présentation finie sur $\cO_{X}$, reste à montrer que  $\Omega^{1}_{X/X_{\al}}$ est de Mittag-Leffler, ce qui montrera également l'assertion pour fortement pro-lisse puique que l'on prendra $X_{\al}=Y$. D'après \ref{def-eq2}, il s'agit de voir que $(f^{*}_{\beta}\Theta_{X_{\beta}/Y})_{\beta\geq\al}$ est un pro-ensemble de Mittag-Leffler  où $\Theta_{X_{\beta}/Y}$ est le dual de $\Omega_{X_{\beta}/Y}$. Soit $\beta\geq\al$, alors $f_{\al\beta}:X_{\beta}\ra X_{\al}$ est lisse d'où une surjection: 
\[\Theta_{X_{\beta}/Y}\dhra f_{\al\beta}^{*}\Theta_{X_{\al}/Y},\]
ce qui conclut.

(iii) D'après \ref{plac}, on sait déjà que $f$ admet une présentation fortement pro-lisse, il s'agit de vérifier que pour tout $\al$, $f_{\al}:X\ra X_{\al}$ est formellement lisse.
On sait déjà que $L_{X/X_{\al}}$ est plat concentré en degré zéro et quasi-isomorphe à $\Omega^{1}_{X/X_{\al}}$.
On en déduit un isomorphisme canonique:
\[\Omega^{1}_{X/Y}/f_{\al}^{*}\Omega^{1}_{X_{\al}/Y}\cong\Omega^{1}_{X/X_{\al}}.\]
Or, comme $f$ est formellement lisse  et $X_{\al}\ra Y$ est de présentation finie, d'après \cite[Tag. 0D0L]{Sta}, $\Omega^{1}_{X/Y}$ est localement projectif et donc d'après \ref{2-quot2}, $\Omega^{1}_{X/X_{\al}}$ est localement 2-presque projectif. Comme il est aussi plat, il est localement projectif d'après \ref{2nis} et $f_{\al}$ est formellement lisse d'après \cite[Tag. 0D10]{Sta}.
\epf

\subsubsection{Modules de Mittag-Leffler-Tate}\label{mod-mlt}
La section suivante est entièrement due à Drinfeld.
\bdefi
Un $A$-pro-module $M$ est un module de Mittag-Leffler-Tate élémentaire s'il s'inscrit dans la catégorie des $A$-pro-modules dans une suite exacte:
\begin{equation}
0\ra Q^{\vee}\ra M\ra P\ra 0,
\label{exa-ml}
\end{equation}
où $Q$ est plat de Mittag-Leffler et $P$ est projectif.
Un $A$-pro-module $M$ est un module de Mittag-Leffler-Tate si c'est un facteur direct d'un module de Mittag-Leffler-Tate élémentaire.
\edefi
\brem\label{split-ml}
Si $P$ et $Q$ sont des colimites filtrantes de modules projectifs de type fini $P_{\al}$ et $Q_i$, alors la suite exacte \eqref{exa-ml} correspond à un système compatible de suites exactes:
\[0\ra Q_i^{\vee}\ra M\ra P_{\al}\ra 0.\]
On en déduit donc que $\Ext(P,Q^{\vee})=\Ext(P\otimes_{A}Q,A)=0$ si $P$ et $Q$ sont tous deux projectifs.
En revanche si $Q$ est de Mittag-Leffler, on a des contre-exemples même pour $A=\bZ$.
\erem
On commence par rappeler un énoncé fondamental sur les modules plats de Mittag-Leffler \cite[p. 73-74]{RG}:
\bthm\label{RGcont}[Raynaud-Gruson]
Les assertions suivantes sont équivalentes:
\benumr
	\remi
	$M$ est un $A$-module plat de Mittag-Leffler.
	\remi
	Pour tout sous-ensemble dénombrable $S\subset M$, il existe un sous-module projectif $Q\subset M$, dénombrablement engendré tel que $M/Q$ est plat.
\eenum
\ethm

On a la propriété suivante:
\bprop\label{mlt-al}
Un $\aleph_0$-$A$-module de Mittag-Leffler-Tate est de Tate.
\eprop
\bpf
On considère une suite exacte telle que dans \eqref{exa-ml}. Supposons qu'un facteur $M'$ de $M$ est $\aleph_0$, alors la composée $Q^{\vee}\hra M\dhra M'$ se factorise en $Q_0^{\vee}\hra M\dhra M'$
pour un sous-module dénombrablement engendré $Q_0\subset Q$ et il résulte de \ref{RGcont} que l'on peut supposer que $Q_0$ est projectif.
On considère alors la suite, obtenue en poussant-en-avant \eqref{exa-ml} par le morphisme $Q^{\vee}\ra Q_0^{\vee}$:
\begin{equation}
0\ra Q_{0}^{\vee}\ra M_0\ra P\ra 0.
\label{exa-ml2}
\end{equation}
On a des morphismes canoniques $M'\ra M_0\ra M'$, dont la composée est $\Id_{M'}$, donc $M'$ est un facteur direct de $M_0$. Maintenant, d'après \ref{split-ml}, la suite \eqref{exa-ml2} se scinde et donc $M_0$ est de Tate, ainsi que $M'$.
\epf

On peut maintenant formuler le théorème principal de cette section:
\bthm\label{car-mlt}[Drinfeld]
Soit un $A$-pro-module $M$ pro-projectif, dénombrablement engendré, qui peut être représenté par un système projectif filtrant de $A$-modules $M_{\al}$ avec des morphismes de transition surjectifs $f_{\beta\al}:M_{\beta}\ra M_{\al}$ dont les noyaux sont de type fini, alors $M$ est de Tate.
\ethm
Avant de montrer le théorème, on commence par deux résultats préliminaires.
\bprop\label{6.6}
Soit $M$ un $A$-pro-module qui vérifie les hypothèses de \ref{car-mlt}, alors les $M_{\al}$ sont 2-presque projectifs.
\eprop
\bpf
Comme $M$ est pro-projectif, pour toute paire $(\al,\beta)$, il existe $P$ libre tel que $f_{\beta\al}$ se factorise en $M_{\beta}\stackrel{g}{\rightarrow} P\ra M_{\al}$, soit $N=\Ker(f_{\beta\al})$, l'isomorphisme canonique $M_{\beta}/N\cong M_{\al}$ admet une factorisation $M_{\beta}/N\ra P/g(N)\ra M_{\al}$, donc $M_{\al}$ est un facteur direct de $P/g(N)$, lequel est une somme directe d'un module de présentation finie et d'un module projectif d'après \ref{2-quot2}.
\epf
\blem\label{6.5}
Soit $M_0$ un $A$-module 2-presque projectif, alors il existe une suite de surjections:
\begin{equation}
\dots\stackrel{\phi_3}{\rightarrow}M_{2}\stackrel{\phi_2}{\rightarrow}M_1\stackrel{\phi_1}{\rightarrow}M_0.
\label{16.2}
\end{equation}
telle que pour tout $i$, $\Ker(\phi_i)$ est contenu dans un sous-module de type fini de $M_i$ et $\phi_i$ se factorise en $M_i\ra P_i\ra M_{i-1}$ avec $P_i$ projectif.
\elem
\bpf
Il suffit de construire $(M_1,\phi_1)$. Soit $P$ un $A$-module projectif et $N$ un $A$-module de présentation finie, tel que $M_0$ soit un facteur direct de $P\oplus N$. On peut écrire que $M_0=\Ima(\pi)$ pour $\pi\in\End(P\oplus N)$ avec $\pi^2=\pi$.
Soit $\phi:M_0\ra M_0$ obtenue comme la composée:
\[M_0\hra P\oplus N\ra P\stackrel{\pi_{\vert P}}{\rightarrow} M_0\]
 Alors, $\Ima(\phi-\Id)$ est contenu dans un sous-module de type fini $L\subset M_0$ et $\phi$ admet une factorisation $M_0\ra P\ra M_0$ avec $P$ projectif. On choisit une surjection $F\dhra L$ avec $F$ projectif de type fini.
On pose alors $M_1=M_0\oplus F$ et on définit $\phi_1:M_1=M_0\oplus F\ra M_0$ avec $\phi_{1}(m,f)=\phi(m)+\pi(f)$. Comme $\Ima(\phi-\Id)\subset L=\Ima(\pi)$, on voit que $\phi_1$ est surjective.
Si $\phi_{1}(m,f)=0$, alors:
\[m=(m-\phi(m))-\pi(f)\in\Ima(\phi-\\Id)+L\subset L,\]
 donc $\Ker\phi_{1}\subset L\oplus F$, comme voulu.
\epf
On passe maintenant à la preuve de \ref{car-mlt}:
\bpf
Soit $M$ un $A$-pro-module qui vérifie les hypothèses du théorème, on va construire $M'$ tel que $M\oplus M'$ soit de Mittag-Leffler-Tate élémentaire. Puis il suffit d'utiliser \ref{mlt-al} pour en déduire que $M$ est de Tate.
On choisit alors un indice $\al_1$, d'après \ref{6.6}, $M_{\al_1}$ est presque projectif donc il existe $M_0$ tel que $M_{\al_1}\oplus M_0$ soit la somme directe d'un $A$-module projectif $P$ et d'un $A$-module de type fini. En appliquant \ref{6.5} à $M_0$, on obtient un système projectif qui vérifie \eqref{16.2}. Soit $M'$ le pro-module correspondant, il reste à voir que $M\oplus M'$ est de Mittag-Leffler-Tate élémentaire. Or, le noyau de la composée $M\oplus M'\ra M_{\al_{1}}\oplus M_0\ra P$ satisfait \ref{flat-ML}.(iii), donc est dual d'un $A$-module plat de Mittag-Leffler.
\epf
On a alors l'application immédiate suivante:
\bthm\label{ind-Tate}
Soit $X$ un $k$-$\aleph_0$-ind-schéma ind-placide formellement lisse, alors $\Omega^{1}_{X/k}$ est de Tate.
\ethm
\bpf
D'après \ref{fl-gr} et \ref{etfl}, $\Omega^{1}_{X/k}$ est pro-projectif, comme $X$ est $\aleph_0$ et ind-placide, $\Omega^{1}_{X/k}$ est bien de Mittag-Leffler et les morphismes de transition ont des noyaux qui sont des faisceaux cohérents de présentation finie, en particulier, d'après \ref{car-mlt}, $\Omega^{1}_{X/k}$ est de Mittag-Leffler-Tate et comme il est $\aleph_0$, il est de Tate d'après \ref{mlt-al}.
\epf

\subsection{Ind-schémas formellement lisses}

On commence par étendre la classe des ind-schémas. 

\bdefi
Un \og ind-schéma\fg~ est un foncteur $X:(\Aff_{\bZ})^{op}\ra\Ens$ qui admet une présentation comme colimite filtrante $X\simeq\varinjlim X_i$ où les $X_i$ sont des schémas et les flèches de transition sont des immersions fermées.
On dit que $X$ est un \og schéma formel\fg~ si $X_{red}$ est un schéma. Il est dit raisonnable, si de plus les immersions fermées sont de présentation finie.
\edefi
\brem
La présence des guillemets vient du fait que l'on n'impose pas que les $X_i$ soient qcqs.
\erem
 D'après \cite[Tag. 0738]{Sta}, pour tout schéma qcqs $T$, on a 
\begin{equation}
\Hom(T,X)\simeq\colim\Hom(T,X_i).
\label{calco}
\end{equation}
 L'énoncé suivant a été noté par Richarz \cite[Lem. 1.4]{Ric}:
\blem\label{fpqc}
Soit un \og ind-schéma\fg~ $X$, pour tout schéma $T$ et un recouvrement fpqc $(T_i\ra T)_{i\in T}$, la suite d'ensembles:
\[\Hom(T,X)\ra\prod\limits_{i\in I}\Hom(T_i,X)\ra\prod\limits_{i,j\in I}\Hom(T_i\times_{T}T_{j},X),\]
est exacte, en particulier  $X$ est un faisceau fpqc.
\elem
\bpf
Comme $\Hom(-,X)$ transforme les colimites en limites et que $T\cong\varinjlim_{U\subset T} U$ où $U$ parcourt les ouverts affines de $T$, on se ramène d'abord au cas où $T$ affine et ensuite, en raffinant le recouvrement, au cas où les $T_i$ sont aussi affines et en nombre fini.
On applique ensuite \eqref{calco} et le fait que les colimites filtrantes commutent aux limites finies, pour se ramener au cas d'un schéma et cela se déduit de \cite[Tag. 023Q]{Sta}.
\epf

On rappelle qu'un faisceau  quasi-cohérent $\cF$ sur un schéma $Y$ est localement 2-presque projectif si pour tout ouvert affine $U\subset Y$, $\cF(U)$ est un $\Gm(U,\cO_{U})$-module 2-presque projectif.
On a le théorème suivant dû à Gabber et à Drinfeld dans le cas formellement lisse standard \cite[Thm. 6.2]{Dr}:
\bthm\label{fond-tate}
Soit un morphisme de \og ind-schémas\fg~ formellement lisse $X\ra S$ avec $S$ un schéma et $X$ ind-placide, alors pour tout sous-schéma fermé de type fini $Y\subset X$, $\Omega^{1}_{Y/S}$ est localement 2-presque projectif.
\ethm
\bpf
La notion de formelle lissité ne dépend que de la complétion formelle de $Y$ le long $X$, de telle sorte que l'on peut supposer que $X$ est un \og schéma formel\fg~ ind-placide.
D'après \ref{2-quot3}, il suffit de vérifier l'assertion pour tout ouvert affine $U\subset Y$. On peut donc supposer que $Y$ est affine et donc $X$ ind-affine d'après \cite[Tag. 0DE8]{Sta}. 
On considère une résolution:
\begin{equation}
0\ra \cI\ra F\ra\Omega^{1}_{Y/S}\ra 0, 
\label{Cplx}
\end{equation}
où $F$ est un $\cO_{Y}$-module libre.
D'après \cite[Ch.III, Thm. 1.2.3]{Ill}, on a un isomorphisme:
\[\Ext^{1}_{\cO_{Y}}(L_{Y/S},\cI)\cong\Exalcom_{\cO_{S}}(\cO_{Y},\cI),\]
où $\Exalcom_{\cO_{S}}(\cO_{Y},\cI)$ consiste en les extensions de $\cO_{Y}$ par $\cI$ qui sont des $\cO_{S}$-algèbres avec $\cI^2=0$.
En utilisant l'augmentation $L_{Y/S}\ra\Omega^{1}_{Y/S}$, l'extension \eqref{Cplx} fournit une classe dans $\Ext^{1}_{\cO_{Y}}(L_{Y/S},\cI)$, soit une extension infinitésimale $Y=\Spec(A_0)\hra Z=\Spec(A)$. Par formelle lissité de $X$, il existe $Z'=\Spec(B)\ra Z$ affine fidèlement plat et un fermé de type fini $Y\subset Y'\subset X$ tel qu'on ait un diagramme commutatif:
$$\xymatrix{\Spec(B_0)=Y\times_{Z}Z'\ar[r]\ar[d]&Z'\ar[d]\\Y\ar[r]&Y'}$$
Soit $I=\Ker(A\dhra A_0)$, comme tous les schémas sont affines, on forme $W$ le pushout qui s'inscrit dans le diagramme:
$$\xymatrix{Y\times_{Z}Z'\ar[r]\ar[d]&Z'\ar[d]\\Y'\ar[r]&W},$$
qui est une extension de $\cO_{Y'}$ par $I\otimes_{A}B=I\otimes_{A_0}B_0$.
Maintenant, l'existence de $Z\ra Y'$ correspond au fait que la classe dans $\Ext^{1}(\Omega^{1}_{Y'/S},I\otimes_{A_0}B_0)$ qui provient de \eqref{Cplx} est triviale, donc admet un scindage, soit comme $A_0\ra B_0$ est plat:
$$\xymatrix{0\ar[r]& I\otimes_{A_0}B_0\ar[r]&F\otimes_{A_0}B_0\ar[r]&\Omega^{1}_{Y/S}\otimes_{A_0}B_0\ar[r]&0\\&&&(\Omega^{1}_{Y'/S}\otimes_{\cO_{Y'}}\cO_{Y})\otimes_{A_0}B_0\ar[u]^{\phi}\ar@{.>}[ul]^\psi }$$
Comme $\Ker\phi$ est de type fini, d'après \ref{2-quot2},
$(F\otimes_{A_0}B_{0})/\psi(\Ker\phi)$ est 2-presque projectif et 
$\Omega^{1}_{Y/S}\otimes_{A_0}B_0$ en est un facteur direct, donc est aussi 2-presque projectif. On conclut alors par descente à l'aide de \ref{2-quot3}.
\epf

\bthm[Gabber]\label{gab-fl}
Soit un morphisme de \og ind-schémas\fg~ formellement lisse $X\ra S$ avec un schéma $S$  et $X$ ind-placide, alors il est formellement lisse standard.
\ethm
\bpf
La notion de formelle lissité pour $X$ est équivalente à la formelle lissité pour tous les complétés formels le long d'un sous-schéma fermé de type fini de $X$. En particulier, on peut supposer que $X$ est un \og schéma formel\fg~ raisonnable.
Soient $T_0=\Spec(B/J)\ra T=\Spec(B)$, pour un idéal $J$ de carré nul et $u:T_0\ra X$.
Comme $X$ est formellement lisse, il existe $T'=\Spec(B')\ra T$ fidèlement plat tel que $T'\times_{T} T_0\ra T_0\stackrel{u}{\rightarrow} X$ se relève en $\ti{u}:T'\ra X$. D'après \eqref{calco}, soit $Y$ un sous-schéma fermé placide de $X$ par lequel $u$ et $\ti{u}$ se factorisent.
On considère alors l'ensemble des extensions de $u:T\ra Y$. C'est un pseudo-torseur sous $\cH=\cHom(u^{*}\Omega^{1}_{Y/S},\cJ)$ d'après \cite[Tag. 061D]{Sta} où $\cJ$ est le faisceau d'idéaux associé à $J$. Il admet une section pour la topologie fpqc, donc c'est un torseur fpqc. Il s'agit de montrer qu'il est trivial; il suffit de voir que $H^{1}_{fpqc}(T,\cH)=0$ (au sens de Cech). D'après \ref{fond-tate}, $u^{*}\Omega^{1}_{Y/S}$ est 2-presque projectif. L'énoncé se déduit alors de \ref{1}.
\epf
\blem\label{1}
Soit un schéma affine $Z=\Spec(R)$, $\cF_1$, $\cF_2$ des faisceaux quasi-cohérents, supposons $M_1:=H^{0}(Z,\cF_1)$ est 2-presque projectif, alors $H_{fpqc}^{1}(Z,\cHom(\cF_1,\cF_2))=0$.
\elem
\bpf
Comme $M_1$ est 2-presque projectif, $\cF_{1,Zar}$ est un facteur direct d'un faisceau somme directe d'un $\cO_{Z_{zar}}$-module libre et d'un $\cO_{Z_{zar}}$-module cohérent de présentation finie. D'après \cite[Tag. 03OJ]{Sta}, c'est aussi le cas comme $\cO_{Z_{fpqc}}$-module.
Il suffit donc de traiter séparément le cas d'un $\cO_{Z}$-module libre et d'un $\cO_{Z}$-module cohérent de présentation finie.
Si $\cF_1$ est de présentation finie, d'après \cite[Tag. 03M1]{Sta}, $\cHom(\cF_1,\cF_2)$ est quasi-cohérent et d'après \cite[Tag. 03DW]{Sta}, on a $H^{1}_{fpqc}(Z,\cHom(\cF_1,\cF_2))=0$. Si $\cF_1=\cO_{Z}^{(I)}$ est libre, $\cHom(\cF_{1},\cF_{2})=\cF_{2}^{I}$
et d'après \cite[Tag. 060L]{Sta} $H^{1}(Z,\cF_2^{I})\hra \prod\limits_{i\in I} H^{1}(Z,\cF_2)$ et le membre est nul à nouveau d'après \cite[Tag. 03DW]{Sta}.  
\epf

\address{ 
  \bigskip
  \footnotesize

  (A.\ Bouthier) \textsc{IMJ-PRG, Sorbonne Université, 4 Place Jussieu, Paris, 75005 France}\par\nopagebreak
  \textit{E-mail address}: \texttt{alexis.bouthier@imj-prg.fr}}

\begin{thebibliography}{}
\bibitem{AGKRRV}
D. Arinkin, D. Gaitsgory, D. Kazhdan, S. Raskin, N. Rozenblyum, Y. Varshavsky.
\newblock
The stack of local systems with restricted variation and geometric Langlands theory with nilpotent singular support, arXiv:2010.01906.

\bibitem{Bart}
O. Bartlett.
\newblock The singular support of an $\ell$-adic sheaf, https://arxiv.org/pdf/2309.02587.pdf.

\bibitem{Bei}
A. Beilinson.
\newblock
Constructible sheaves are holonomic.
\newblock
Selecta Mathematica, vol. 22,pp. 1797-1819 (2016).

\bibitem{Bor}
Armand Borel.
\newblock Linear algebraic groups 2nd ed.
\newblock Graduate Texts in Mathematics, vol. 126, Springer Verlag, New York, 1991.

\bibitem{BBki}
N. Bourbaki.
\newblock \'{E}léments de mathématique: Algèbre, Chapitres I à III.
\newblock Berlin, Hermann, 1970 (réimpr. 2007).

\bibitem{B}
A. Bouthier.
\newblock 
Faisceaux caractères sur les espaces de lacets d'algèbres de Lie, https://arxiv.org/abs/2108.12663.


\bibitem{BC}
A. Bouthier, K. Cesnavicius.
\newblock Torsors on loop groups and the Hitchin fibration.
\newblock{\em Ann. Sci. de l'ENS}, Tome 55, n°3, pp. 791-864, 2022.


\bibitem{BKV}
A. Bouthier, D. Kazhdan, Y. Varshavsky.
\newblock
Perverse sheaves on infinite-dimensional stacks, and affine
Springer Theory.
\newblock{\em Adv in Maths,} Vol 408, Part A, Oct. 2022, 

\bibitem{SGA4-1/2}
P. Deligne.
\newblock Cohomologie étale SGA 4 1/2.
\newblock Lect. Notes in Math. 569, Springer-Verlag, pp. 233-251, 1977.

\bibitem{Dr}
V. Drinfeld.
\newblock
Infinite-dimensional vector bundles in algebraic geometry: an introduction.
\newblock The Unity of Mathematics, pp. 263-304, Birkhäuser, Boston 2006.


\bibitem{Dr2}
V. Drinfeld.
\newblock
Infinite-dimensional vector bundles in algebraic geometry, version étendue. Notes personnelles.


\bibitem{Ev-Mir}
S. Evens, I. Mirkovic.
\newblock Characteristic cycles for the loop grassmannian and nilpotent orbits.
\newblock{\em Duke Math. J.}, vol. 97, n°1, pp. 109-126 1999.

\bibitem{GaR}
O. Gabber, L. Ramero.
\newblock Almost ring theory.
\newblock Lecture Notes in Mathematics, vol. 1800, Springer-Verlag, Berlin, 2003.

\bibitem{Gai}
D. Gaitsgory.
\newblock Space of rational maps.
\newblock{\em Inventiones Math.}, vol. 191, pp. 91-196, 2013.

\bibitem{GR}
D. Gaitsgory, N. Rozenblyum.
\newblock
DG-indschemes.
\newblock Contemporary Mathematics 610, pp. 139-251, 2014.

\bibitem{GRI}
D. Gaitsgory, N. Rozenblyum.
\newblock A study in derived algebraic geometry, vol. I.
\newblock Mathematical Surveys and Monographs, vol. 221, AMS, Providence, RI, 2017.

\bibitem{GRII}
D. Gaitsgory, N. Rozenblyum.
\newblock A study in derived algebraic geometry, vol. II.
\newblock Mathematical Surveys and Monographs, vol. 221, AMS, Providence, RI, 2017.

\bibitem{GKM}
M. Goresky, R. Kottwitz, R. McPherson.
\newblock
Codimension of root valuation strata.
\newblock Pure and Applied Mathematics Quarterly 5 , 1253-1310, 2009.

\bibitem{GKM2}
M. Goresky, R. Kottwitz, R. McPherson.
\newblock Purity of equivalued affine Springer fibers.
\newblock{\em Representation Theory}, vol. 10, pp. 130-146, (2006).

\bibitem{HR}
T. Haines, T. Richarz.
\newblock
The Test Function Conjecture for Local Models of Weil-restricted groups.
\newblock{\em Compositio Mathematica} 156, 1348-1404, 2020.

\bibitem{Hen}
B.Hennion.
\newblock
Higher dimensional formal loops.
\newblock{\em Ann. Sci. ENS} vol. 50 (4), pp. 609-663, 2017.


\bibitem{Hum}
J. Humphreys.
\newblock Conjugacy classes in semisimple algebraic groups.
\newblock Mathematical Surveys and Monographs, vol. 43, American Mathematical Soc., 2011.

\bibitem{Hum2}
J. Humphreys.
\newblock Linear algebraic groups, 2e édition.
\newblock Graduate Text in Mathematics, vol. 21, Springer-Verlag, 1981.


\bibitem{Gov}
 V. E. Govorov. 
\newblock On flat modules (Russian).
\newblock{\em Siberian Math. J.}, vol. 6 , pp. 300-304, (1965).

\bibitem{EGAIV}
A. Grothendieck, J. Dieudonné.
\newblock Éléments de géométrie algébrique. IV: Étude locale des schémas et des morphismes de schémas
Quatrième partie.
\newblock{\em Publ. Math. IHES}, vol. 32, (1967).

\bibitem{SGA1}
A. Grothendieck.
\newblock Séminaire de Géométrie Algébrique du Bois Marie: Revêtements étales et groupe
fondamental.
\newblock Lecture notes in mathematics 224, Springer-Verlag (1971).

\bibitem{He}
B. Hennion.
\newblock  
\newblock{\em Annales Sci. de l'ENS}  vol. 50 (4), pp. 609-663, 2017.


\bibitem{Ill}
L. Illusie.
\newblock Complexe cotangent et déformations I.
\newblock Lecture Notes in Mathematics, Vol. 239, Springer-Verlag, Berlin.

\bibitem{Ill2}
L. Illusie.
\newblock Complexe cotangent et déformations  II.
\newblock Lecture Notes in Mathematics, Vol. 283, Springer-Verlag, Berlin.

\bibitem{Sta}
A. J. de Jong et al.
\newblock The Stacks Project. Available at http://stacks.math.columbia.edu.

\bibitem{Kapl}
I. Kaplansky.
\newblock Projective modules. 
\newblock{\em Annals of Math.}, vol. 68, pp. 372-377, (1958).

\bibitem{KS}
M. Kashiwara, P. Schapira.
\newblock Sheaves on manifolds.
\newblock  Grundlehren der Mathematischen Wissenschaften, vol. 292, Springer-Verlag, Berlin, 1990.

\bibitem{Lau}
G. Laumon.
\newblock Correspondance de Langlands géométrique pour les corps de fonctions.
\newblock{\em Duke Math. J.} 54 (no. 2), pp. 309-360, (1987).

\bibitem{Laz}
D. Lazard.
\newblock Autour de la platitude.
\newblock{\em Bull. SMF}, vol. 97, pp. 81-128, (1969).

\bibitem{Let}
E. Letellier.
\newblock
Fourier transforms of invariant functions on finite reductive Lie algebras.
\newblock Lecture Notes in Mathematics, vol. 1859, Springer-Verlag, Berlin, 2005.

\bibitem{LZ1}
Y.Liu, W. Zheng.
\newblock Enhanced six operations and base change theorem for sheaves on Artin
stacks. preprint, arXiv:1211.5948.


\bibitem{LZ2}
Y.Liu, W. Zheng.
\newblock
Enhanced adic formalism and perverse t-structures for higher Artin stacks.
preprint, arXiv:1404.1128.

\bibitem{LuZ}
Q. Lu, W. Zheng.
\newblock
Duality and nearby cycles over general bases.
\newblock{\em Duke Math. J.}, Vol. 168, n°16, pp. 3135-3213, 2019.

\bibitem{Lu1}
J. Lurie.
\newblock
Higher Topos theory.
\newblock Annals of Mathematics Studies, vol. 170, Princeton University Press, Princeton, NJ, 2009.


\bibitem{Lu2}
J. Lurie.
\newblock
Higher algebra. 
\newblock disponible sur https://www.math.ias.edu/~lurie/papers/HA.pdf.


\bibitem{Lu3}
J. Lurie.
\newblock Spectral algebraic geometry.
\newblock  disponible sur https://www.math.ias.edu/~lurie/papers/SAG-rootfile.pdf.

\bibitem{Lus1}
G. Lusztig.
\newblock Green polynomials and singularities of conjugacy classes.
\newblock{\em Adv. in Maths. 42}, pp. 169-178, 1981.

\bibitem{Lus2}
G. Lusztig.
\newblock Intersection cohomology complexes on a reductive group.
\newblock{\em Invent. Math. 75}, pp. 205-272, 1984.

\bibitem{Mau}
D. Maulik.
\newblock Stable pairs and the HOMFLY polynomial.
\newblock{\em Invent. Math. 204}, no.2 , pp. 787-831, (2016).

\bibitem{Mig}
L. Migliorini.
\newblock HOMFLY polynomials from the Hilbert schemes of a planar curve. 
\newblock{\em Séminaire Bourbaki 71e année}, no. 1160, pp 355-389, Mars 2019, SMF Paris.

\bibitem{MV}
I. Mirkovic, K. Vilonen.
\newblock Characteristic varieties of character sheaves.
\newblock{\em Invent. Math. 93}, 405-418, (1988).


\bibitem{YunIII}
A. Oblomkov, Z. Yun.
\newblock Geometric representations of graded and rational Cherednik algebras.
\newblock{\em Advances in Math.} no. 292, pp. 601-706, (2016).

\bibitem{RG}
M. Raynaud, L. Gruson.
\newblock Critères de platitude et de projectivité.
\newblock{\em Inventiones Math. 13}, pp. 1-89, (1971).

\bibitem{Ri}
S. Riche.
\newblock Kostant section, universal centralizer, and a modular derived Satake equivalence.
\newblock{\em Math. Z.}, vol. 286, no. 1-2, 223-261, (2017).

\bibitem{Ric}
T. Richarz.
\newblock Basics on affine grassmannians, $https://timo-richarz.com/wp-content/uploads/2020/02/BoAG_02.pdf.$

\bibitem{Roz}
N. Rozenblyum.
\newblock
Filtered colimits of $\infty$-categories,
\newblock disponible sur http://www.math.harvard.edu/~gaitsgde/GL/colimits.pdf.

\bibitem{Sai}
T. Saito.
\newblock The characteristic cycle and the singular support of a constructible sheaf.
\newblock{\em Invent. Math. 207}, vol. 207, pp. 597-695, (2017).

\bibitem{Sh}
V. Shende.
\newblock Hilbert schemes of points on a locally planar curve and the
Severi strata of its versal deformation.
\newblock{\em Compositio Mathematica 148}, 2, pp. 531-547, (2012).

\bibitem{Tei}
B.Teissier.
\newblock The hunting of invariants in the geometry of discriminants.
\newblock{\em Real and complex singularities, Proc. Ninth Nordic Summer School}, pp.565-678, 1976.

\bibitem{To}
B. Toën.
\newblock Derived algebraic geometry.
\newblock{\em EMS Surv. Math. Sci. 1}, no. 2, pp. 153-245, (2014).

\bibitem{TV}
B. Toën, G. Vezzosi.
\newblock Homotopical algebraic geometry. II: Geometric stacks and
applications.
\newblock Mem. Amer. Math. Soc. 193, no. 902, (2008).


\bibitem{Zhu}
X. Zhu.
\newblock An introduction to affine Grassmannians and the geometric Satake equivalence.
\newblock IAS/Park City Mathematics Series, Vol. 24., pp. 59-154, (2017).
\end{thebibliography}
\end{document}